\documentclass[11pt]{article}

\usepackage[sort&compress,square,comma,authoryear]{natbib}
\usepackage{array}
\usepackage{graphicx}
\usepackage{amsmath}
\usepackage{amssymb}
\usepackage{hhline}
\usepackage{latexsym}
\usepackage{bm}
\usepackage{bbm}
\usepackage{cancel}
\usepackage[framed,numbered]{matlab-prettifier}

 \title{\bf Study of instability of the Fourier split-step method for 
 the massive Gross--Neveu model}

\author{ T.I. Lakoba\footnote{tlakoba@uvm.edu, \ 1 (802) 656-2610} 
 \vspace{0.5cm} \\
  Department of Mathematics and Statistics, \\
 University of Vermont, Burlington, VT 05405, USA}

\newcommand{\noi}{\noindent}
\newcommand{\und}{\underline}

\newcommand{\be}{\begin{equation}}
\newcommand{\ee}{\end{equation}}
\newcommand{\bsube}{\begin{subequations}}
\newcommand{\esube}{\end{subequations}}
\newcommand{\ba}{\begin{array}}
\newcommand{\ea}{\end{array}}
\newcommand{\To}{\rightarrow}

\newcommand{\bea}{\begin{eqnarray}}
\newcommand{\eea}{\end{eqnarray}}

\newcommand{\dt}{\Delta t}
\newcommand{\dx}{\Delta x}
\newcommand{\dtthresh}{\dt_{\rm thresh}}
\newcommand{\kmax}{k_{\max}}
\newcommand{\dk}{\Delta k}
\newcommand{\szer}{\bm{\sigma}_0}
\newcommand{\sone}{\bm{\sigma}_1}
\newcommand{\stwo}{\bm{\sigma}_2}
\newcommand{\sthr}{\bm{\sigma}_3}
\newcommand{\bP}{\bm{P}}
\newcommand{\bQ}{\bm{Q}}
\newcommand{\bma}{\bm{a}}
\newcommand{\bmb}{\bm{b}}
\newcommand{\dkp}{\delta k^{(\ge 0)}}
\newcommand{\dkm}{\delta k^{(< 0)}}
\newcommand{\wap}{\widehat{a}_{(+)}}
\newcommand{\wam}{\widehat{a}_{(-)}}
\newcommand{\wbp}{\widehat{b}_{(+)}}
\newcommand{\wbm}{\widehat{b}_{(-)}}
\newcommand{\pap}{a_{(+)}}   
\newcommand{\pam}{a_{(-)}}   
\newcommand{\pbp}{b_{(+)}}
\newcommand{\pbm}{b_{(-)}}

\newcommand{\bme}{\bm{e}}
\newcommand{\mbbi}{\mathbb{I}}

\newcommand{\mbbm}{\mathbb{M}}
\newcommand{\nss}{{\{n\}}}  

\begin{document}

\baselineskip 18 pt

\maketitle

\begin{center}
 {\bf Abstract}
\end{center}

Stability properties of the well-known Fourier split-step method used to simulate
a soliton and similar solutions of the nonlinear Dirac equations, known as the 
Gross--Neveu model, are studied numerically and analytically. 
Three distinct types of numerical instability that can occur in this case, are 
revealed and explained. 
While one of these types can be viewed as being related to the numerical instability
occurring in simulations of the nonlinear Schr\"odinger equation, 
the other two types have not been studied
or identified before, to the best of our knowledge.
These two types of instability
are {\em unconditional}, i.e. occur for arbitrarily small values of the time step. 
They also persist in the continuum limit, i.e. for arbitrarily fine spatial discretization.
Moreover, one of them persists in the limit of an infinitely large computational domain. 
It is further demonstrated that similar  instabilities also occur
for other numerical methods applied to the Gross--Neveu soliton, as well as to 
certain solitons of another relativistic field theory model, the massive Thirring.

\vskip 1.1 cm

\noi
{\bf Keywords}: \ Fourier split-step method; numerical instability; 
 solitary waves; massive Gross--Neveu model; nonlinear Dirac equations. 

\bigskip

\newpage

\section{Introduction}

The split-step, or operator-splitting, method (SSM)  is one of the
most widely used numerical tools to model evolution of linear and 
nonlinear waves. It is explicit (hence straightforward to code), 
has a number of desirable 
structure-preserving properties (e.g., preserves the $L_2$-norm exactly
in Hamiltonian systems 
and is symplectic), is easy to implement with the 2nd-order accuracy in
time \cite{68_Strang}, and allows algorithmic extensions for higher-order
accurate implementations \cite{90_Suzuki,90_Yoshida,91_Glassner,93_Bandrauk}. 
The idea of the SSM for wave equations
is to account for linear terms with spatial derivatives,
on one hand, and for all other terms, on the other hand,
 in separate substeps, where each of these substeps can be 
performed more efficiently  than a step corresponding to
the full evolution. A popular method to implement 
the substep accounting for linear terms with spatial derivatives
is via (fast) discrete Fourier transform and its inverse; hence the 
name `Fourier' (or `spectral') in the corresponding version of the SSM. 
A vast body of literature exists on just the Fourier version of the
SSM, not to mention its other (e.g., finite-difference) versions. 
In what follows 
{\em we will consider only the Fourier SSM and therefore will
omit the modifier `Fourier'}
 unless a different version of the SSM will be referred to.  
A (far from complete) list of application of this numerical method, focusing  
{\em only on nonlinear Hamiltonian} systems, includes: 
the nonlinear Schr\"odinger equation (NLS) 
\citep[]{73_HardinTappert, 84_AblowitzTaha}
\be
iu_t + u_{xx} + |u|^2 u = 0\,;
\label{e1add1_01}
\ee
Gross--Pitaevskii equation (i.e., NLS with a potential term)
with  a magnetic field term \cite[]{06_Bao};
Vlasov--Poisson equations \citep[]{76_VlasovP}; 
nonlinear Dirac equations in one spatial dimension \citep[]{89_DiracSSM}
\be
\ba{l}
 \psi_{1,\,t} + \psi_{2,\,x}  = 
   \;\; i\,( |\psi_1|^2 - |\psi_2|^2 - 1 )\psi_1, \vspace{0.1cm}\\
 \psi_{2,\,t} + \psi_{1,\,x}  =  
  -i\,( |\psi_1|^2 - |\psi_2|^2 - 1 )\psi_2 \,;
\ea
\label{e1add1_02}
\ee
generalized Zakharov equations \citep[]{03_ZakharovEqs, 04_ZakharovEqs};
and Korteweg--de Vries equation
\citep[]{73_HardinTappert}.

Since the SSM is an explicit method, it can be only conditionally stable.
Its numerical (in)stability was extensively studied for the 
NLS in one spatial dimension.
In \citep[]{86_WH}, the von Neumann analysis was applied to the SSM
simulating a solution close to the plane wave, and the 
(in)stability threshold and the instability growth rates were found.
Specifically, the instability threshold in this case is:
\be
\dtthresh = \dx^2\,/\,\pi,
\label{e1_01}
\ee
where $\dt$ and $\dx$ are discretization steps in time and space.
The quadratic dependence in \eqref{e1_01} is a consequence of the
``resonance", or ``phase matching", condition 
\bsube
\be
\omega \dt = \pi n, \qquad n\in \mathbb{N},
\label{e1_02a}
\ee
and the dispersion relation of the linear part of the NLS:
\be
\omega=k^2,
\label{e1_02b}
\ee
where $\omega$ is the frequency and
$k\in[-k_{\max},\,k_{\max})$ is the wavenumber; 
\be
k_{\max}=\pi/\dx.
\label{e1_02c}
\ee
\label{e1_02}
\esube
%
Long-term behavior of the solution
obtained by the SSM for the NLS with small and localized 
initial data was studied
in \cite{10_GaucklerLubich} by the modulated Fourier expansion
and in \cite{10_Faouetal} by the Birkhoff normal form analysis.
In \cite{12_ja}, stability of the near-soliton (i.e., localized, not small-norm)
solution was considered by a modified equation technique. 
While the instability threshold for the near-soliton
solution is still given by \eqref{e1_01}, the instability growth rate
was found to be significantly smaller than that for the plane-wave case.

Given the wide popularity of the SSM for NLS-type models, it is not 
surprising that it was also used extensively in studying solutions
of the Gross--Neveu model \eqref{e1add1_02}. In numerical simulations
of \cite{13_JCP_methods,16_SciChina,17_NMPDEs}, SSM's performance
for this model was favorably compared to that of other methods. In recent
studies \cite{17_CMS,18_ESAIM}, it was shown that the SSM is capable of
resolving distinctly different scales, i.e., of efficiently obtaining 
highly-oscillatory solutions, that occur in the non-relativistic regime
of \eqref{e1add1_02}. This method was also studied in \cite{18_preprintSchratz}
for the massive Thirring model, which is closely related to \eqref{e1add1_02}.
In Section 7 we will show that there are close similarities not only in the
analytical, but also in the numerical, solution of the two models. 
The equations of Bragg solitons in nonlinear fiber optics, 
which are an extension of the Thirring model and
thus are also mathematically related to the Gross--Neveu model, have been
also often simulated by the SSM (see, e.g., 
\cite{98_JOSAB_SSMforGS,04_Mak_SSMforGS,17_PRE_SSMforGS,17_SciRep_SSMforGS}). 
It is relevant to note that this method is also widely used in simulations of
linear Dirac--Maxwell systems; see the original paper \cite{04_DiracMaxwell}
and a recent study \cite{17_DiracMaxwell}.

Despite this significant activity, stability of the SSM for the Gross--Neveu
(or any other nonlinear spinor) model was studied analytically only in
\cite{89_DiracSSM}. However, the instability threshold obtained 
therein appears to contradict both our analysis and numerical evidence,
presented in Sections 4 and 3.1, respectively. On the other hand, in
\cite{13_JCP_methods,16_SciChina,17_NMPDEs} it was stated without a proof
that the SSM is unconditionally stable for the Gross--Neveu model. While this
may be correct for some small-norm solutions, it is {\em not} correct for 
essentially nonlinear solutions, such as the soliton.

In this work, we show analytically and verify by extensive numerical simulations
that the 2nd-order 
SSM applied to the soliton of nonlinear spinor models, including the
Gross--Neveu model \eqref{e1add1_02}, can exhibit {\em unconditional} numerical instability
(NI) via two distinctly different mechanisms. Not only does this NI persist
as $\dt\To 0$, but it is also {\em not} 
suppressed by taking
 $\dx\To 0$. Thus, no matter
how fine a discretization in both space and time one uses, the numerical solution
of \eqref{e1add1_02} by the SSM will eventually become invalid. Moreover, for one
of the two mechanisms, the NI is also not 
suppressed
 by taking an increasingly
large computational window. For a certain range of the soliton parameters, 
this NI is strong enough to destroy the solution 
even
over moderately (as opposed 
to very) long times. We will also show that this NI can occur for other nonlinear
spinor models, as well as for other numerical methods applied to such models.

It should be noted that numerical instabilities of the Gross--Neveu soliton have
previously been reported \cite{14_NLDE_numerics,15_NLDE_numerics} by methods
other than the SSM. The main differences between our work and those studies are
in that: \ (i) we reveal the mechanisms (i.e., the forms of coupling between  
different Fourier modes) via which NI can occur, and \ (ii) we {\em systematically}
study the dependence of NI's characteristics on the computational parameters,
such as $\dt$, $\dx$, and the computational domain length $L$. Further differences
will be discussed in Section 8.

Our analysis of the unconditional NI mechanisms for the Gross--Neveu solitons
will be semi-numerical. That is,
while we will be able to identify {\em how} coupling among different
Fourier modes can drive them towards NI, the corresponding equations will 
turn out to be too complicated to allow analytical treatment 
and therefore have to be solved numerically.  We believe that this
is still valuable for at least two reasons other than the conceptual reason 
of understanding NI mechanisms, which can also appear in other models and numerical
methods. First, our analysis allows us to obtain the dependence of the most
important NI parameter --- its growth rate --- on the parameters of the simulated 
solution and of the numerical scheme in a matter of seconds or at most minutes,
whereas direct numerical simulations by the SSM may take hours (when one considers
limits such as $\dt\To 0$, $\dx\To 0$, or $L\To \infty$). Second, knowing for 
what Fourier modes the NI arises, and how its growth rate depends on the problem
parameters, may help one 
to tell 
a NI from a true physical instability.
The latter does not occur for the soliton of \eqref{e1add1_02} \cite{12_BerkCom},
but does occur in models with nonlinearity higher than cubic 
(see, e.g., \cite{14_NLDE_numerics,15_NLDE_numerics}).

The main part of this paper is organized as follows. In Section 2 we will
summarize relevant facts about the soliton solution of the Gross--Neveu model.
In Section 3 we will 
present simulation results which illustrate manifestations of NI
and thus will provide motivation for the analyses in later sections.
We will begin, in Section 4, by analyzing the possibility of a {\em conditional} NI, 
not mentioned above.  For that, we will derive a counterpart of an 
instability threshold
\eqref{e1_01} for the Gross--Neveu model. However, unlike the NLS case,
here we will show that even when the time step exceeds such a ``threshold",
{\em high}-$k$ Fourier modes near the corresponding ``resonant" value $k_{\pi}$
will {\em not} grow exponentially unless the soliton is subject to a
{\em low}-$k$ instability (whether physical or numerical).

In Section 5, we will analyze the unconditional NI that occurs near
the edges of the computational spectral window. 
In Section 6 we will analyze a different kind of unconditional NI,
whereby {\em all modes} sufficiently remote from the spectral window's edges
and from $k\approx 0$ become exponentially unstable. 
In order to carry out the analyses in these two sections, 
we have to assume that the background solution is a single soliton, perturbed
infinitesimally. 
In Section 7 we will demonstrate that the 
unconditional NIs occur in a significantly wider
variety of applications than the specific problem analyzed in those sections.
Namely, we will show that these NIs can be observed when the
background solution is of a more general form than a single soliton, when 
the SSM is applied to models similar to, but other than, the Gross--Neveu,
and that they are also observed in other popular numerical methods applied
to such models. We also briefly discuss how such NIs can be suppressed. 
 In Section 8 we present conclusions of this study and discuss the relation of its
findings with those of \cite{14_NLDE_numerics,15_NLDE_numerics}. 
{\em The reader who is not interested in technical details, may read only
Sections 2, 3, and 8.}

Since our numerical results pertaining to 
unconditional NIs of the SSM may appear counter-intuitive to the reader,
in Appendix A we present a short Matlab code that can be used to reproduce 
all of the reported results for the SSM and the Gross--Neveu soliton. 
In Appendix B we discuss in some detail the issue of ``fragility" of the Gross--Neveu
soliton, introduced in Section 2. Appendices C and D contain technical results 
related to the analysis in Sections 5 and 6, respectively.


\section{Soliton solution of the Gross--Neveu model and its stability}

\begin{figure}[!ht]
\begin{minipage}{7.5cm}
\hspace*{-0.1cm} 
\includegraphics[height=5.6cm,width=7.5cm,angle=0]{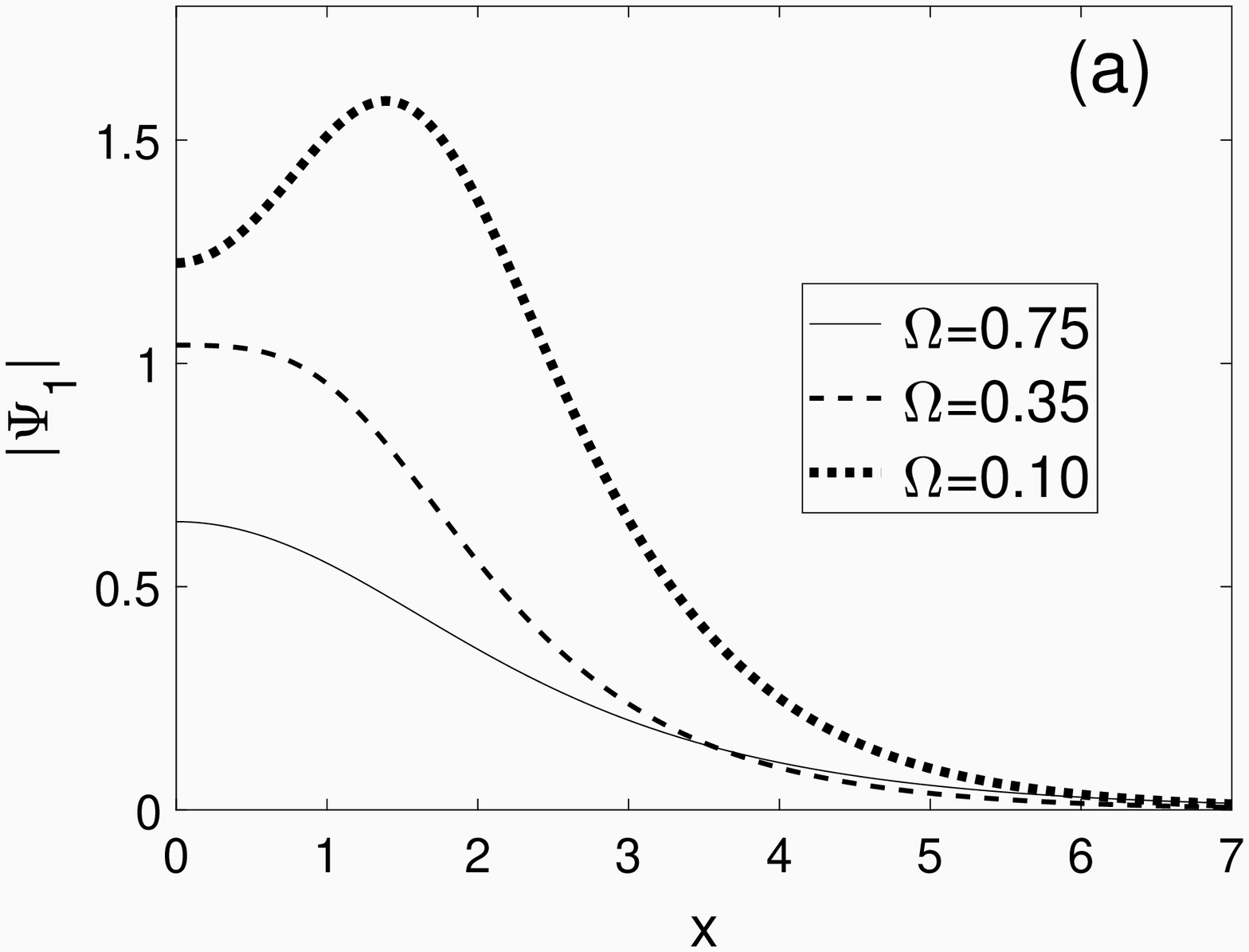}
\end{minipage}
\hspace{0.5cm}
\begin{minipage}{7.5cm}
\hspace*{-0.1cm} 
\includegraphics[height=5.6cm,width=7.5cm,angle=0]{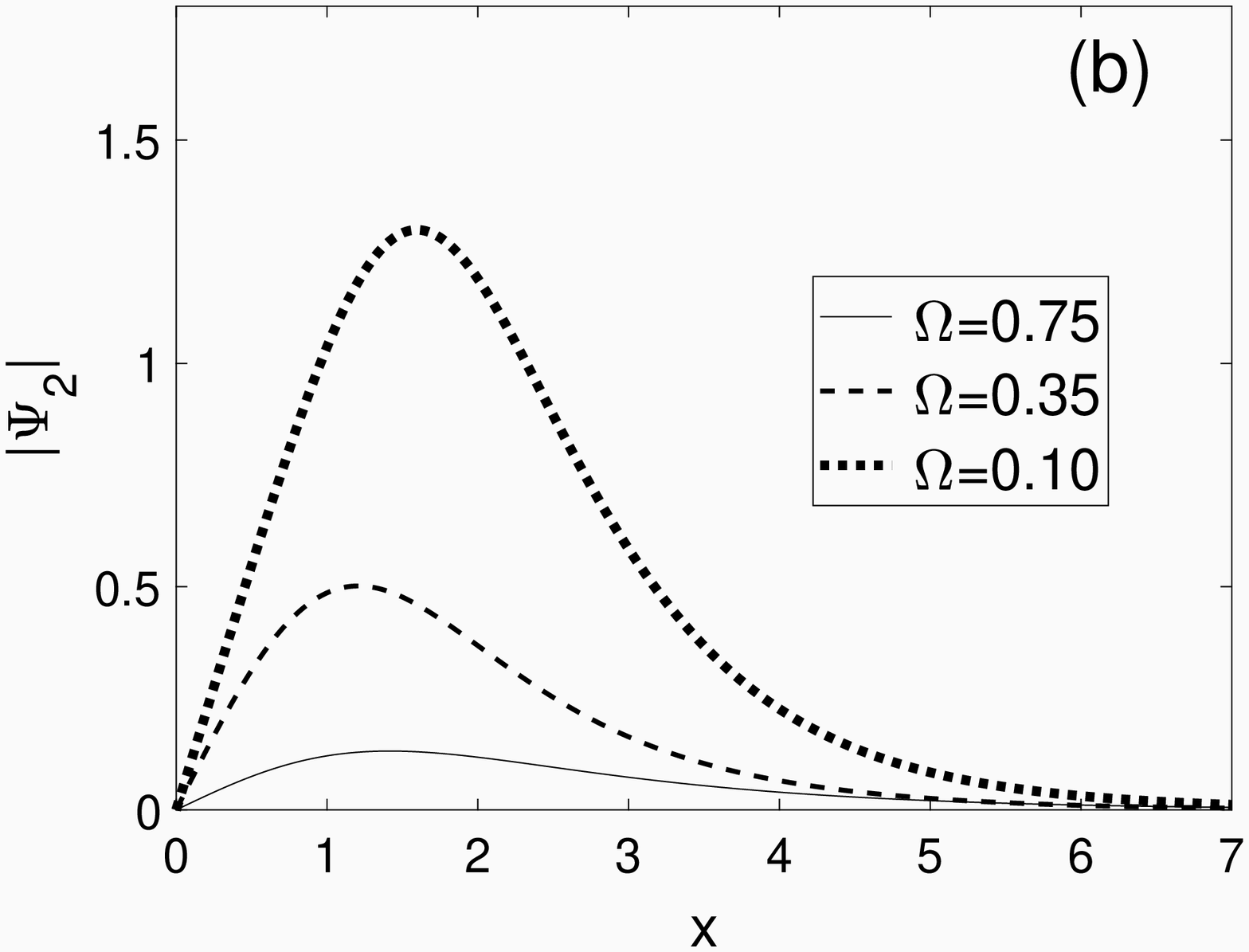}
\end{minipage}
\caption{Components $\Psi_1$ (a) and $\Psi_2$ (b) of soliton \eqref{e2_01}. 
By symmetry, $\Psi_1(-x)=\Psi_1(x)$ and $\Psi_2(-x)=-\Psi_2(x)$. 
}
\label{fig_1}
\end{figure}

The standing 
soliton solution of the Gross--Neveu model \eqref{e1add1_02} is \citep[]{75_Exact}:
\bsube
\be
\psi_1 \,=\, \Psi_1(x)\,e^{-i\Omega t}, 
\quad
\psi_2 \,=\, \Psi_2(x)\,e^{-i\Omega t},
\qquad \Omega \in (0,1);
\label{e2_01a}
\ee
\be
\Psi_1(x)\, = \, \frac{ \sqrt{2(1-\Omega)}\, \cosh(\beta x) }
				{ \cosh^2(\beta x) - \mu^2\,\sinh^2(\beta x) };
				\qquad
	\Psi_2(x) = i\mu\,\tanh(\beta x)\,\Psi_1(x);
\label{e2_01b}
\ee
\label{e2_01}
\esube
with $\beta=\sqrt{1-\Omega^2}$ and $\mu=\sqrt{(1-\Omega)/(1+\Omega)}$. 
Representative members of this family for different values of $\Omega$
are shown in Fig.~\ref{fig_1}.

A soliton moving with velocity $V\in(-1,1)$ is obtained from \eqref{e2_01} 
by a Lorentz transformation (see, e.g., \cite{81_2solColl}):
\bsube
\be
\left( \ba{c} \Psi_{1,\,\rm mov}(x,t) \\ \Psi_{2,\,\rm mov}(x,t)  \ea \right) = 
 \frac1{\sqrt{2} }
\left( \ba{cc} \sqrt{\Gamma+1} & \sqrt{\Gamma-1} \\ 
               \sqrt{\Gamma-1} & \sqrt{\Gamma+1} \ea \right) 
\left( \ba{c} \Psi_1(x_{\rm mov},t_{\rm mov}) \\ 
               \Psi_2(x_{\rm mov},t_{\rm mov})  \ea \right),
\label{e2_02a}
\ee
where
\be
\Gamma = 1/\sqrt{1-V^2}, \qquad 
x_{\rm mov} = \Gamma\,(x-x_0-Vt), \quad
t_{\rm mov} = \Gamma\,(t-V(x-x_0))\,.
\label{e2_02b}
\ee
\label{e2_02}
\esube

For future reference we present the Gross--Neveu equations linearized
on the background of the soliton \eqref{e2_01}. They are obtained by
substitution of 
\be
\psi_{1,2} = \left(\Psi_{1,2}(x) + \widetilde{\psi}_{1,2}(x,t)\right)\,
  e^{-i\Omega t}, \qquad 
	\left| \widetilde{\psi}_{1,2} \right| \ll 1
\label{e2_03}
\ee
into Eqs.~\eqref{e1add1_02} and discarding terms nonlinear in $\widetilde{\psi}_{1,2}$.
Defining a vector
$\bm{\widetilde{\psi}}=(\widetilde{\psi}_{1}, \widetilde{\psi}_{2})^T$,
one can write the linearized Eqs.~\eqref{e1_02} in the form:
\bsube
\be
\bm{\widetilde{\psi}}_t - i\Omega\, \bm{\widetilde{\psi}} + 
 \sone \bm{\widetilde{\psi}}_x \,=\, 
 i\,\bP(x) \bm{\widetilde{\psi}} + i\,\bQ(x) \bm{\widetilde{\psi}}^*
\label{e2_04a}
\ee
where Pauli matrices are:
\be
\szer = \left( \ba{cc} 1 & 0 \\ 0 & 1 \ea \right), \qquad
\sone = \left( \ba{cc} 0 & 1 \\ 1 & 0 \ea \right), \qquad
\stwo = \left( \ba{cc} 0 & -i \\ i & 0 \ea \right), \qquad
\sthr = \left( \ba{cc} 1 & 0 \\ 0 & -1 \ea \right), \qquad
\label{e2_04b}
\ee
and $\bP\equiv \sum_{j=0}^3 \bm{\sigma}_j P_j$, \ \ 
$\bQ\equiv \sum_{j=0}^3 \bm{\sigma}_j\, Q_j$ \ with:
\be
P_0 = \frac12\left(|\Psi_1|^2 + |\Psi_2|^2\right), \qquad
P_1 = 0, \qquad
P_2 = {\rm Im}\,\left( \Psi_1\Psi_2^* \right), \qquad
P_3 = \frac32\left(|\Psi_1|^2 - |\Psi_2|^2\right) \, -1\,;
\label{e2_04c}
\ee
\be
Q_0 = \frac12\left(\Psi_1^2 + \Psi_2^2\right), \qquad
Q_1 = -\Psi_1\Psi_2, \qquad
Q_2 = 0, \qquad 
Q_3 = \frac12\left(\Psi_1^2 - \Psi_2^2\right)\,.
\label{e2_04d}
\ee
\label{e2_04}
\esube
In \eqref{e2_04c} we have used the fact that Re$(\Psi_1^*\Psi_2)=0$. 
Equation \eqref{e2_04a} can be rewritten in another form which will be
convenient to refer to in what follows:
\be
\left( \ba{c} \bm{\widetilde{\psi}} \vspace{0.1cm} \\
              \bm{\widetilde{\psi}}^* \ea \right)_t \;=\; 
			{\mathcal L} \, \left( \ba{c} \bm{\widetilde{\psi}} \vspace{0.1cm} \\
              \bm{\widetilde{\psi}}^* \ea \right),
							\qquad 
	{\mathcal L} = \left( \ba{cc} i\szer\Omega - \sone\partial_x & 0  \\
	                   0 & -i\szer\Omega - \sone\partial_x \ea \right) 
			\,+\, 
			\left( \ba{cc} i{\bm P}(x) & i{\bm Q}(x) \\ -i\bQ^*(x) & -i{\bm P}^*(x) \ea \right)\,.
\label{e2_05}
\ee
Let us note that in Eqs.~\eqref{e2_04} and \eqref{e2_05} and everywhere below,
a lower-case and upper-case boldface letters indicate $(2\times 1)$ vectors and
$(2\times 2)$ matrices, respectively.

An important issue for the study of the stability of the numerical method,
which we will undertake in subsequent sections, is the {\em physical}
(as opposed to numerical) stability of the soliton \eqref{e2_01}. 
By this we mean stability of the soliton of the original partial differential
equation \eqref{e1add1_02} (i.e., regardless of any numerical scheme) to small
perturbations of its initial profile. This is often referred to as linear,
or spectral, stability.
To that end, we note that, on one hand, the soliton was semi-analytically
proved to be spectrally stable \cite[]{12_BerkCom}.
On the other hand, numerical simulations by various
methods detected a slow instability of an unknown origin 
\citep[]{14_NLDE_numerics} for $\Omega\lesssim 0.56$ 
in \eqref{e2_01}. 
In \cite{17_jaNLDE} we 
numerically demonstrated, by a non-Fourier SSM, that the Gross--Neveu soliton
is linearly stable, as predicted semi-analytically
in \cite{12_BerkCom}, for $\Omega\ge 0.01$
(and we had no evidence to suspect that it could be unstable for smaller $\Omega$).
However, we
also confirmed the observation of \citep{14_NLDE_numerics,15_NLDE_numerics} that 
for sufficiently small $\Omega$ the soliton gets increasingly (as $\Omega$ decreases)
 ``fragile" with respect to small perturbations.
   According to numerical evidence, such a perturbation may be either due to a 
	 permanently active driving source, such as the discretization error of the 
	 numerical scheme, or due to the initial condition not being the exact soliton.
	 (In the latter case, the soliton ``sheds" the excess field into dispersive
	 radiation, which affects the soliton by passing through it.)
The origin of this empirically observed
fragility is not presently understood. 
In Appendix B we present cursory evidence that it is related to
some instability of certain low-$k$ Fourier harmonics, which may occur due
to the finite size of the computational domain. 
In this work we do not further discuss this issue. 
However, two  clarifications are in order.

First, we will refer to solitons that exhibit such a fragile behavior
as {\em fragile}, as opposed to unstable. The latter term may trigger a
confusion with linear instability, which would be incorrect, given that
the soliton is linearly stable (see above). The term `fragile' is, admittedly,
not conventional in the mathematical literature, but it does accurately
describe the empirically observed soliton's behavior.
Namely, a fragile soliton may be destroyed relatively quickly by the presence
of a perturbation that is several orders of magnitude smaller than
the soliton. 

Second, it should be noted that the value $\Omega \approx 0.56$ reported
in \cite{14_NLDE_numerics} is not to be regarded as a
sharp threshold between fragile and non-fragile behaviors.
In Appendix B we present numerical evidence that weakly fragile 
behavior occurs for the soliton with $\Omega=0.75$. We believe that 
there is no sharp boundary between non-fragile and fragile behaviors, and
whether one observes signs of the latter depends on $\Omega$,
$L$, and the computational time. We therefore arbitrarily set $\Omega=0.6$
as a mark separating the two types of behavior. For $\Omega=0.6$, one
does not observe distinct signs of fragile behavior for $L\gtrsim 100$
and $t\lesssim 10,000$, and thus we refer to the corresponding solitons
as non-fragile. Solitons with $\Omega<0.6$ will then be referred to as
fragile.


\section{Numerical examples motivating subsequent analysis}

In the three subsections of this section we will provide results of numerical
simulations that demonstrate three different types of behavior of the SSM
applied to the Gross--Neveu soliton \eqref{e2_01} 
which have not been observed for the NLS soliton.

In all simulations reported below we used the 2nd-order SSM,
for which the evolution over one time step has the form:
\be
\left( \ba{c} \psi_1 \\ \psi_2 \ea \right)(x,t+\dt) = 
\exp\big[i(\dt/2){\mathcal D}\big] \; 
\exp\big[i\dt{\mathcal N}\big] \; 
\exp\big[i(\dt/2){\mathcal D}\big] \; 
\left( \ba{c} \psi_1 \\ \psi_2 \ea \right)(x,t), 
\label{e3add1_01}
\ee
where ${\mathcal D}$ is the spatial-derivative operator on the 
left-hand side (l.h.s.) of \eqref{e1add1_02} and ${\mathcal N}$
is the operator on the r.h.s. of that system. 
As stated in the Introduction, here ${\mathcal D}$ is implemented via
discrete Fourier transform \eqref{e4_02}.
The initial condition is taken
as the soliton plus a white noise of size $\sim 10^{-12}$ in the $x$-domain.


\subsection{Linear numerical instability is {\em not}
  inevitable for $\dt > ``\dt_{\rm thresh}"$}


To put the forthcoming results in context,
we note that the same logic that led one from Eqs.~\eqref{e1_02} to \eqref{e1_01},
along with the asymptotically (for large $k$) linear dispersion relation
\be
\omega = \pm k
\label{e1_03}
\ee
of \eqref{e1add1_02}, suggests that the NI threshold for the Gross--Neveu model
should be
\be
\dt_{\rm thresh}=\dx\,.
\label{e1_04}
\ee
Below we will show that this is {\em not} always the case.

\begin{figure}[!ht]
\begin{minipage}{7.5cm}
\hspace*{-0.1cm} 
\includegraphics[height=5.6cm,width=7.5cm,angle=0]{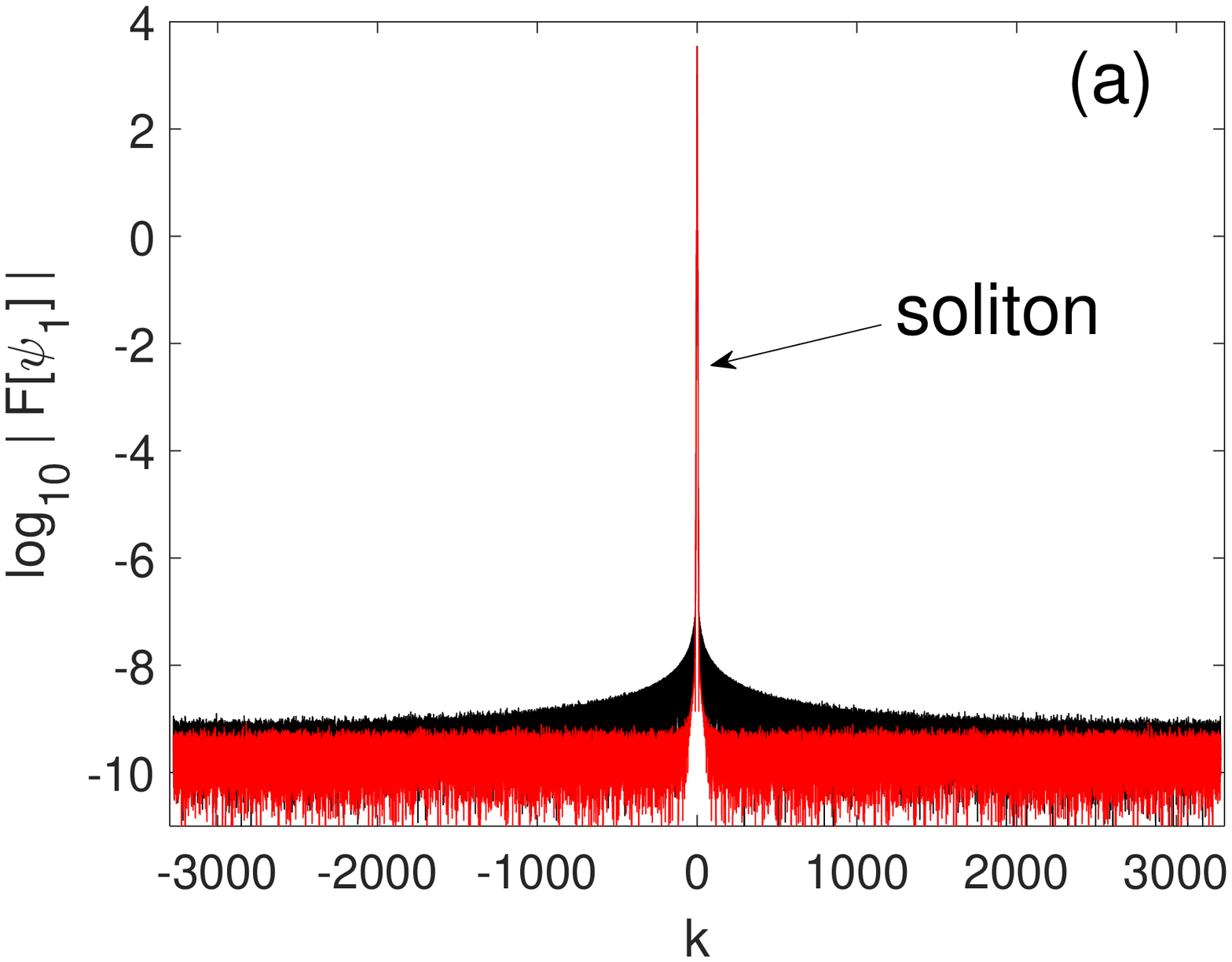}
\end{minipage}
\hspace{0.5cm}
\begin{minipage}{7.5cm}
\hspace*{-0.1cm} 
\includegraphics[height=5.6cm,width=7.5cm,angle=0]{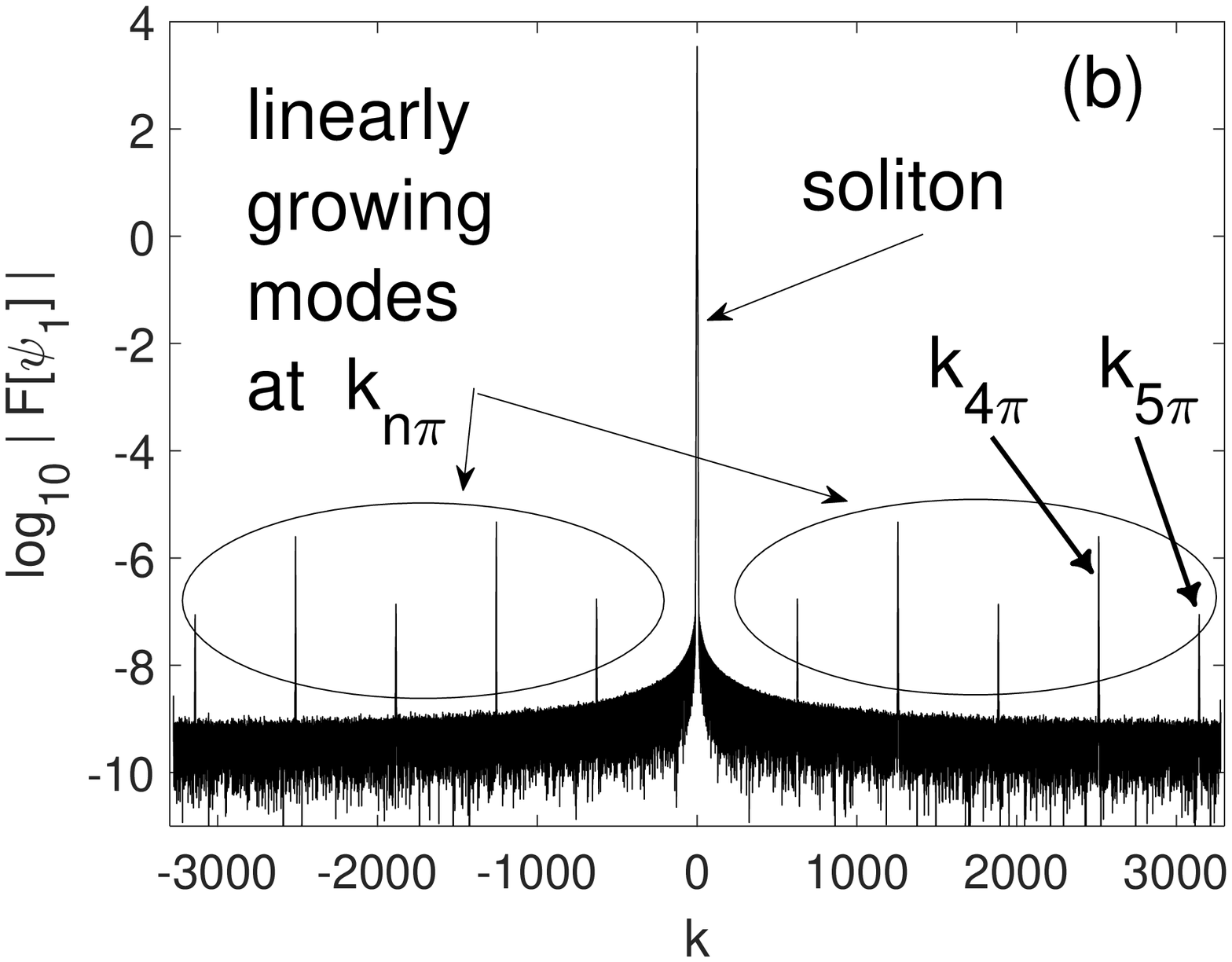}
\end{minipage}

\bigskip

\begin{minipage}{7.5cm}
\hspace*{-0.1cm} 
\includegraphics[height=5.6cm,width=7.5cm,angle=0]{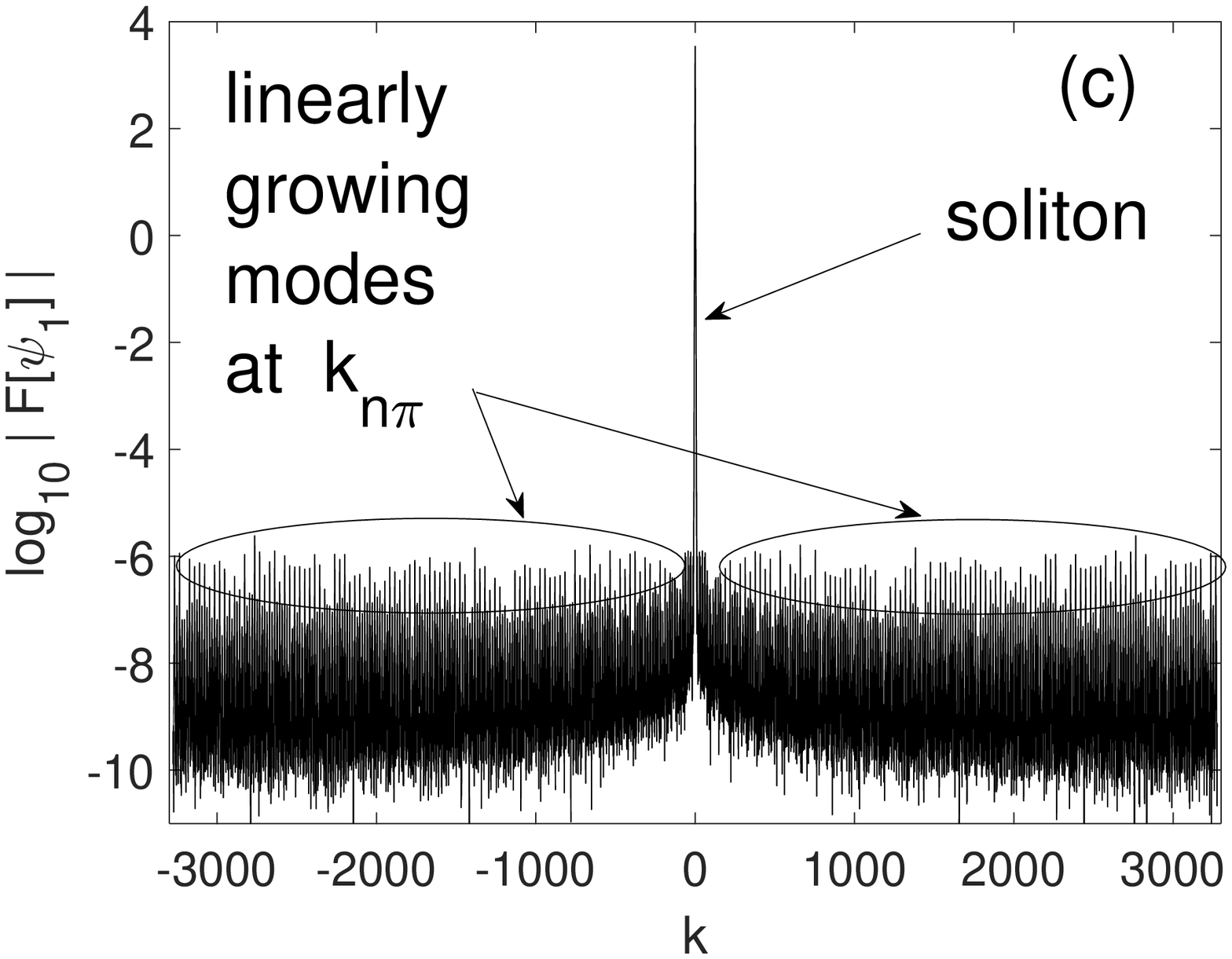}
\end{minipage}
\hspace{0.5cm}
\begin{minipage}{7.5cm}
\hspace*{-0.1cm} 
\includegraphics[height=5.6cm,width=7.5cm,angle=0]{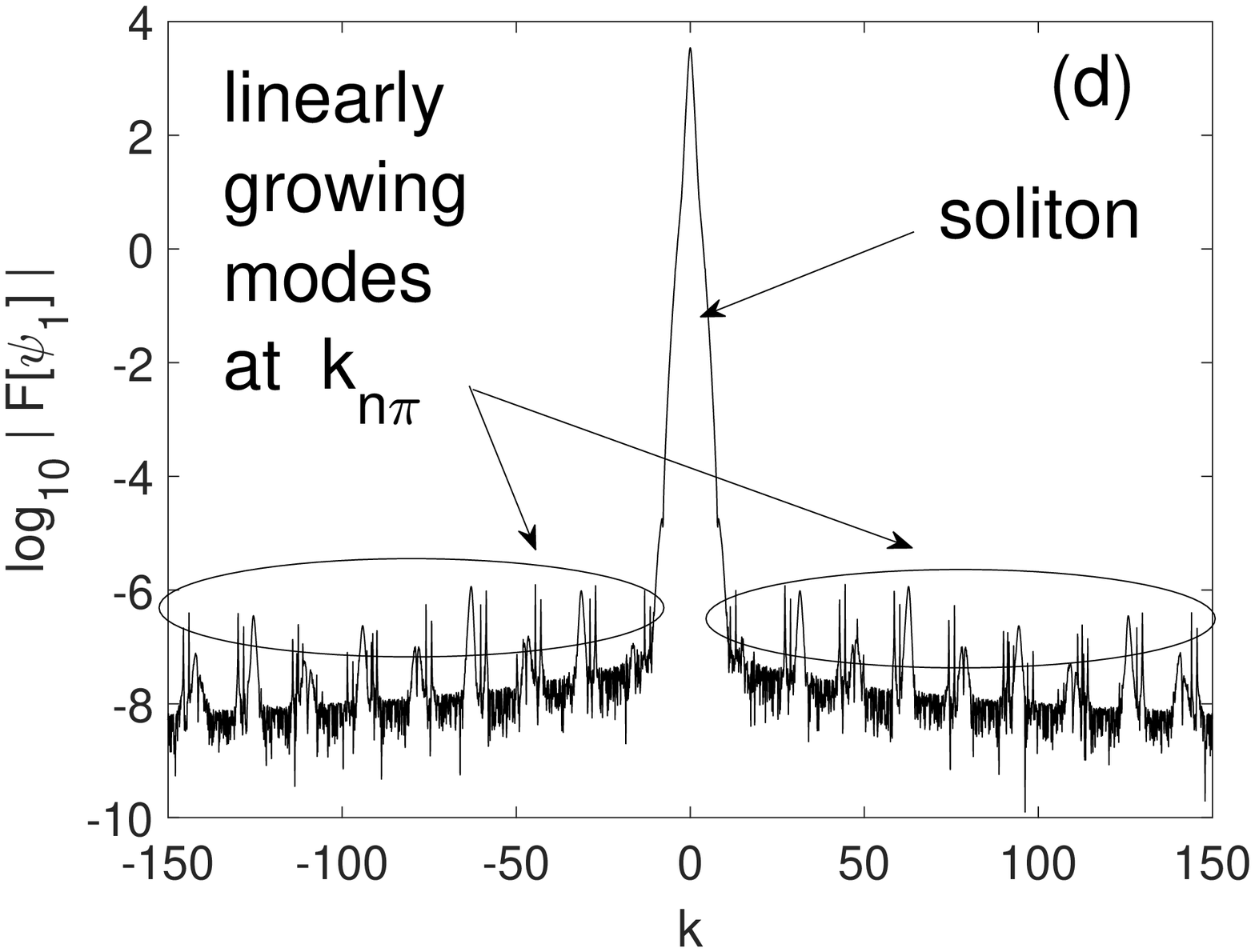}
\end{minipage}
\caption{Numerical solution's spectrum for the initial soliton
with $\Omega=0.75$ (non-fragile). 
Simulation parameters are listed in the text, 
and: \ (a) $\dt=0.9\dtthresh\equiv 0.9\dx$; \ 
(b) $\dt=0.005 \gtrsim 5\dtthresh$; \ (c) $\dt=0.2 \gtrsim 200\dtthresh$.
Panel (d) is a magnified view of (c), presented in order
 to show details of the spectrum near the soliton's spectral location. 
In (a) only, the red curve shows the spectrum of the initial condition.
In all figures in Section 3.1 we do not show results for the $\psi_2$-component
 because they are similar to those for the $\psi_1$-component. 
}
\label{fig_2}
\end{figure}

In Fig.~\ref{fig_2} we report simulation results for a non-fragile 
soliton with $\Omega=0.75$ (see Fig.~\ref{fig_1}).
Other parameters are: domain length $L=20\pi$; number of grid points
$N=2^{16}$, so that $\dx \approx 9.6\cdot 10^{-4}$; and simulation time 
$t_{\max}=10,000$. 
(An explanation for using such a small $\dx$ is found in the next paragraph.)
Figure \ref{fig_2}(a) shows that when $\dt<\dx$,
no sign of NI is observed. When we increase $\dt$ to fall in the interval
$(n\dtthresh, \, (n+1)\dtthresh)$, where $\dtthresh$ is defined by \eqref{e1_04}
and $n\in \mathbb{N}$, we observe (see Figs.~\ref{fig_2}(b)--(d))
 $2n$ spectral peaks around wavenumbers 
$\pm k_{\pi}, \ldots, \pm k_{n\pi}$, where 
\be
k_{n\pi} = n\pi/\dt;
\label{e3_01}
\ee
see Section 4 for a derivation. 
However, these peaks grow {\em not} exponentially in time, but only linearly.
We verified this fact for several different sets of $(\Omega,L,\dx,\dt)$. 
We were unable to explain the origin of these peaks and therefore will
not discuss them below. However, we emphasize that due to their linear,
as opposed to exponential, growth, they do not present a real problem
for any but extremely long-time (on the order of many million time units)
simulations.

%


The purpose of 
Fig.~\ref{fig_2}(c) is to present evidence that the stability threshold
$\dt_{\rm thresh}=O(\sqrt{\dx})$, stated in Corollary 3.5 of \cite{89_DiracSSM},
is likely incorrect. Indeed, for $\dt=0.2$ used to obtain Fig.~\ref{fig_2}(c),
$\dt/\sqrt{\dx}\approx 6.5$, but no trace of exponentially growing harmonics
is seen. (It is this ability to have the ratio $\dt/\sqrt{\dx}$ to 
significantly exceed 1 that made us choose $\dx$  to be 
as small as above. Had one not needed to ensure that $\dt/\sqrt{\dx} > 1$, 
then a much greater $\dx$ would have sufficed to make
the approximation error of the discrete Fourier transform in \eqref{e3add1_01}
negligible compared to the splitting error $O(\dt^2)$ of the 2nd-order SSM.)
While the numeric constant in front of $\dt/\sqrt{\dx}$
was not specified in \cite{89_DiracSSM}, it appears to be unlikely that
it would be greater than $6.5$, as in the case reported
in Fig.~\ref{fig_2}(c).

\begin{figure}[!ht]
\begin{minipage}{5.2cm}
\hspace*{-0.1cm} 
\includegraphics[height=4.2cm,width=5.1cm,angle=0]{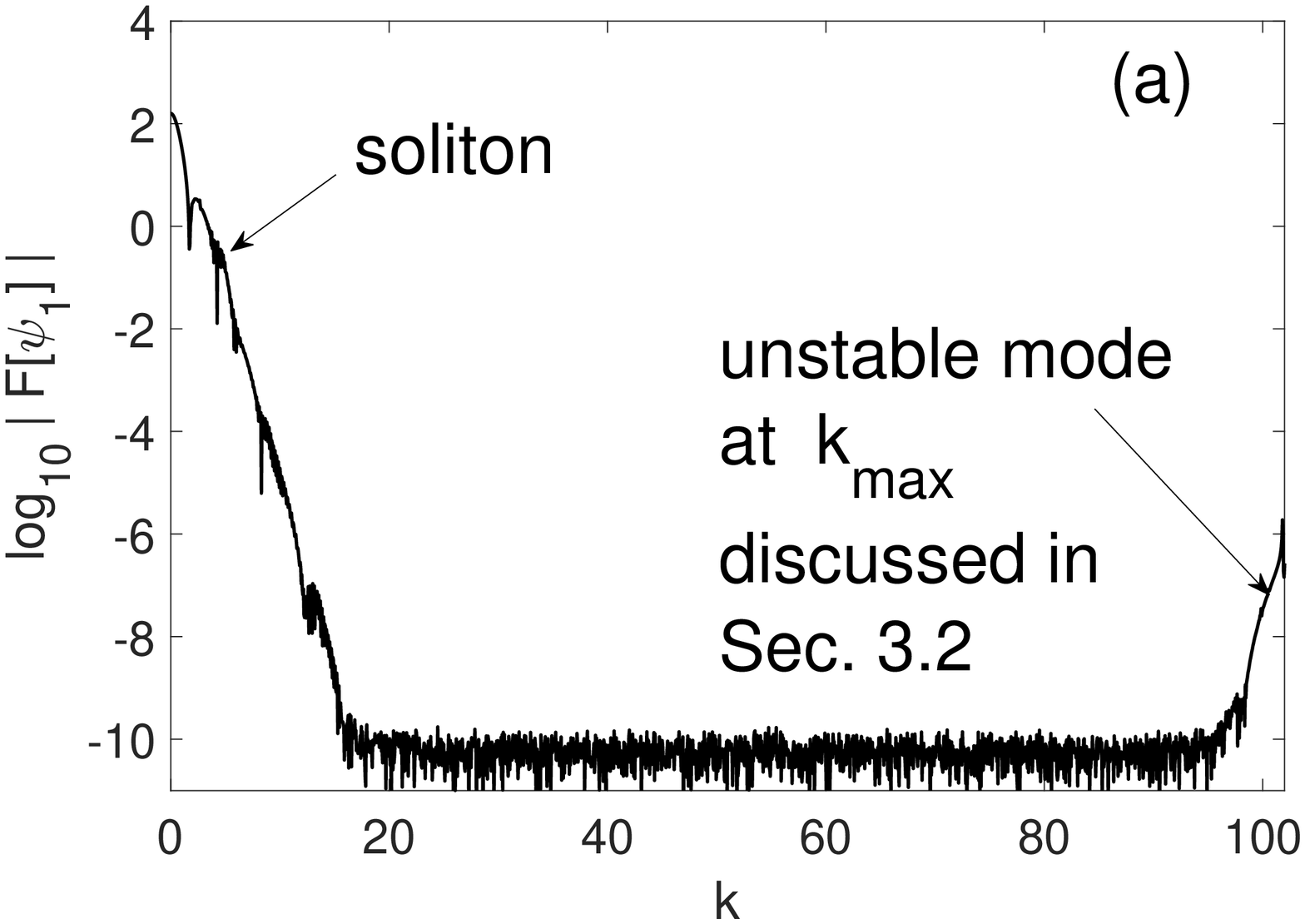}
\end{minipage}
\hspace{0.2cm}
\begin{minipage}{5.2cm}
\hspace*{-0.1cm} 
\includegraphics[height=4.2cm,width=5.1cm,angle=0]{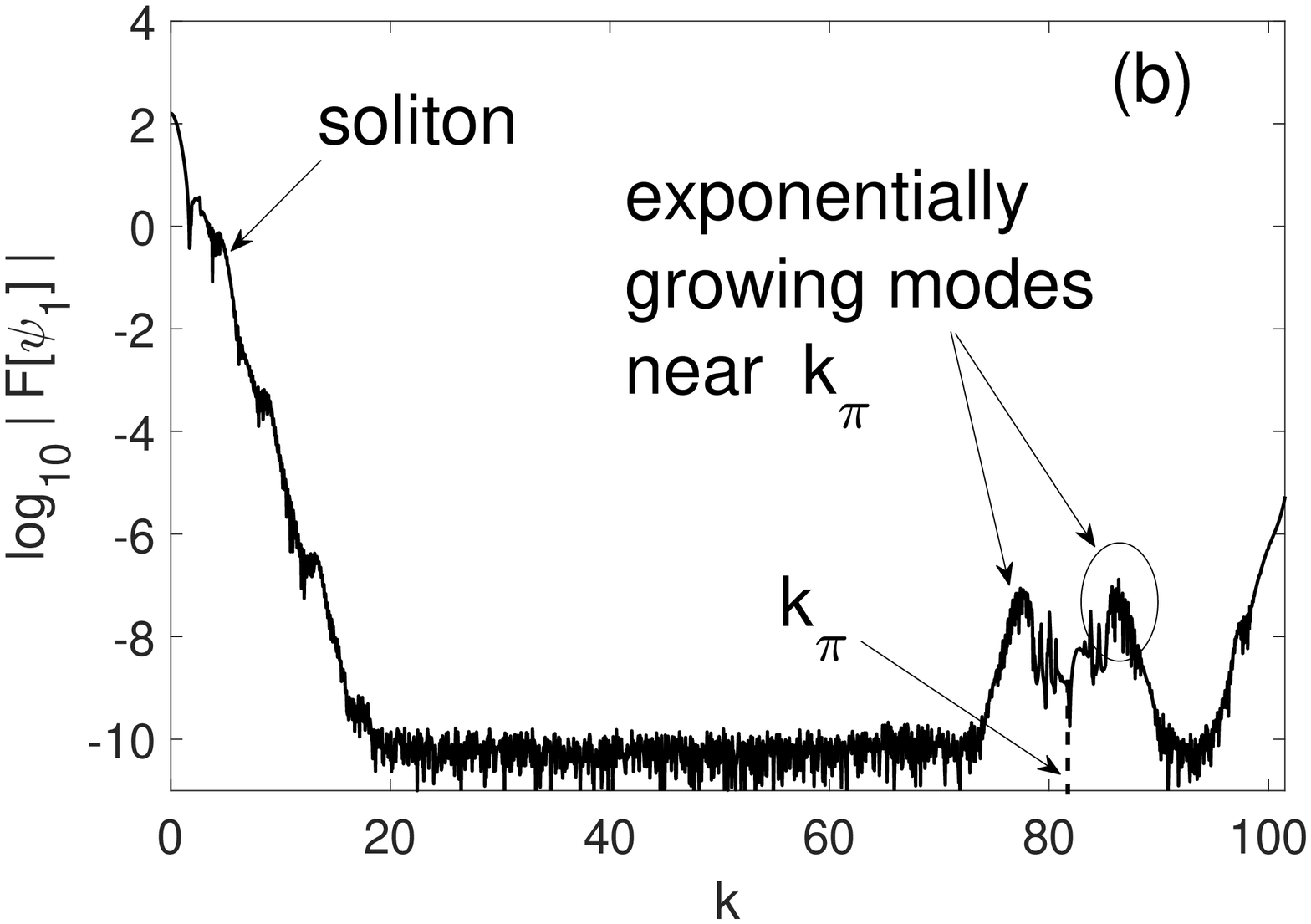}
\end{minipage}
\hspace{0.2cm}
\begin{minipage}{5.2cm}
\hspace*{-0.1cm} 
\includegraphics[height=4.2cm,width=5.1cm,angle=0]{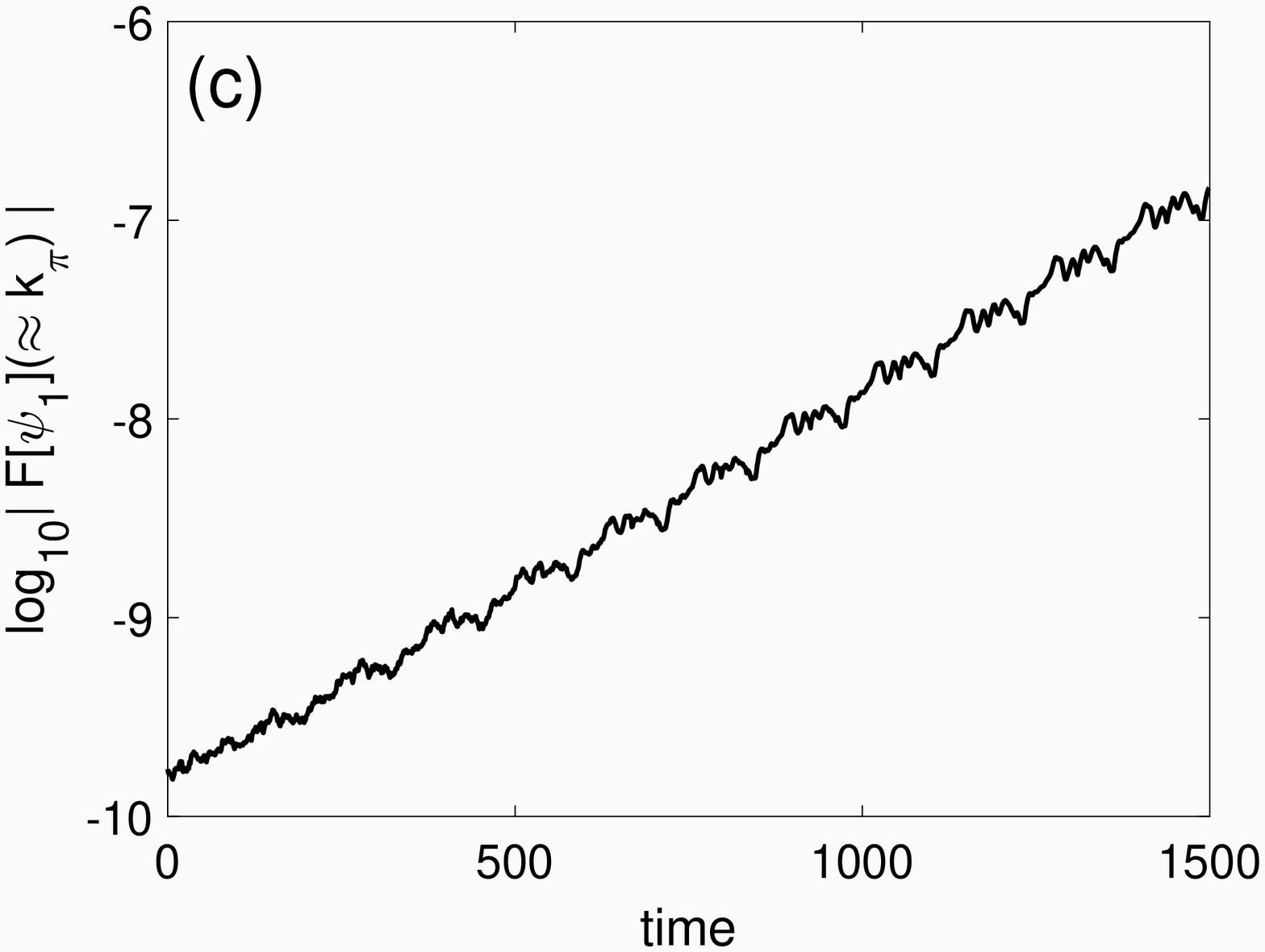}
\end{minipage}
\caption{(a,b): \ Numerical solution's spectrum for the initial soliton
with $\Omega=0.35$ (fragile). (Only the $k>0$ part of the spectrum
is shown for better visibility of details; the spectrum is 
symmetric about $k=0$.)
Simulation parameters are: $L=40\pi$;
$N=2^{12}$, so that $\dx\approx 0.031$; $t_{\max}=1,500$; 
and: \ (a) $\dt=0.9\dtthresh\equiv 0.9\dx$; \ 
(b) $\dt=1.25\dtthresh$.
\ (c) Evolution of the amplitude of the largest mode in the
group circled in panel (b).  
}
\label{fig_4}
\end{figure}

We will now show that NI near $k_{n\pi}$ {\em does} occur for $\dt>\dtthresh$ if the
background soliton is fragile. We illustrate this in 
Fig.~\ref{fig_4} for the soliton with $\Omega=0.35$ (see Fig.~\ref{fig_1}). 
In Fig.~\ref{fig_4}(a), where $\dt < \dtthresh$, the resonant wavenumber
$k_{\pi}$ is outside the computational domain, and hence no NI can exist
near it. The modes which appear at $\kmax$ (and which we verified to
grow exponentially in time) occur there for {\em any} $\dt$, no matter
how small, and hence are not related to the condition $\dt < \dtthresh$.
We will discuss them in more detail in Section 3.2. 
In Fig.~\ref{fig_4}(b), where $\dt > \dtthresh$ and hence $k_{\pi}<\kmax$, 
one observes two groups of modes on both sides of $k_{\pi}$ which grow
exponentially in time: see Fig.~\ref{fig_4}(c). We discuss their relation to
(weak) fragility of solitons in Appendix B.

%


\subsection{Unconditionally unstable modes near edges of spectral domain} 

\begin{figure}[!ht]
\begin{minipage}{7.5cm}
\hspace*{-0.1cm} 
\includegraphics[height=5.6cm,width=7.5cm,angle=0]{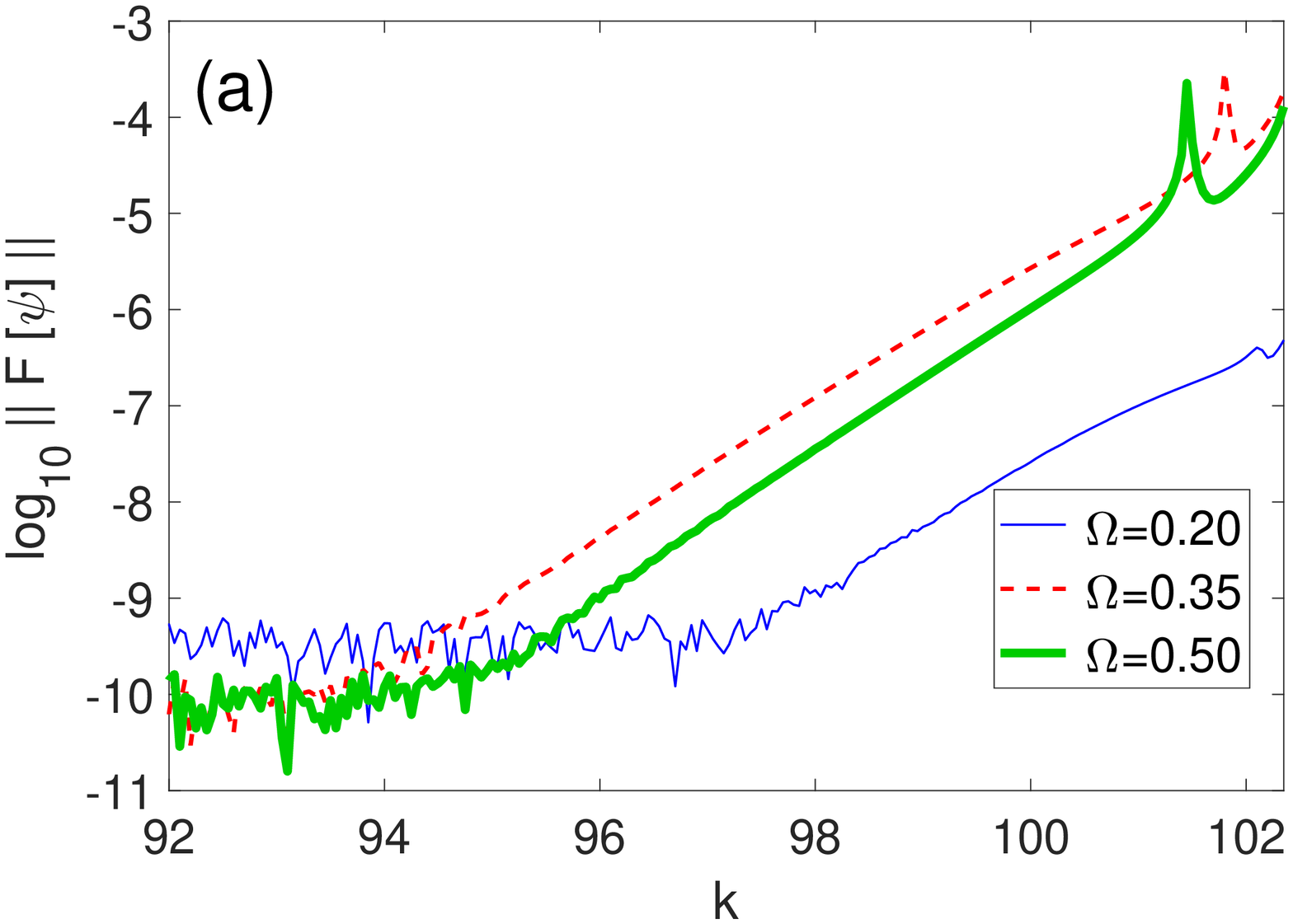}
\end{minipage}
\hspace{0.5cm}
\begin{minipage}{7.5cm}
\hspace*{-0.1cm} 
\includegraphics[height=5.6cm,width=7.5cm,angle=0]{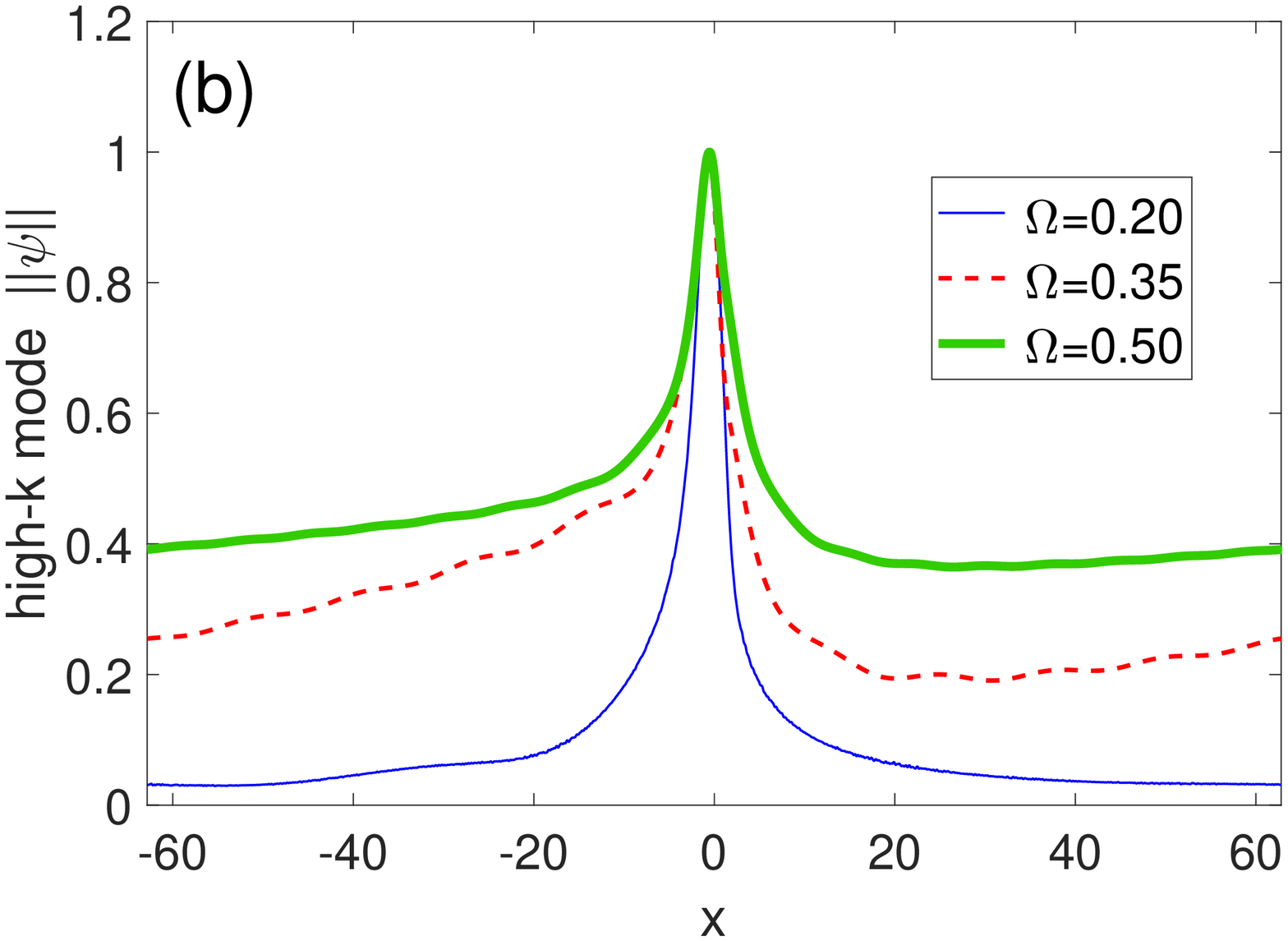}
\end{minipage}
\caption{(a) \ Close-up on the spectrum near $k_{\max}$. Note that 
simulation times were different, as follows: $t=150$ for $\Omega=0.2$;
$t=1500$ for $\Omega=0.35$, and $t=5000$ for $\Omega=0.5$.
Other simulation parameters are: $L=40\pi$, $N=2^{12}$ 
($\dx\approx 0.031$), and $\dt=0.01$. \ 
(b) \ The shape of the unstable modes in the $x$-domain, extracted
from the results reported in (a) with a high-pass filter that extends
from $k=80$ to $k=k_{\max}$. 
The modes' amplitudes are normalized to one. 
The notation $\|\cdots \|$ stands for the $\ell^2$-norm of the corresponding
two-component vector.
}
\label{fig_5}
\end{figure}

As we mentioned in the previous paragraph, the modes referred to in the title 
of this subsection can be observed when the soliton is fragile. 
We have verified that, for the parameters of Fig.~\ref{fig_4}(a),
such modes exist for $\dt$ as small as $10^{-4}$, i.e. $\dt < \dx/300$.
In the analysis in 
Section 5 we will show that, indeed, such modes persist for $\dt\To 0$.
As $\Omega$ increases, the growth rate of these modes decreases,
and vice versa. This is illustrated in Fig.~\ref{fig_5}(a); note the
different simulation times. 

Applying gentle absorbing boundary conditions, as described in 
\cite{17_jaNLDE}, was found to suppress these unstable modes for 
$\Omega=0.35$ and $0.50$, but not for $0.20$. This observation is explained
by Fig.~\ref{fig_5}(b), which shows that this mode becomes localized
(albeit with slowly decaying `tails') as $\Omega$ decreases. Clearly,
absorbing boundary conditions can have only a relatively small effect on 
such a localized mode.

Similarly to the NI for the NLS soliton \cite{12_ja}, 
relatively small variations of $L$ lead to {\em substantial} changes of 
the modes' growth rate. Moreover, they also lead to 
another unexpected behavior, which, along with the changes mentioned in
the previous sentence, is described in the next subsection.


\subsection{Unconditionally unstable ``noise floor"}

The spectrum of the numerical solution in Figs.~\ref{fig_2} and 
\ref{fig_4} 
away from the soliton 
looks like a ``floor". It appears approximately level 
because so is the spectrum of the white noise, which is added to the
initial condition in all our simulations. Thus, we will refer to this
part of the spectrum as the ``noise floor". In this subsection we
present numerical evidence that the ``noise floor" 
{\em as a whole} can also become unconditionally unstable.

\begin{figure}[!ht]
\begin{minipage}{7.5cm}
\hspace*{-0.1cm} 
\includegraphics[height=5.6cm,width=7.5cm,angle=0]{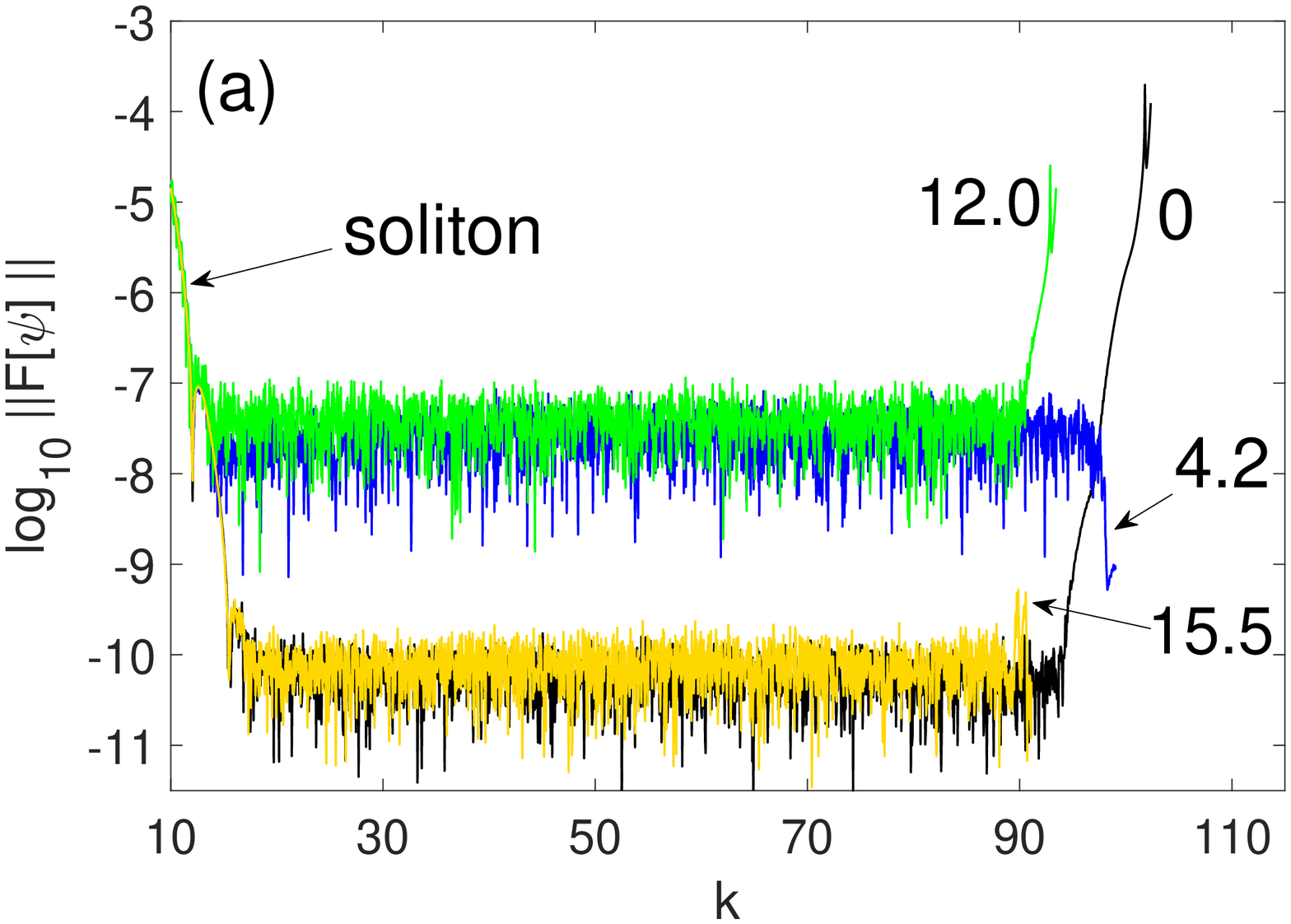}
\end{minipage}
\hspace{0.5cm}
\begin{minipage}{7.5cm}
\hspace*{-0.1cm} 
\includegraphics[height=5.6cm,width=7.5cm,angle=0]{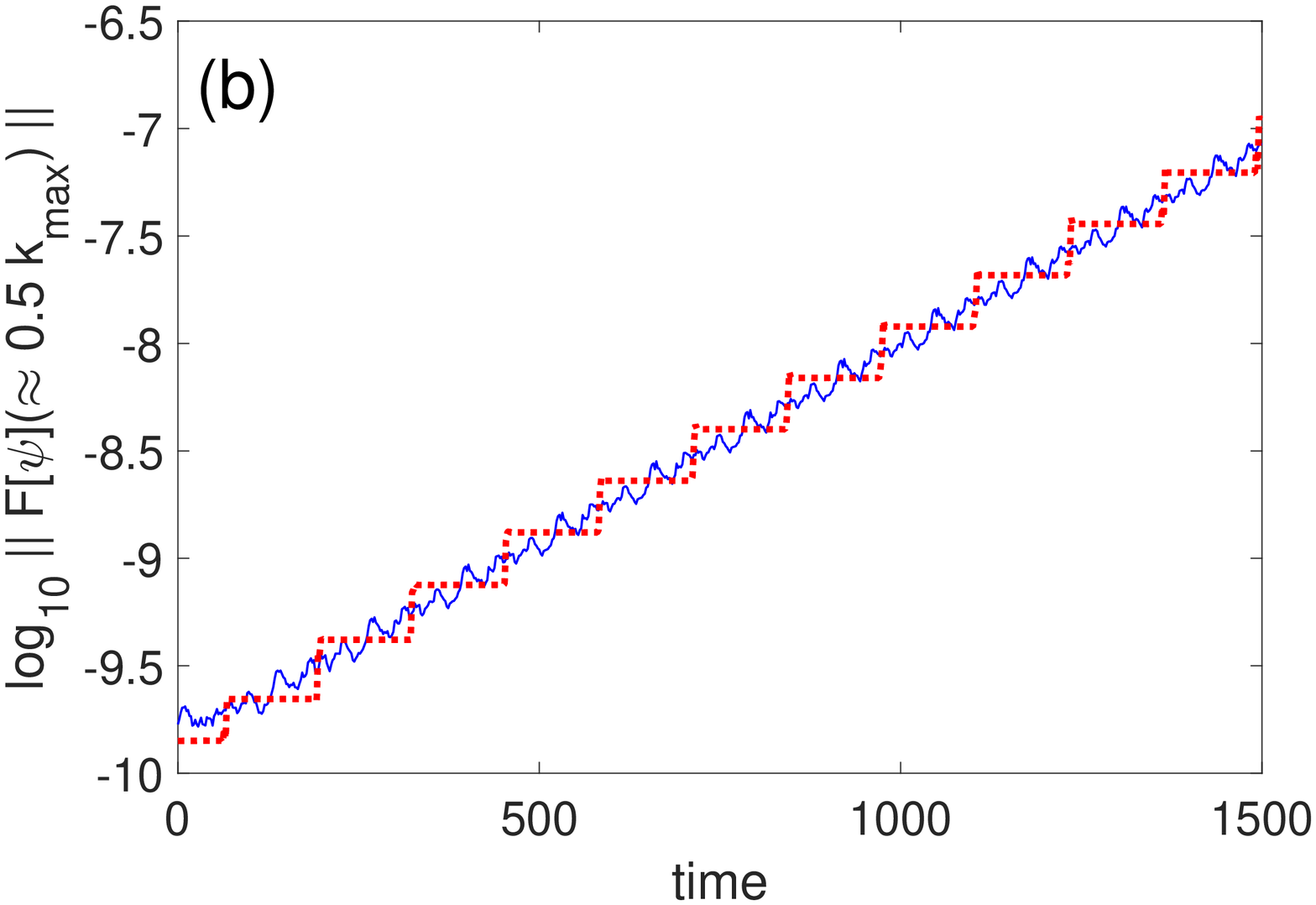}
\end{minipage}
\caption{(a) \ Part of the numerical solution's spectrum focusing on the 
``noise floor". Simulation parameters are listed in the text, and the
length of the computational domain is $L=40\pi +\delta L$, with the 
value of $\delta L$ labeling the corresponding curve. Four possible 
combinations, with the ``noise floor" and the modes near $k_{\max}$ being
(almost) stable or unstable, are shown. \ 
(b) \ Evolution of the amplitude of the Fourier modes near $k=60$ for
the curve labeled with `4.2' in panel (a). Solid and dashed lines 
correspond to a small white noise or a constant, respectively, being added
to each Fourier mode in the initial condition. 
}
\label{fig_6}
\end{figure}

In Fig.~\ref{fig_6}(a) we show part of the spectrum of the numerical 
solution for parameters that are similar to those in Fig.~\ref{fig_5}:
$\Omega=0.35$, $N=2^{12}$, $\dt=0.01$ ($<\dtthresh/3$), $t=1500$,
and four values of $L$ in the vicinity of $40\pi$. These results 
demonstrate that relatively small variation of $L$ can suppress the 
growth of unstable modes near $\pm k_{\max}$ and/or make the ``noise floor"
unstable. This NI is unconditional: reducing $\dt$ to $0.001$ left
the results shown in Fig.~\ref{fig_6} essentially unchanged. It is also
not affected by varying $N$ (for the same $L$). 

In Fig.~\ref{fig_6}(b) we demonstrate that the ``noise floor" NI is
exponential. For the purposes of analysis in Section 6 we also show the 
growth of the same Fourier mode when instead of adding to the initial 
condition a small white noise, one adds a small constant 
of the same order of magnitude to each 
Fourier mode. The analysis forthcoming  in Section 6
will explain the observed staircase 
structure of the dashed curve in Fig.~\ref{fig_6}(b).


\section{Setup of analysis, and dynamics of Fourier modes near $k_{m\pi}$}

The main results of this section are Eqs.~\eqref{e4_11} and \eqref{e4_12}. They
govern the evolution of the Fourier harmonics of a small and spectrally localized
perturbation of the soliton with wavenumbers being far from (Eqs.~\eqref{e4_11})
and near (Eqs.~\eqref{e4_12}) the resonant wavenumbers \eqref{e3_01}. Equations
\eqref{e4_12} explain the numerical observations reported in Section 3.1.
Moreover, the derivation steps carried out in this section are later used in the
analyses of Sections 5 and 6 for different types of perturbations.

Stability of the second-order SSM \eqref{e3add1_01} is the same as that of the
first-order one, 
\be
\left( \ba{c} \psi_1 \\ \psi_2 \ea \right)(x,t+\dt) = 
\exp\big[i(\dt){\mathcal D}\big] \; 
\exp\big[i\dt{\mathcal N}\big] \;  
\left( \ba{c} \psi_1 \\ \psi_2 \ea \right)(x,t), 
\label{e4_01}
\ee
because in the bulk of the calculation,
the first and last ${\mathcal D}$ $\dt/2$-substeps in \eqref{e3add1_01}
merge into one ${\mathcal D}$ $\dt$-substep in \eqref{e4_01}.
Therefore, below we study stability of the SSM \eqref{e4_01}. 
Stability of a higher-order SSM can differ from that of the first- and
second-order SSM \citep[]{12_ja}; however, this issue is outside the scope
of this study.

We will require the definition of an $N$-point 
discrete Fourier transform and its inverse:
\be
{\mathcal F}[f(x)] \equiv \widehat{f}(k_l) \, = \,
  \sum_{j=-N/2}^{N/2-1} f(x_j) \,e^{-ik_l x_j}\,;
	\qquad 
{\mathcal F}^{-1}[\widehat{f}(k)] \equiv f(x_j) \,=\, \frac1N\,
  \sum_{l=-N/2}^{N/2-1}  \widehat{f}(k_l) \,e^{ik_l x_j}\,,
\label{e4_02}
\ee
with $k_l =l\Delta k \equiv 2\pi l/L$. 
Following \cite{12_ja}, we consider stability of a numerical perturbation
whose spectral content is localized near wavenumbers $\pm k_0$ for
some $k_0\gg 1$, \ $k_0 \, \cancel{\approx}  \,\kmax$. 
Thus, the solution has the form \eqref{e2_03}, where the
vector $\bm{\widetilde{\psi}}_{\{n\}}(x)\equiv \bm{\widetilde{\psi}}(x,n\dt)$, 
defined after that equation, is sought in the form:
\be
\bm{\widetilde{\psi}}_{\{n\}}(x) = 
 \bm{\alpha}_\nss(x) e^{ik_0 x} + \bm{\beta}_\nss(x) e^{-ik_0 x}.
\label{e4_03}
\ee
Here we assume that the spectral width of $\bm{\alpha}$ and $\bm{\beta}$
is of order one and hence is much smaller than $k_0$; therefore,
the two terms in \eqref{e4_03} are well separated in the Fourier space. 
Note that here and everywhere below, the subscript in curly brackets, 
such as $\{n\}$, denotes the time level $n\dt$, whereas subscripts without
curly brackets denote either indices of Fourier harmonics, as in 
\eqref{e4_02}, or partial differentiation, as in \eqref{e1add1_01} or
\eqref{e1add1_02}, depending on the context.

Substitution of \eqref{e2_03} and \eqref{e4_03} into \eqref{e4_01},
subsequent linearization, and taking the Fourier transform result in: 
\bsube
\be
|k-k_0|=O(1): \qquad 
{\mathcal F} 
 \left[ {\bm{\alpha}}_{\{n+1\}}\,e^{ik_0 x} \right] \, e^{-i\Omega\dt} \,=\, 
e^{-i\sone k\dt} \, {\mathcal F} 
 \left[ \bm{\alpha}_\nss e^{ik_0 x} + 
   i\dt\left( \bP\bm{\alpha}_\nss + \bQ\bm{\beta}^*_\nss \right)e^{ik_0 x}\,
 \right]\,,
\label{e4_04a}
\ee
\vspace{0.1cm}
\be
|k+k_0|=O(1): \qquad 
{\mathcal F} 
 \left[ {\bm{\beta}}_{\{n+1\}}\,e^{-ik_0 x} \right]\, e^{-i\Omega\dt} \,=\, 
e^{-i\sone k\dt} \, {\mathcal F} 
 \left[ \bm{\beta}_\nss e^{-ik_0 x} + 
  i\dt\left( \bP\bm{\beta}_\nss + \bQ\bm{\alpha}^*_\nss \right)e^{-ik_0 x}\,
 \right]\,,
\label{e4_04b}
\ee
\label{e4_04}
\esube
where $\bP$ and $\bQ$ are defined in \eqref{e2_04}, and 
\be
e^{-i\sone k\dt} \equiv \szer\cos(k\dt) - i\sone\sin(k\dt).
\label{e4_05}
\ee
%
Next, in \eqref{e4_04a}, 
one writes $\exp[-i\sone k\dt]=\exp[-i\sone k_0\dt]\,\exp[-i\sone (k-k_0)\dt]$
and takes the inverse Fourier transform to obtain:
\bsube
\be
\bm{\alpha}_{\{n+1\}} \, e^{-i\Omega\dt} \,=\, e^{-i\sone k_0\dt}
\left( \bm{\alpha}_{\nss} - \dt \sone \bm{\alpha}_{\nss,\,x} + 
   i\dt(\bP\bm{\alpha}_{\nss} + \bQ\bm{\beta}_{\nss}^*)\,  \right) +O(\dt^2).
\label{e4_06a}
\ee
Indeed, since $|k-k_0|\dt \ll 1$, we have approximated
${\mathcal F}^{-1}\big[\, \exp[ - i\sone (k-k_0)\dt]\, 
 {\mathcal F}\left[ \bm{\alpha}_{\nss} \exp[ik_0x]\, \right]\;\big]$ by
$$
{\mathcal F}^{-1}\left[ (\szer -i\sone (k-k_0)\dt)\, 
 {\mathcal F}\left[ \bm{\alpha}_{\nss} e^{ik_0x}\, \right]\,\right] = 
\bm{\alpha}_\nss e^{ik_0x} - \sone 
    (\partial_x - ik_0)\left(\bm{\alpha}_{\nss} e^{ik_0x} \right) = 
\left(\bm{\alpha}_\nss - \sone \bm{\alpha}_{\nss,\,x}\right) e^{ik_0x}
$$
and $\dt\,\exp[-i\sone (k-k_0)\dt]$ by $\dt$ with accuracy $O(\dt)$. 
Similarly, and omitting $O(\dt^2)$ terms, one obtains from \eqref{e4_04b}:
\be
\bm{\beta}_{\{n+1\}} \, e^{-i\Omega\dt} \,=\, e^{i\sone k_0\dt}
\left( \bm{\beta}_{\nss} - \dt \sone \bm{\beta}_{\nss,\,x} + 
   i\dt(\bP\bm{\beta}_{\nss} + \bQ\bm{\alpha}_{\nss}^*)\,  \right).
\label{e4_06b}
\ee
\label{e4_06}
\esube

Since $k_0\gg 1$, the factors $\exp[\mp i\sone k_0\dt]$ dominate the
evolutions of $\bm{\alpha}_n$ and $\bm{\beta}_n$. Accordingly, we
use the standard perturbation theory approach and 
seek solutions of \eqref{e4_06} in the form:
\bsube
\be
\bm{\alpha}_\nss = e^{-in\sone k_0\dt}\,
 \left( \bma_\nss^{(0)} + \bma_\nss^{(1)} \right)\, ,  
\qquad 
\bm{\beta}_\nss = e^{in\sone k_0\dt}\,
 \left( \bmb_\nss^{(0)} + \bmb_\nss^{(1)} \right) \,, 
\label{e4_07a}
\ee
where $\bma_\nss^{(0)},\,\bmb_\nss^{(0)}$ vary with $n$ slowly compared to 
$\exp[\mp in\sone k_0\dt]$, and
\be
\left\{ \big|\bma_\nss^{(1)}\big|,\,\big|\bmb_\nss^{(1)}\big| \right\} = 
O\left(\dt\cdot 
  \left\{ \big|\bma_\nss^{(0)}\big|,\,\big|\bmb_\nss^{(0)}\big| \right\}
	\right) 
	\qquad \forall n\,.
\label{e4_07b}
\ee
\label{e4_07}
\esube
Substituting \eqref{e4_07} into \eqref{e4_06a}, multiplying both sides by
$\exp[i(n+1)\sone k_0\dt]$,  and again discarding $O(\dt^2)$ terms, one has: 
\bea
\bma_{\{n+1\}}^{(0)}(1-i\Omega\dt) + \bma_{\{n+1\}}^{(1)}  \; = \; 
  \bma_{\nss}^{(0)} - \dt \sone \bma^{(0)}_{\nss,\,x} + \bma_{\nss}^{(1)} \; + 
\label{e4_08} 
\vspace{0.3cm}
\\
	i\dt \left[\,
	 \left( \bP_{01} + e^{2i(n+1)\sone k_0\dt}\,\bP_{23} \right) \bma^{(0)}_\nss 
+  \left( \bQ_{01} + e^{2i(n+1)\sone k_0\dt}\,\bQ_{23} \right) (\bmb^{(0)}_\nss)^*
 \,\right].  
\nonumber 
\eea
Here $\bP_{jl}\equiv \bm{\sigma}_jP_j + \bm{\sigma}_l P_l$ and similarly for $\bQ_{jl}$,
and we have used the identity
\be
e^{i\sone s}\,\bm{\sigma}_j\, e^{-i\sone s} = e^{2i\sone s}\,\bm{\sigma}_j\,, 
\qquad j=2,3,
\label{e4add1_01}
\ee
which follows from the anti-commutation of $\sone$ with $\bm{\sigma}_{2,3}$. 
A similar equation holds for $\bmb^{(0),(1)}$.

Subsequent analysis of these equations depends on whether $k_0\approx k_{m\pi}$
for some $m\in \mathbb{N}$. To set up analyses in Sections 5 and 6, we
first consider the case where $k_0\neq k_{m\pi}$ (i.e.,
$|k_0-k_{m\pi}|\gg 1$). Then terms in \eqref{e4_08}
which are proportional to $\bma^{(0)}$ and $(\bmb^{(0)})^*$ split into 
two groups: fast-oscillating (proportional to $\exp[2i(n+1)\sone k_0\dt]$)
and those varying slowly with $n$. If the latter terms do not all cancel out,
they will drive the $\bma^{(1)}_\nss$-term, which therefore will grow, and
condition \eqref{e4_07b} will eventually be violated. To prevent this
from occurring, one requires that all these slowly varying terms cancel out,
yielding: 
\bsube
\be
\bma_{\{n+1\}}^{(0)}(1-i\Omega\dt) \,=\,
(\szer - \dt \sone\partial_x) \bma_{\nss}^{(0)}  + 
  i\dt  \left( \bP_{01} \bma^{(0)}_\nss 
+  \bQ_{01} (\bmb^{(0)}_\nss)^* \right).
\label{e4_09a}
\ee
The remaining terms provide an equation for $\bma^{(1)}$:
\be
\bma_{\{n+1\}}^{(1)} - \bma_\nss^{(1)} = \dt\,e^{2i(n+1)\sone k_0\dt} \bm{c}_\nss,
\qquad 
\bm{c}_\nss \equiv \bP_{23}\bma_{\nss}^{(0)} +  \bQ_{23}(\bmb^{(0)}_\nss)^*;
\label{e4_09b}
\ee
\label{e4_09}
\esube
note that $\bm{c}_\nss$ varies slowly with $n$. 
Due to this fact and the presence of a fast-varying exponential on the
r.h.s.~of \eqref{e4_09b}, the solution of that equation does not grow with $n$,
 and hence condition \eqref{e4_07b} holds. 
Thus, since the terms $\bma^{(1)},\,\bmb^{(1)}$ remain small 
compared to $\bma^{(0)},\,\bmb^{(0)}$ at all times, we will no 
longer consider the former terms and will focus on the latter.
Taking into account
the slow dependence of $\bma^{(0)}_\nss$ on $n$, we approximate 
\be
\bma^{(0)}_{\{n+1\}} = \bma^{(0)}_\nss + \dt (\partial_t\bma^{(0)})_\nss + O(\dt^2).
\label{e4_10}
\ee
Finally, substituting \eqref{e4_10} into \eqref{e4_09a} and 
omitting $(\dt^2)$ terms, we obtain:
\bsube
\be
\bma_t + \sone \bma_x - i\Omega \bma = i\bP_{01}\bma + i\bQ_{01}\bmb^*.
\label{e4_11a}
\ee
Here and below we will omit the superscipt $(0)$ of $\bma$ and $\bmb$.
Similar calculations yield the equation for $\bmb$:
\be
\bmb_t + \sone \bmb_x - i\Omega \bmb = i\bQ_{01}\bma^* + i\bP_{01}\bmb.
\label{e4_11b}
\ee
\label{e4_11}
\esube
We will return to the derivation in this paragraph in Sections 5 and 6.

Now, in the case $k_0=k_{m\pi}$, since $2k_0\dt=2m\pi$,
all terms in \eqref{e4_08} vary slowly with $n$. 
(This also holds for $|k_0-k_{m\pi}|=O(1)$, since we can use the
freedom in \eqref{e4_03}, which requires only that $\bm{\alpha}$
and $\bm{\beta}$ vary in $x$ on the scale of order one.)
Following the derivation in the previous paragraph, one obtains,
instead of \eqref{e4_11}, equations:
\bsube
\be
\bma_t + \sone \bma_x - i\Omega \bma = i\bP\bma + i\bQ\bmb^*.
\label{e4_12a}
\ee
\be
\bmb_t + \sone \bmb_x - i\Omega \bmb = i\bQ\bma^* + i\bP\bmb.
\label{e4_12b}
\ee
\label{e4_12}
\esube
Thus, vector $(\,\bma^T,\,(\bmb^*)^T\,)^T$, representing a 
perturbation that is spectrally localized near $k_{m\pi}$,
 satisfies the same equation \eqref{e2_05} that governs the evolution
of a low-$k$ perturbation to the soliton.
Therefore, modes near $k_{m\pi}$ have the same stability properties as the 
low-wavenumber modes. In other words, NI near $k_{m\pi}$ occurs if and
only if low-$k$ perturbations of the soliton grow exponentially. Thus, the above analysis
has explained the numerical observations of Section 3.1 and Appendix B.


\section{Unconditional numerical instability near $|k|=\kmax$} 

In this section we will present a theory explaining the numerical results
of Section 3.2. This will be accomplished by numerically solving the eigenvalue 
problem \eqref{e5_07}.

A numerical perturbation whose spectral content is concentrated near the
edges of the computational domain is sought in the form \eqref{e4_03},
where now $k_0=\kmax$ and, in addition, the Fourier transform of 
$\bm{\alpha}(x)$ (of $\bm{\beta}(x)$) contains harmonics with 
only negative (respectively, nonnegative) wavenumbers:
\be
\bm{\alpha}(x) = \sum_{l=1}^M \widehat{\bm{\alpha}}_l \,
    e^{-ik_l x}, 
		\qquad 
\bm{\beta}(x) = e^{-i\dk x}\,\sum_{l=1}^M \widehat{\bm{\beta}}_l \,
    e^{ik_l x},
\label{e5_01}
\ee
where $1\ll M \ll N/2$. The last strong inequality holds because the
spectral content of the perturbation is concentrated near the edges
of the spectral domain.  
Recall that the computational spectral window is 
$k_l \in [-\kmax, \, \kmax-\Delta k]$, where $\kmax$ is defined in
\eqref{e1_02c} and $k_l$, $\Delta k$ are defined after \eqref{e4_02}.

When we substitute Eqs.~\eqref{e4_03} with \eqref{e5_01}
 into \eqref{e4_01}, then, similarly to 
the r.h.s.~of \eqref{e4_04}, we obtain terms like 
$\bP\bm{\alpha}\exp[i\kmax x]$ etc.. Note, however, that the coefficient
$\bP(x)\bm{\alpha}(x)$ had Fourier harmonics with wavenumbers 
of both signs, due to
$\bP(x)$ having such harmonics. Therefore, such terms are to be split
into two groups:
\bsube
\be
\bP\bm{\alpha} \,\equiv\, \left[ \bP\bm{\alpha} \right]^{(<0)} + 
    \left[ \bP\bm{\alpha} \right]^{(\ge 0)},
\label{e5_02a}
\ee
where the superscript indicates what Fourier harmonics the term has. 
When multiplied by $\exp[i\kmax x]$, the former group of terms 
on the r.h.s.~of \eqref{e5_02a} will have the spectral content near
the right edge of the spectral domain, while the latter group's 
spectral content will ``spill over" to the left edge due to aliasing. 
Thus,
\be
\bP\bm{\alpha}\,e^{i\kmax x}  \,=\, 
 \left[ \bP\bm{\alpha} \right]^{(<0)}\,e^{i\kmax x} + 
    \left[ \bP\bm{\alpha} \right]^{(\ge 0)}\,e^{-i\kmax x}.
\label{e5_02b}
\ee
\label{e5_02}
\esube
With this observation in mind, the counterparts of Eqs.~\eqref{e4_04}
near the edges of the spectral domain become:
\bsube
\be
 \widehat{\bm{\alpha}}_{\{n+1\}} \, e^{i\sone (\kmax+\dkm)\dt-i\Omega\dt} 
 \,=\, 
 \widehat{\bm{\alpha}}_\nss + 
   i\dt\, {\mathcal F}\left[ 
	       \left[ \bP\left( \bm{\alpha}_\nss  + \bm{\beta}_\nss \right) + 
	        \bQ\left( \bm{\alpha}^*_\nss + \bm{\beta}^*_\nss \right) \right]^{(<0)} 
		\, \right] \,,
\label{e5_03a}
\ee
\vspace{0.1cm}
\be
\widehat{\bm{\beta}}_{\{n+1\}}\, e^{i\sone (-\kmax+\dkp)\dt-i\Omega\dt} 
 \, = \,  
 \widehat{\bm{\beta}}_\nss + 
  i\dt \,{\mathcal F} 
  \left[ \left[ \bP\left( \bm{\alpha}_\nss  + \bm{\beta}_\nss \right) + 
	        \bQ\left( \bm{\alpha}^*_\nss + \bm{\beta}^*_\nss \right) \right]^{(\ge 0)}
 \,\right]\,, 
\label{e5_03b}
\ee
\label{e5_03}
\esube
where
$$
\dkm = k-\kmax < 0, \qquad \dkp = k+\kmax \ge 0.
$$
Note that all Fourier transforms in \eqref{e5_03a} (in \eqref{e5_03b})
are evaluated at $\dkm$ (respectively, at $\dkp$). 

In Section 4 we were able to take the inverse Fourier transform of 
Eqs.~\eqref{e4_04}, which are counterparts of Eqs.~\eqref{e5_03}, and
proceed with the analysis in the $x$-space. In the case of Eqs.~\eqref{e5_03},
this would result not in differential equations, such as \eqref{e4_11}
or \eqref{e4_12}, but in integro-differential ones, due to the 
separation of positive and negative wavenumbers on the r.h.s.~of
\eqref{e5_03}. Quantitative analysis of such integro-differential 
equations would be more
difficult than analysis of the original Fourier-space 
equations \eqref{e5_03}. Therefore, below we will proceed with the latter
analysis.

Following the derivation of Eqs.~\eqref{e4_11}, except for not taking 
the inverse Fourier transform, one obtains from \eqref{e5_03}:
\bsube
\be
\widehat{\bma}_t \,=\, i \left( \szer \Omega - \sone \dkm \right) \widehat{\bma}
 + i\, {\mathcal F} \left[ 
	       \left[ \bP_{01} \bma  + \bP_{23} \bmb + 
	        \bQ_{23} \bma^* + \bQ_{01} \bmb^* 
					 \right]^{(<0)} 		\, \right] \,, 
 \label{e5_04a}
\ee
\be
\widehat{\bmb}_t \,=\, i \left( \szer \Omega - \sone \dkp \right) \widehat{\bmb}
 + i\, {\mathcal F} \left[ 
	       \left[ \bP_{23} \bma  + \bP_{01} \bmb + 
	        \bQ_{01} \bma^* + \bQ_{23} \bmb^* 
					 \right]^{(\ge 0)} 		\, \right] \,.
 \label{e5_04b}
\ee
\label{e5_04}
\esube
Using the same justification as in \eqref{e4_11}, 
we have neglected small terms $\widehat{\bma}^{(1)}$, 
$\widehat{\bmb}^{(1)}$ and omitted the superscript $(0)$. 
System \eqref{e5_04} determines stability of a $8M$-dimensional vector:
$M$ harmonics in each of the two-component vectors $\widehat{\bma}$ 
and $\widehat{\bmb}$ are coupled with as many harmonics of their
complex conjugates. It is possible to halve the size of the involved
vectors (i.e., from $8M$ to $4M$) by means of the following substitution:
\be
\widehat{\bma} = \wap \bme_{(+)} + \wam \bme_{(-)}, \qquad
\widehat{\bmb} = \wbp \bme_{(+)} + \wbm \bme_{(-)}, \qquad
\bme_{(\pm)} \equiv \left( \ba{r} 1 \\ \pm 1 \ea \right);
\label{e5_05}
\ee
the same decomposition also holds for $\bma$ and $\bmb$.
Vectors $\bme_{(\pm)}$ satisfy the following relations:
\be
\szer \bme_{(\pm)} = \bme_{(\pm)}, \quad 
\sone \bme_{(\pm)} = \pm\bme_{(\pm)}, \quad
\stwo \bme_{(\pm)} = \mp i \bme_{(\mp)}, \quad
\sthr \bme_{(\pm)} = \bme_{(\mp)}.
\label{e5_06}
\ee
The resulting equations, which couple harmonics 
$\wap$, $\wbm$, $\wam^{\,*}$, $\wbp^{\,*}$,
are found in Appendix C. Here we
present only their matrix form, which is needed for a
stability analysis. 
Defining 
$M$-dimensional column vectors, e.g.: 
$\und{\widehat{a}}_{(+)} = \big[ (\wap)_1, \ldots \,, (\wap)_M \big]^T$,
etc., where the numeric subscript denotes the harmonic's number (see
\eqref{e5_01}) and `$T$' denotes the transpose, one can write 
system \eqref{B_01} as:
\bsube
\be
 \und{\mathbbm{s}}_{\;t}
	\;=\; 
	i\, \left( \mathbbm{D} + \mathbbm{C}
	\right)\, \und{\mathbbm{s}} \,, 
\label{e5_07a}
\ee
where:
\be
\mathbbm{D}  = 
	\left( \ba{cccc} 
	       \dk\, \mbbm + \Omega\, \mbbi & 0 & 0 & 0 \vspace{0.1cm} \\
				 0 & \dk (\mbbm - \mbbi) + \Omega\,\mbbi & 0 & 0 \vspace{0.1cm} \\
				 0 & 0 &  \dk\,\mbbm -\Omega\,\mbbi & 0 \vspace{0.1cm} \\	
				 0 & 0 & 0 & \dk (\mbbm - \mbbi) - \Omega\,\mbbi  
				 \ea  \right) \,,
\label{e5_07b}
\ee
\be
   \und{\mathbbm{s}} \,=\, 
   \left( \ba{c} \und{\widehat{a}}_{(+)} \vspace{0.1cm} \\
   \und{\widehat{b}}_{(-)}  \vspace{0.1cm} \\  
	 \und{\widehat{a}}_{(-)}^* \vspace{0.1cm} \\
	 \und{\widehat{b}}_{(+)}^*   \ea  \right), 
	 \qquad
\mathbbm{C}  = 
  \frac1N\,
	\left( \ba{cccc} 
	       \mathbb{C}_{11} & \mathbb{C}_{12} & 
				 \mathbb{C}_{13} & \mathbb{C}_{14} \vspace{0.1cm} \\
				 \mathbb{C}_{21} & \mathbb{C}_{22} & 
				 \mathbb{C}_{23} & \mathbb{C}_{24} \vspace{0.1cm} \\
				 \mathbb{C}_{31} & \mathbb{C}_{32} & 
				 \mathbb{C}_{33} & \mathbb{C}_{34} \vspace{0.1cm} \\
				 \mathbb{C}_{41} & \mathbb{C}_{42} & 
				 \mathbb{C}_{43} & \mathbb{C}_{44} 
				 \ea  \right) \,,
\label{e5_07c}
\ee
\label{e5_07}
\esube
$N$ is the total number of grid points in the computational domain
(see \eqref{e4_02});
$\mbbi$ is the $M\times M$ identity matrix; $\mbbm={\rm diag}\,[1:M]$ 
is the $M\times M$ diagonal matrix with integer entries; and 
the $M\times M$ blocks $\mathbb{C}_{jm}$, which appear from the
convolution-like terms on the r.h.s.~of \eqref{e5_04}, are
written out explicitly in Appendix C. Each of these blocks has
a logically clear structure and can be easily programmed.

\begin{figure}[!ht]
\begin{minipage}{5.2cm}
\hspace*{-0.1cm} 
\includegraphics[height=4.2cm,width=5.1cm,angle=0]{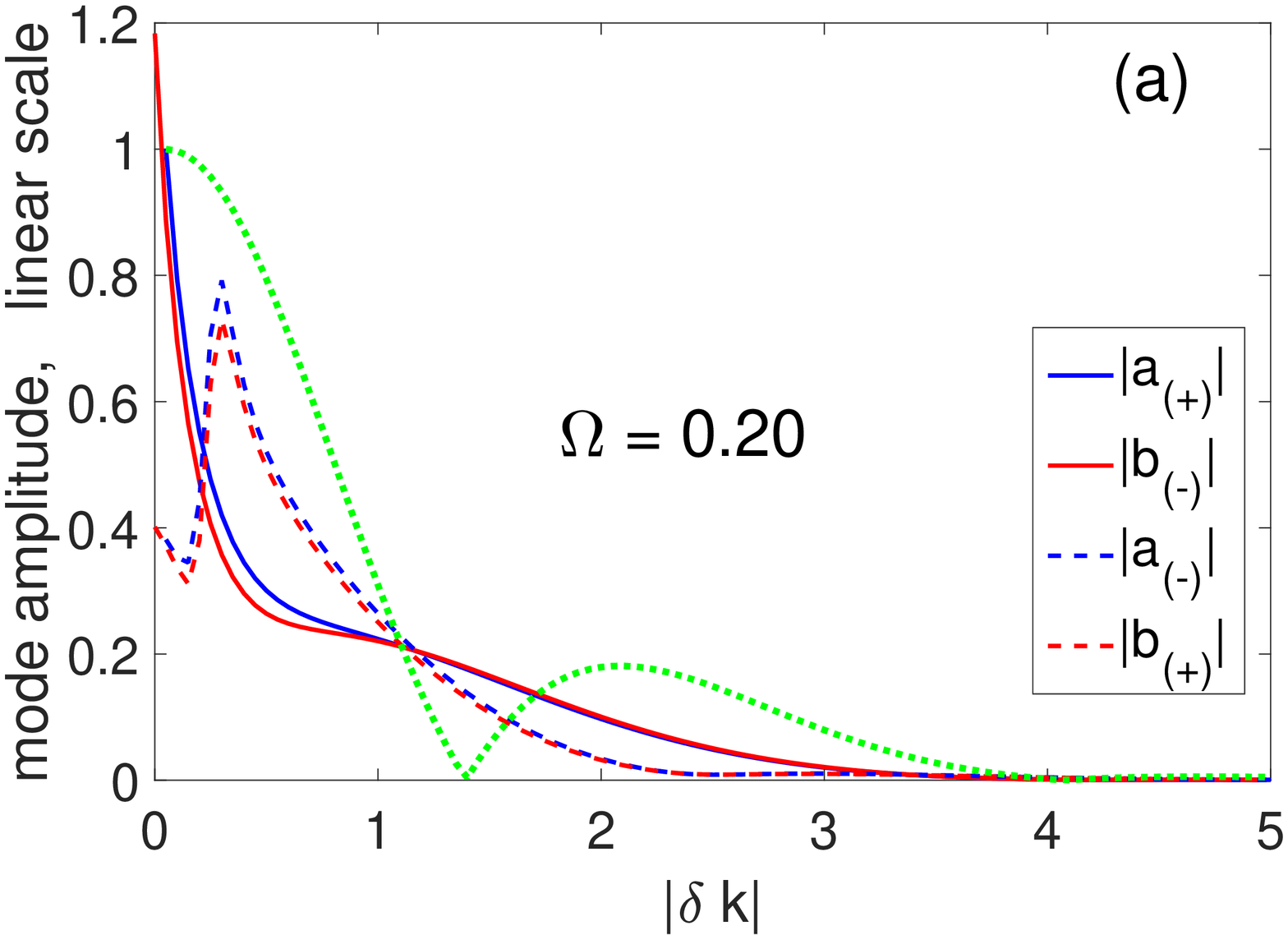}
\end{minipage}
\hspace{0.2cm}
\begin{minipage}{5.2cm}
\hspace*{-0.1cm} 
\includegraphics[height=4.2cm,width=5.1cm,angle=0]{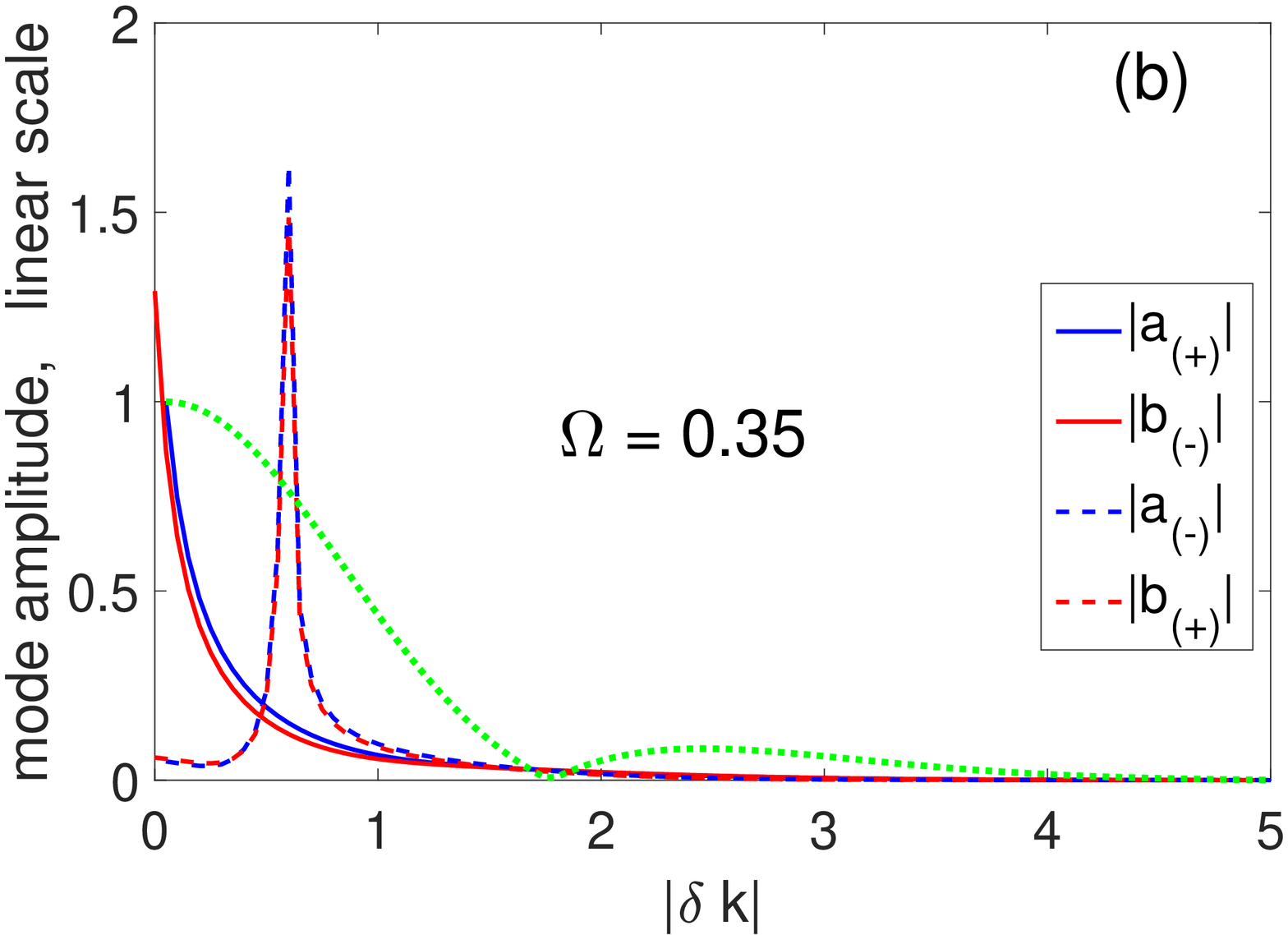}
\end{minipage}
\hspace{0.2cm}
\begin{minipage}{5.2cm}
\hspace*{-0.1cm} 
\includegraphics[height=4.2cm,width=5.1cm,angle=0]{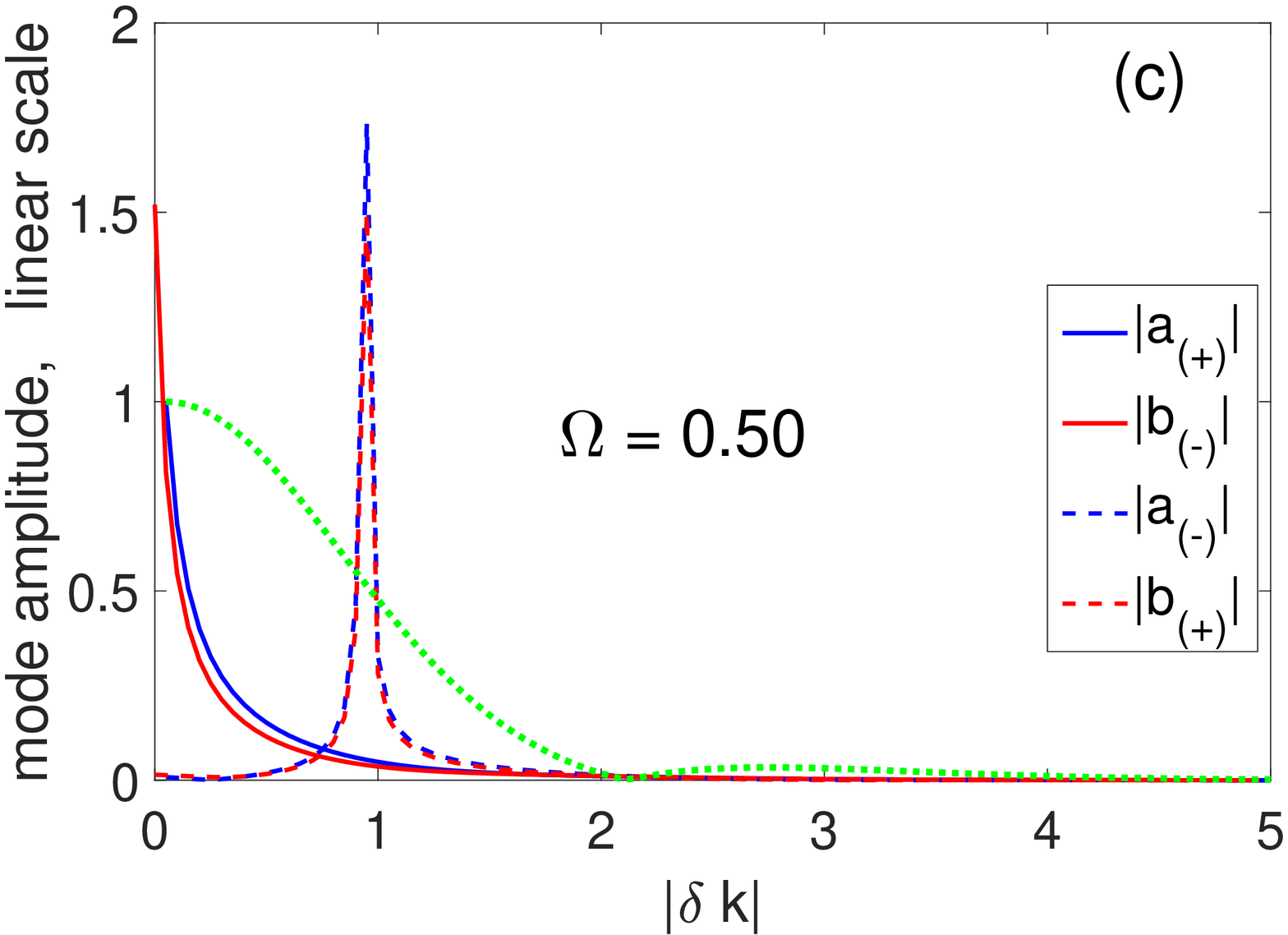}
\end{minipage}

\medskip

\begin{minipage}{5.2cm}
\hspace*{-0.1cm} 
\includegraphics[height=4.2cm,width=5.1cm,angle=0]{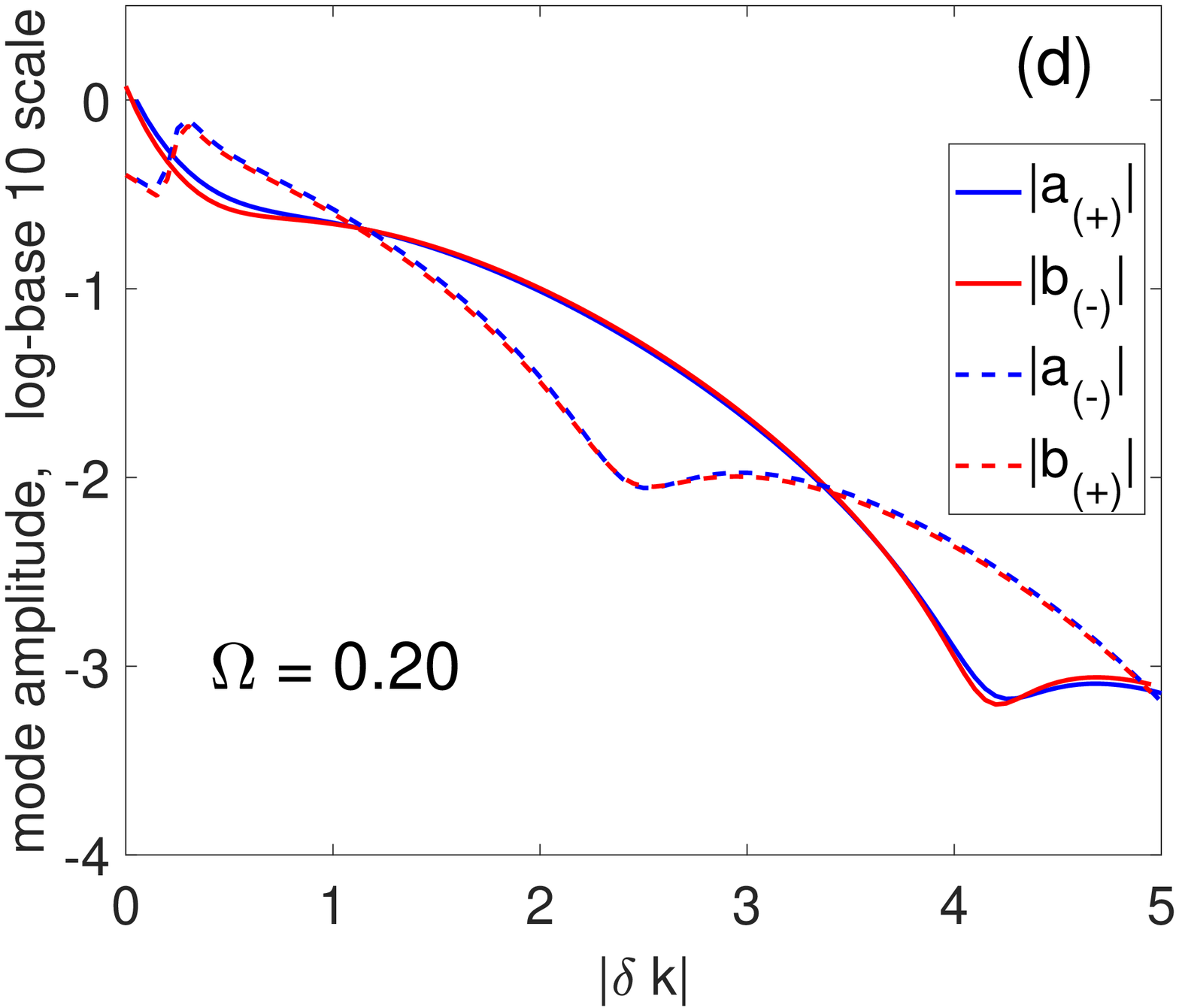}
\end{minipage}
\hspace{0.2cm}
\begin{minipage}{5.2cm}
\hspace*{-0.1cm} 
\includegraphics[height=4.2cm,width=5.1cm,angle=0]{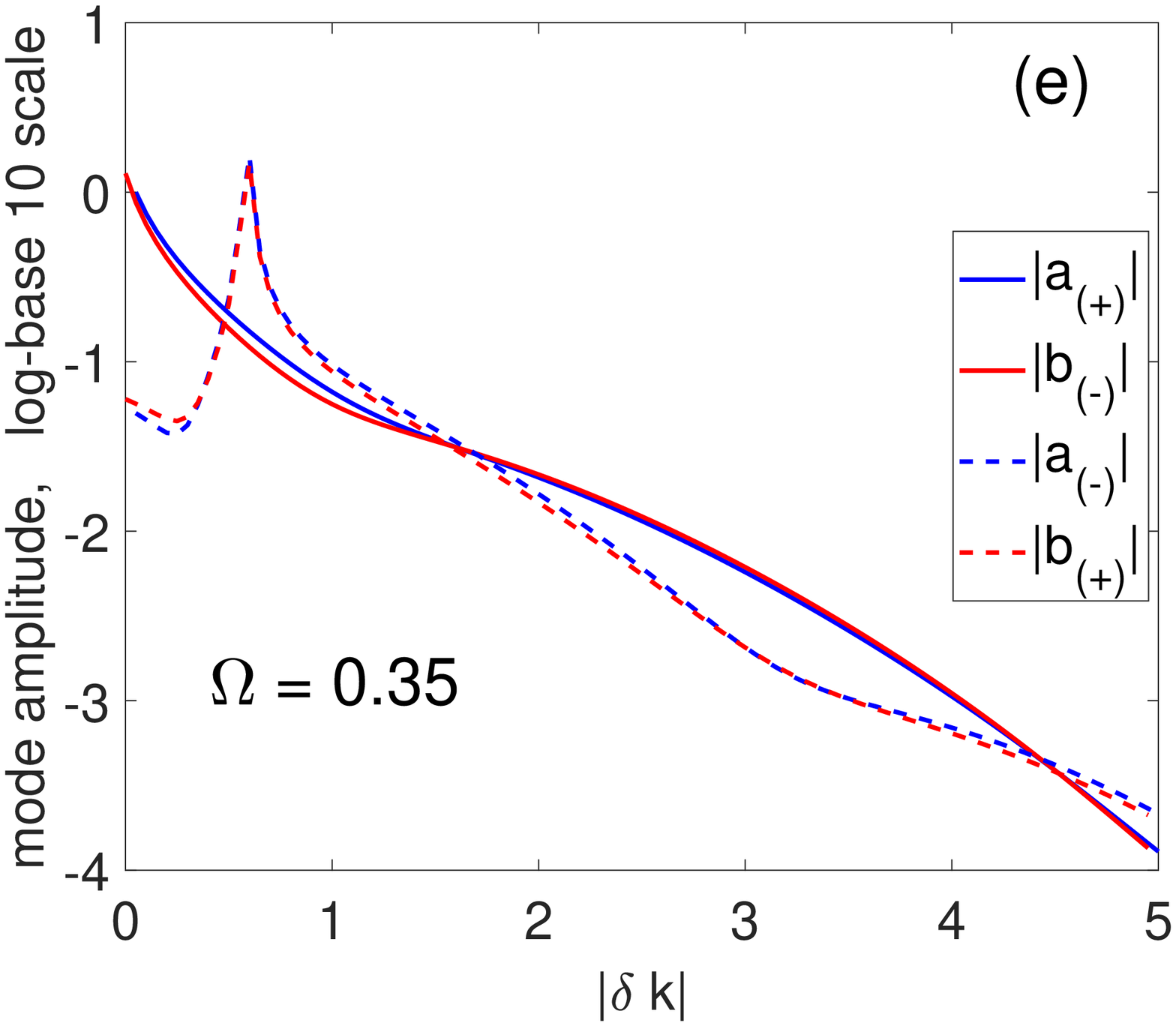}
\end{minipage}
\hspace{0.2cm}
\begin{minipage}{5.2cm}
\hspace*{-0.1cm} 
\includegraphics[height=4.2cm,width=5.1cm,angle=0]{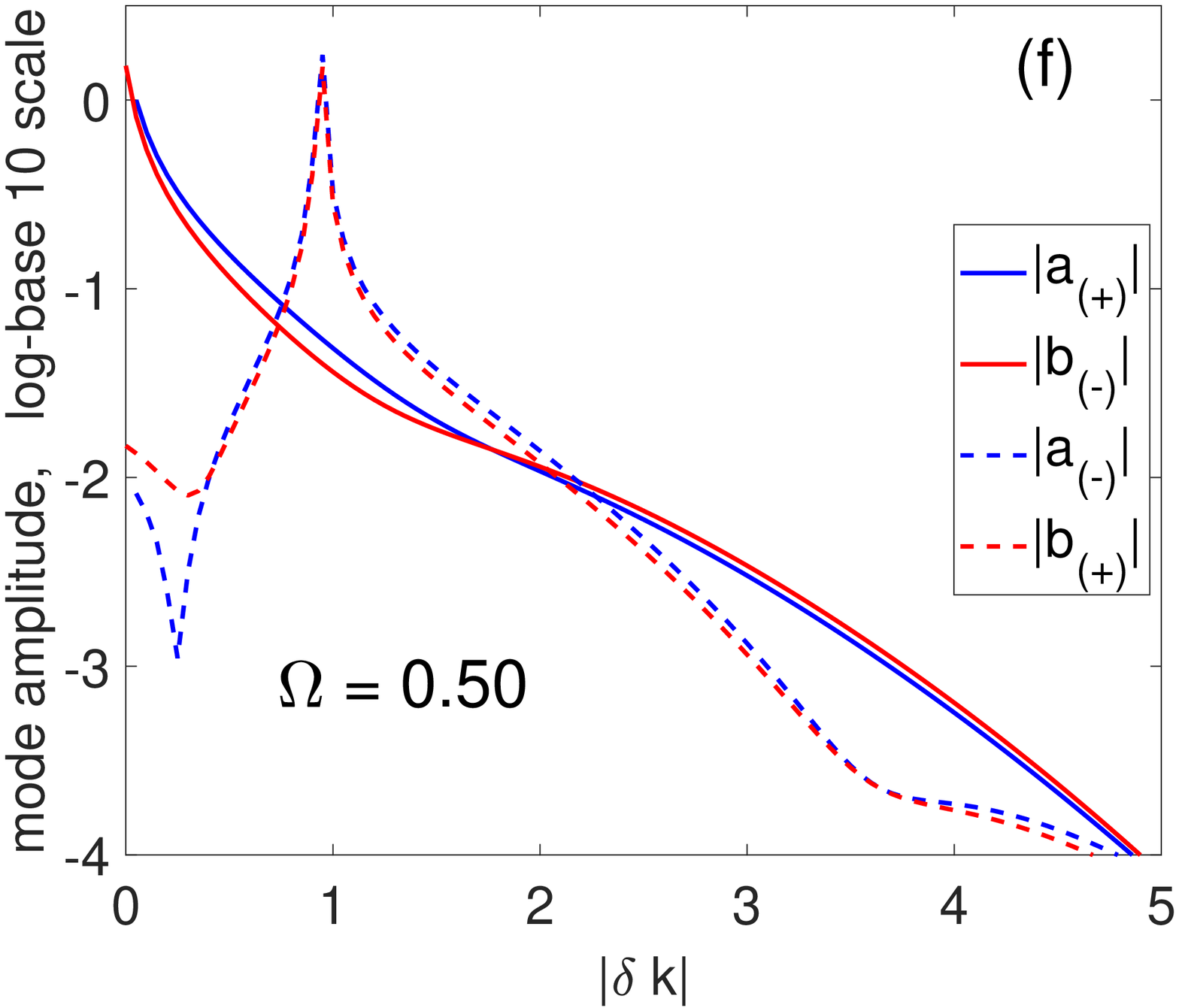}
\end{minipage}
\caption{(Color online) \ Shapes of the most unstable eigenvectors of \eqref{e5_07},
$\big|\und{\widehat{a}}_{(\pm)}\big|$ vs. $\delta k^{(<0)}$ and 
$\big|\und{\widehat{b}}_{(\pm)}\big|$ vs. $\delta k^{(\ge 0)}$,
 for the 
solitons with the same values of $\Omega$ as shown in Fig.~\ref{fig_5}(a).
Panels (a)--(c): linear scale; \ panels (d)--(f): logarithmic scale. 
All eigenvectors are normalized to $\max \big|\und{\widehat{a}}_{(+)}\big|$.
In (a)--(c), the dotted green line shows $\widehat{P_0}(k)$, normalized to 1,
for comparison. 
}
\label{fig_7}
\end{figure}

\begin{figure}[!ht]
\begin{minipage}{5.2cm}
\hspace*{-0.1cm} 
\includegraphics[height=4.2cm,width=5.1cm,angle=0]{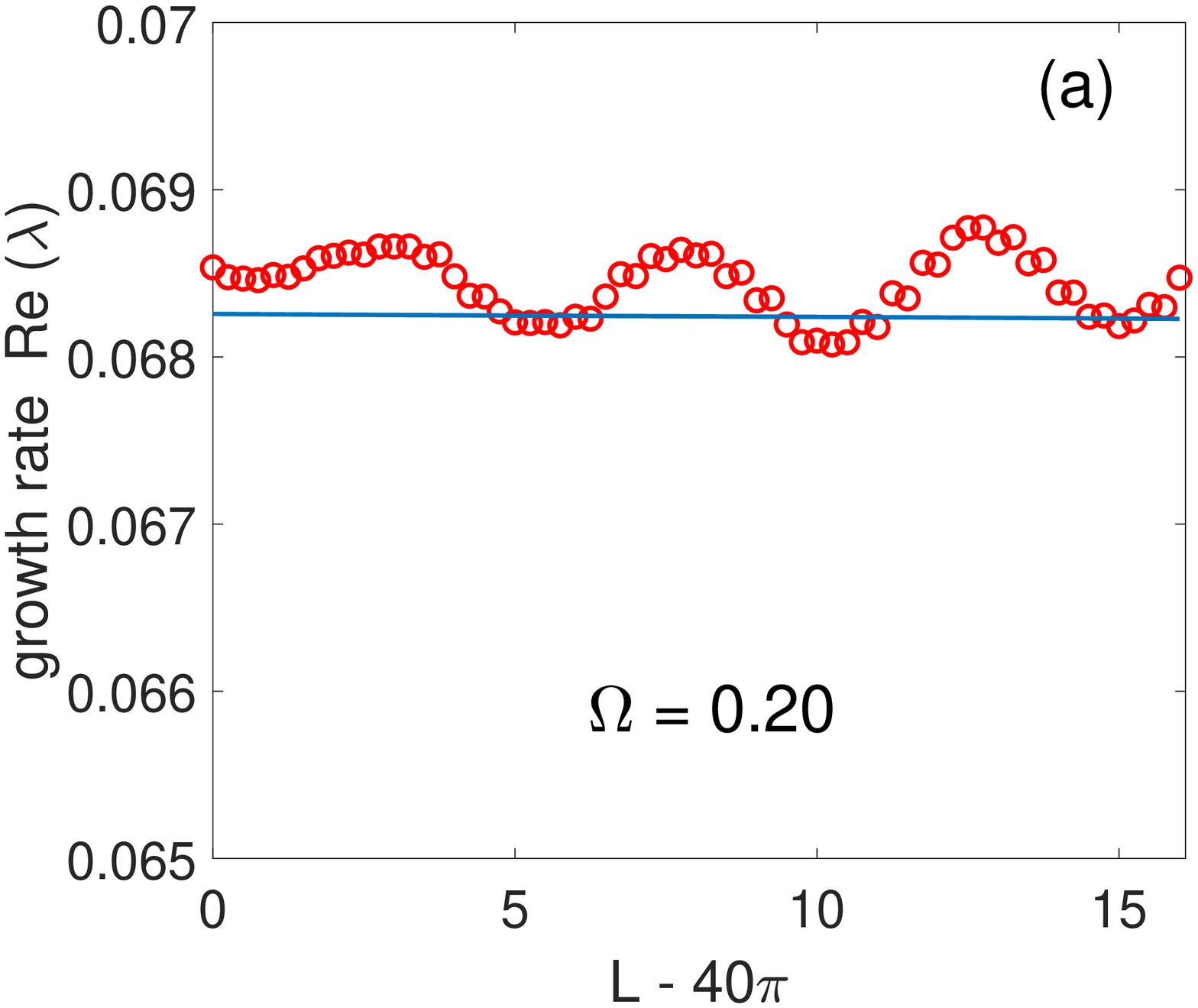}
\end{minipage}
\hspace{0.2cm} 
\begin{minipage}{5.2cm}
\hspace*{-0.1cm} 
\includegraphics[height=4.2cm,width=5.1cm,angle=0]{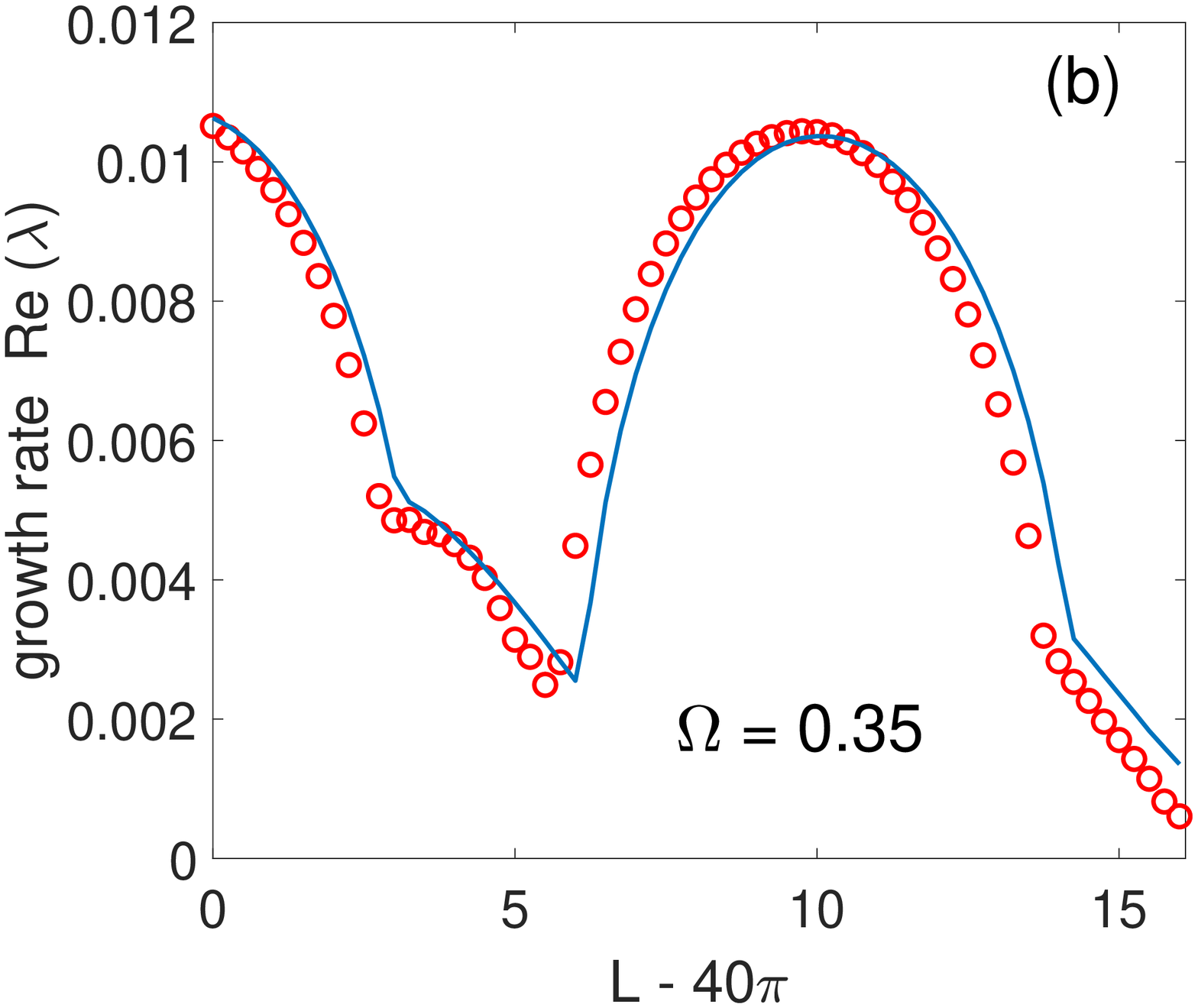}
\end{minipage}
\hspace{0.2cm} 
\begin{minipage}{5.2cm}
\hspace*{-0.1cm} 
\includegraphics[height=4.2cm,width=5.1cm,angle=0]{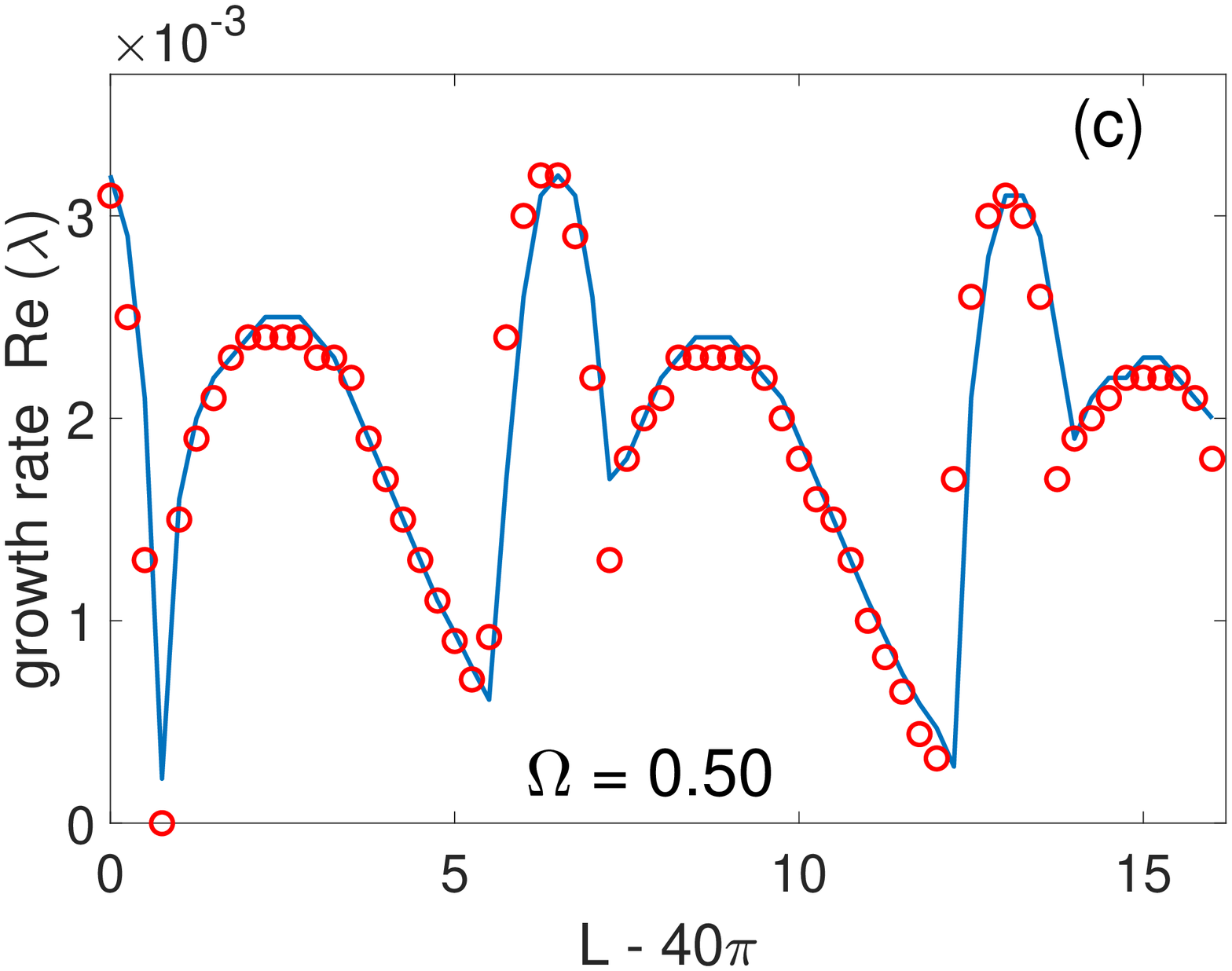}
\end{minipage}
\caption{Growth rate of the most unstable mode of eigenvalue problem
\eqref{e5_07} as the function of length $L$ of the computational domain.
Note that  the vertical scales 
are drastically different in all three panels.
Solid lines correspond to the solution of \eqref{e5_07}, while circles
are the result of simulations by the SSM. 
}
\label{fig_8}
\end{figure}

The stability analysis of system \eqref{e5_07} is straightforward:
one seeks $\und{\mathbbm{s}}$ proportional to $\exp[\lambda t]$, whence
the system becomes a $4M\times 4M$ eigenvalue problem. 
For $M \sim 100$, it is solved by Matlab in about 1 second. 
Eigenvalues
with Re$\lambda>0$ correspond to numerically unstable modes, which,
by design, occur near the edges of the computational spectrum. 
In Fig.~\ref{fig_7} we show the 
unstable eigenvectors of that system corresponding to the solitons
with the three values of $\Omega$ for which graphs are shown in
Fig.~\ref{fig_5}(a). One can see a good agreement between the spectral
profiles of the unstable modes in that figure and their counterparts
in Fig.~\ref{fig_7}(d)--(f). Furthermore, 
in Fig.~\ref{fig_8} we show the NI growth rate
(i.e., \ $\max\,{\rm Re}\lambda$) versus the length $L$ of the 
computational domain for the same three values of $\Omega$. 
The values obtained from the preceding analysis match closely 
the corresponding values measured in direct numerical simulations,
thus {\em validating our analysis}.


Two clarifications about Fig.~\ref{fig_8} are in order. 
First, the dependence of the eigensolutions of problem
\eqref{e5_07} on $L$ occur only via $\dk=2\pi/L$. 
Second, the growth rate deduced from direct numerical simulations was
computed as follows. For a given $\Omega$, the simulations were run
up to the respective time indicated in the caption to Fig.~\ref{fig_5}. 
The temporal evolution of the logarithm of the maximum amplitude of 
Fourier harmonics in some 
vicinity of $k=\kmax$ was plotted, as in Fig.~\ref{fig_4}(c). 
The growth rate was deduced from the slope of the linear part of that graph.

\begin{figure}[!ht]
\hspace*{0.4cm}
\begin{minipage}{7.5cm}
\hspace*{-0.1cm}
\includegraphics[height=5.6cm,width=7.5cm,angle=0]{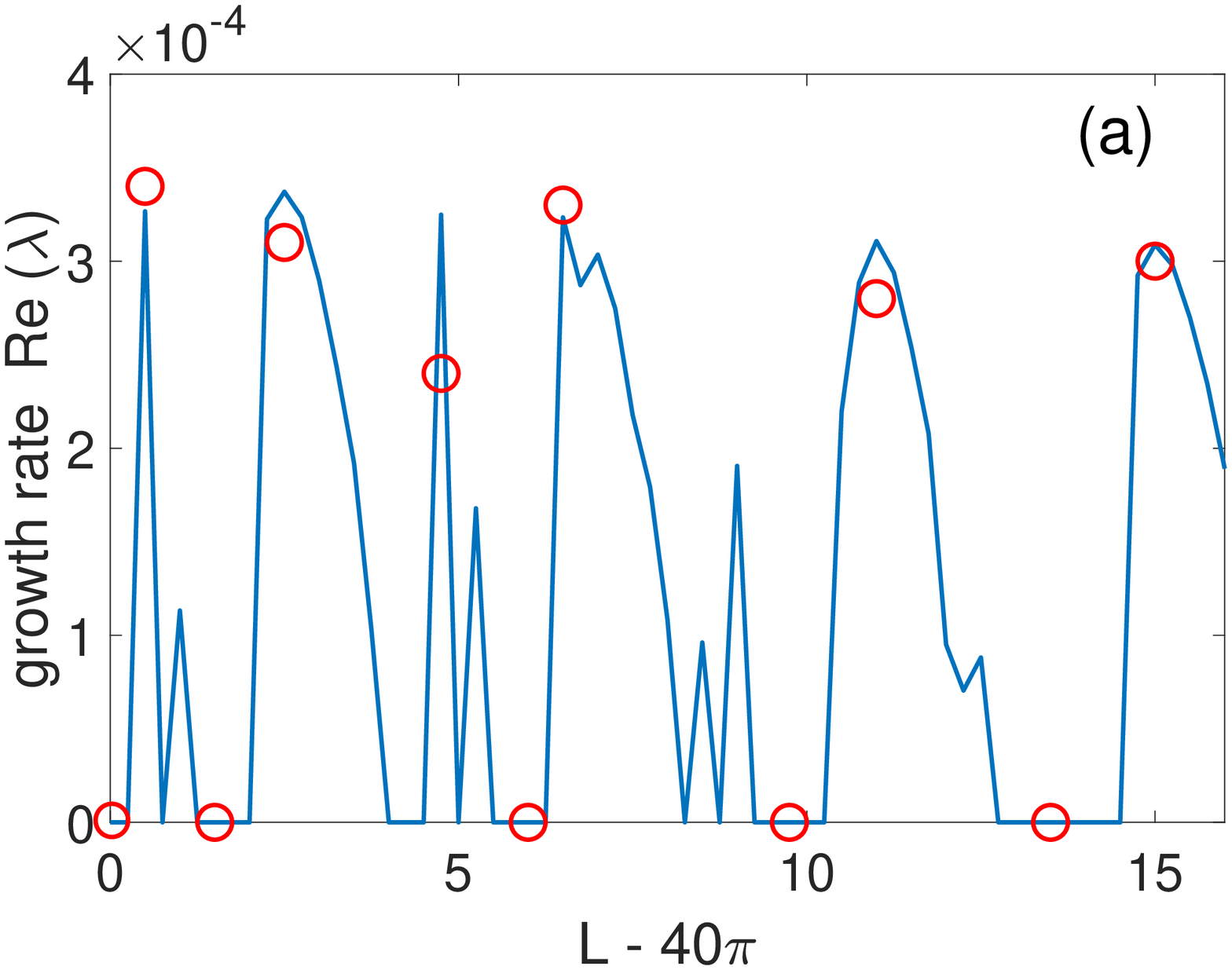}
\end{minipage}
\hspace{0.5cm}
\begin{minipage}{7.5cm}
\hspace*{-0.1cm}
\includegraphics[height=5.6cm,width=7.5cm,angle=0]{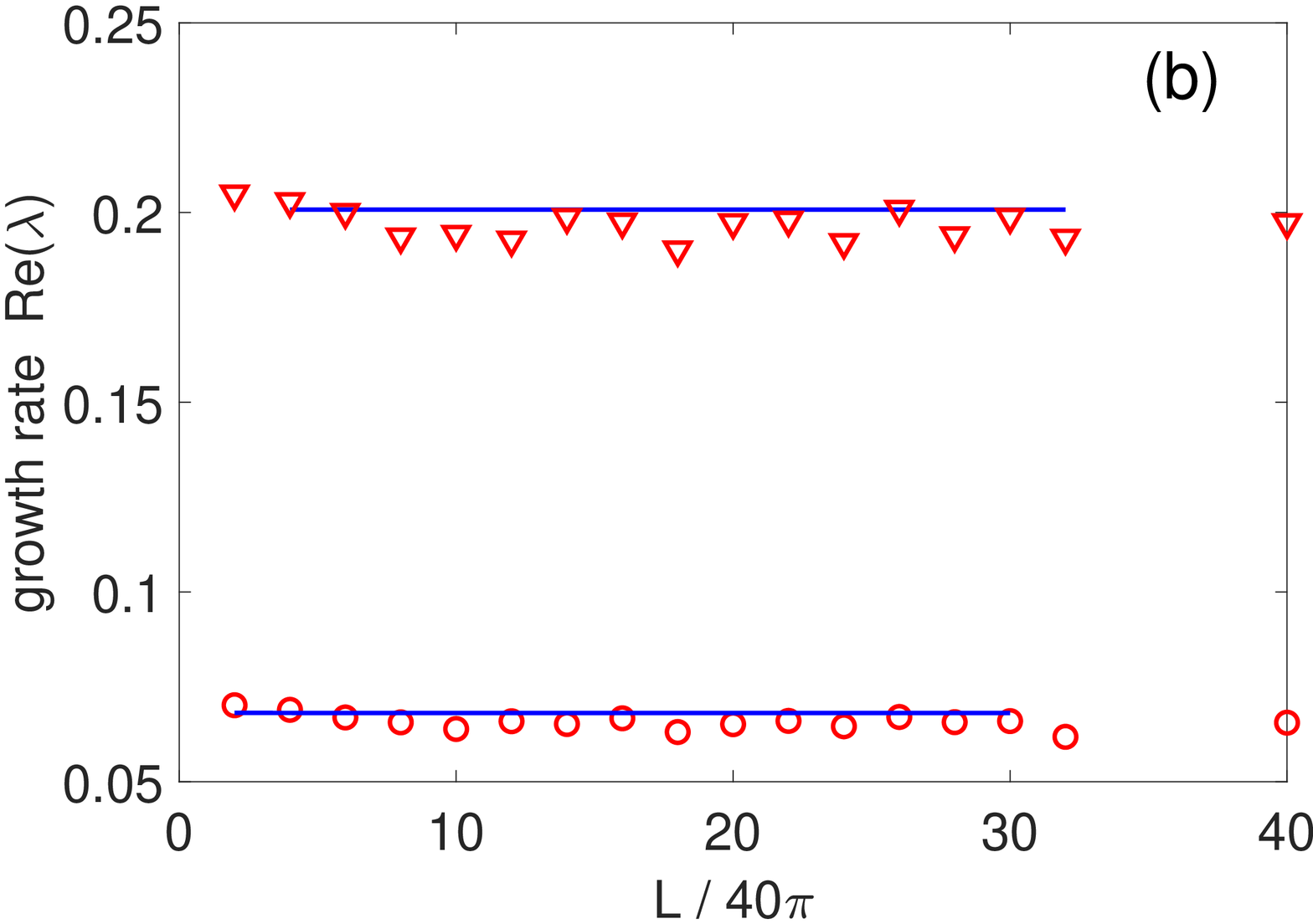}
\end{minipage}
\caption{
Growth rate of the most unstable mode at the spectral edges
as the function of length $L$ of the computational domain.
Note that 
   {\em both} vertical and horizontal scales are
drastically different in panels (a) and (b).
Solid line: solution of \eqref{e5_07}. Symbols: simulations 
by the SSM with $\dt=\dx/5$. \ 
Panel (a): \ $\Omega=0.75$, $N=2^{12}$, and $t=50,000$. Because of the large
simulation time, which in turn is due to a very small growth rate,
simulations were run only for the major ``peaks" and for
the midpoints of the intervals where the theory predicts zero growth rate. \ \ 
Panel (b): \ $\Omega=0.2$ and $0.1$. Simulations by the SSM were run to 
$t=300$ (circles) and $150$ (triangles), respectively. 
Number of grid points was adjusted as
$N=2^{12}\,L/(40\pi)$. 
}
\label{fig_9}
\end{figure}

In Fig.~\ref{fig_9}(a) we show the NI growth rate as a function
of $L$ of a non-fragile soliton with $\Omega=0.75$. 
Again, a good agreement between our analysis and direct numerics is seen.
Thus, the SSM can be unconditionally unstable
even for non-fragile solitons; it is just that in that case, the instability is
so weak that it will not affect simulations for most realistic simulation times. 

Finally, in Fig.~\ref{fig_9}(b) we expand on the trend which may already be noticed
from Fig.~\ref{fig_8}(a). Namely, for sufficiently small $\Omega$, the growth rate
of the unstable modes near $\pm \kmax$ is nearly independent of the length of the 
computational domain. We will discuss this further in Section 8. 

\section{Unconditional instability of the ``noise floor"}

The main result of this section is based on Eq.~\eqref{e6_10}, which describes the
evolution of the ``noise floor" amplitude. The growth rate of the ``noise floor" 
is found from the spectral radius of the corresponding fundamental solution,
defined in \eqref{e6_11b}.

To study the ``noise floor" NI, described in Section 3.3 and illustrated in Fig.~\ref{fig_6},
one cannot use ansatz \eqref{e4_03}. Indeed, a perturbation described by the ``noise
floor" is not spectrally narrow and hence is not a slowly-modulated plane wave in the
$x$-space; see the sentence after \eqref{e4_03}. Instead, we will consider the most general
form of a perturbation:
\bsube
\be
\bm{\widetilde{\psi}}_{\{n\}}(x) \,=\, 
\sum_{j} 
 \bm{\widehat{\alpha}}_{j,\,\{n\}} e^{ik_j x} + 
 \bm{\widehat{\beta}}_{j,\,\{n\}} e^{-ik_j x},
\label{e6_01a}
\ee
where the double subscript in \eqref{e6_01} stands
for a Fourier harmonic at wavenumber $k_j$ and time level $n\dt$; a similar
notation will be employed below in this section.  
The summation here is assumed over the harmonics of the ``noise floor", i.e.
those with 
\be
k_j\gg 1 \quad {\rm and} \quad (\kmax - k_j)\gg 1. 
\label{e6_01b}
\ee
\label{e6_01}
\esube
The first inequality in \eqref{e6_01b} ensures
that the harmonics are outside the spectrum of the background soliton, whereas
the second one ensures that they are sufficiently far from the edges of the
spectral domain.
While the derivation of the equations predicting unstable dynamics of
$\bm{\widehat{\alpha}}_{j}$ and $\bm{\widehat{\beta}}_{j}$ will generally
follow the derivation of Section 4, we will specifically emphasize places 
where key differences occur.

Substituting \eqref{e2_03} and \eqref{e6_01} into \eqref{e4_01} and
linearizing, one obtains equations analogous to \eqref{e4_04}:
\bsube
\be
\bm{\widehat{\alpha}}_{j,\,\{n+1\}} e^{-i\Omega\dt} = 
 e^{-ik_j\sone\dt} \left( \bm{\widehat{\alpha}}_{j,\,\{n\}} + 
 \frac{i\dt}{N} \sum_l \widehat{\bP}_l \bm{\widehat{\alpha}}_{j-l,\,\{n\}} + 
                       \widehat{\bQ}_l (\bm{\widehat{\beta}}_{j-l,\,\{n\}})^* \,
								   \right),
\label{e6add1_01a}
\ee
\be
\bm{\widehat{\beta}}_{j,\,\{n+1\}} e^{-i\Omega\dt} = 
 e^{ik_j\sone\dt} \left( \bm{\widehat{\beta}}_{j,\,\{n\}} + 
 \frac{i\dt}{N} \sum_l \widehat{\bQ}_l (\bm{\widehat{\alpha}}_{j+l,\,\{n\}})^* + 
                       \widehat{\bP}_l \bm{\widehat{\beta}}_{j+l,\,\{n\}} \,
								   \right)\,.
\label{e6add1_01b}
\ee
\label{e6add1_01}
\esube
The sums on the r.h.s.~of these equations were obtained similarly to those in
\eqref{C_01} and \eqref{C_03}; however, for a reason that will become clear soon,
the indices were switched between the \ $\widehat{\bP},\widehat{\bQ}$ \ and the \ 
$\bm{\widehat{\alpha}},\bm{\widehat{\beta}}$ \ terms
in \eqref{e6add1_01} 
compared to those in \eqref{C_01} and \eqref{C_03}. 
Since the spectral width
of the background soliton is $O(1)$, 
   the sum in \eqref{e6add1_01} contains
   only $O(1)/\dk \ll N$ terms that are significantly different from zero.

Similarly to \eqref{e4_07}, one seeks
\be
\bm{\widehat{\alpha}}_{j,\{n\}} = e^{-in\sone k_j\dt}\,
 \left( \widehat{\bma}_{j,\{n\}}^{(0)} + 
     \widehat{\bma}_{j,\{n\}}^{(1)} \right) \,, 
\qquad 
\bm{\widehat{\beta}}_{j,\{n\}} = e^{in\sone k_j\dt}\,
 \left( \widehat{\bmb}_{j,\{n\}}^{(0)} + 
  \widehat{\bmb}_{j,\{n\}}^{(1)} \right) \,,
\label{e6_02}
\ee
where the quantities with superscript `$(0)$' vary slowly in time and those
with superscript `$(1)$' are small. 
Note, however, that the time dependence explicitly stated in \eqref{e6_02} is 
different from that in \eqref{e4_07}. In the latter case, since the
perturbation $\bm{\widetilde{\psi}}(x)$ was spectrally localized near
some wavenumber $k_0$, it was appropriate to assume that in the main 
order, all harmonics evolved proportionally to the same factor,
either $\exp[-in\sone k_0 \dt]$ or $\exp[in\sone k_0 \dt]$. 
On the other hand, since the perturbation \eqref{e6_01}, considered in
this section, is {\em not} spectrally localized, then the principal evolution
of each Fourier harmonic followed its individual exponential,
$\exp[\mp in\sone k_j \dt]$. A consequence of this difference will appear
in the subsequent derivation.

When one substitutes \eqref{e6_02} into \eqref{e6add1_01}, one obtains,
similarly to \eqref{e4_08}, two distinct groups of oscillating terms. 
For example, the $\widehat{\bP}_l \bm{\widehat{\alpha}}_{j-l,\,\{n\}}$ term yields:
\be
e^{in\sone k_j n\dt}\, \widehat{\bP}_l \,\bm{\widehat{\alpha}}_{j-l,\,\{n\}} = 
\left( \widehat{\bP_{01}}_{\;l} \;e^{in\sone k_l n\dt} + 
  \widehat{\bP_{23}}_{\;l} \;e^{in\sone (k_l-2k_j) n\dt} \right)
	\left( {\widehat{\bma}}_{j-l,\,\{n\}}^{(0)} + {\widehat{\bma}}_{j-l,\,\{n\}}^{(1)}
	 \right)\,.
\label{e6_03}
\ee
The terms in the first group on the r.h.s.~vary on the time scale of order $O(1)$
(see the second sentence after \eqref{e6add1_01}), while the terms in the
second group vary much faster due to the first inequality in 
\eqref{e6_01b}.\footnote{
   Note that since we intend to investigate an unconditional numerical instability,
   one can consider the limit $\dt \To 0$, and thus 
   there can be no resonances like \eqref{e3_01}, because of which any 
	 of the terms from the second group may become slowly varying.}
Then, arguing as in the paragraph surrounding Eqs.~\eqref{e4_09}, one shows 
that the rapidly oscillating terms affect only the small corrections  
$\widehat{\bma}_{j,\{n\}}^{(1)}$, $\widehat{\bmb}_{j,\{n\}}^{(1)}$,
but in the main order do not contribute to the evolution of the 
principal terms 
$\widehat{\bma}_{j,\{n\}}^{(0)}$, $\widehat{\bmb}_{j,\{n\}}^{(0)}$.
Therefore, in what follows we omit those rapidly oscillating terms
and will also omit the superscript $(0)$, as done after \eqref{e4_10}. 
Following the above steps and also approximating the finite differences
in time with time derivatives, as done in obtaining \eqref{e5_04} from 
\eqref{e5_03}, we find:
\bsube
\be
(\widehat{\bma}_{j})_t  = i\Omega \widehat{\bma}_{j} + 
 \frac{i}{N} \sum_{l} 
 \widehat{\bP_{01}}_{\;l}  \,e^{i\sone k_l t }\, \widehat{\bma}_{j-l} + 
 \widehat{\bQ_{01}}_{\;l}  \,e^{i\sone k_l t }\, (\widehat{\bmb}_{j-l})^*\,,
\label{e6_04a}
\ee
\be
(\widehat{\bmb}_{j} )^*_t = -i\Omega (\widehat{\bmb}_{j})^* -  
 \frac{i}{N} \sum_{l} 
 \widehat{\bQ_{01}}^{\,*}_{\;l}  \,e^{-i\sone k_l t }\, \widehat{\bma}_{j+l} + 
 \widehat{\bP_{01}}^{\,*}_{\;l}  \,e^{-i\sone k_l t }\, (\widehat{\bmb}_{j+l})^*\,,
\label{e6_04b}
\ee
\label{e6_04}
\esube
where $t \equiv n\dt$.

Exact solution of system \eqref{e6_04} would be possible only by direct numerical
simulations. Not only would such an approach not be illuminating in any respect, but 
it would also be considerably more difficult than solving \eqref{e5_04} (or,
equivalently, \eqref{e5_07}), despite the latter equation
appearing to have more terms. Indeed, \eqref{e6_04}, unlike \eqref{e5_04},
has time-dependent terms on the r.h.s., which is the consequence of the difference
between \eqref{e6_02} and \eqref{e4_07}, emphasized after \eqref{e6_02}. Even more
importantly, the system in \eqref{e6_04} couples, by virtue of \eqref{e6_01b},
a much greater number of Fourier harmonics, $\lesssim N/2$ instead of $M\ll N/2$. 
For these reasons, below we will use a {\em simplified approach}. It will still require a
numerical solution, but only of a $2\times 2$ system. More
importantly, it will allow us to {\em explain the mechanism} by which harmonics
of the ``noise floor" can become unstable.

The simplification occurs from an observation that none of the coefficients of terms
$\widehat{\bma}$ and $\widehat{\bmb}$ on the r.h.s.~ of \eqref{e6_04} depends on the
harmonic's index $j$. Therefore, these equations will be satisfied by a
$j$-independent ansatz:
\be
\widehat{\bma}_j \equiv \widehat{\bma}, \qquad
\widehat{\bma}_j \equiv \widehat{\bmb} \qquad
\mbox{for all $j$ satisfying \eqref{e6_01b}}.
\label{e6_05}
\ee
We will come back to interpretation of this ansatz later in this Section. 
We now expand $\widehat{\bma}$ and $\widehat{\bmb}$ using \eqref{e5_05},
substitute the result into \eqref{e6_04} and collect the scalar coefficients
of vectors $\bm{e}_{(\pm)}$. For example, the coefficients at $\bm{e}_{(+)}$
yield:
\bsube
\be 
(\wap)_t = i\Omega\wap +\frac{i}{N} \sum_l 
 \widehat{P_{0}}_{\;l}  \,e^{i k_l t }\, \wap + 
 \left( \widehat{Q_{0}}_{\;l}   +  \widehat{Q_{1}}_{\;l}   \right)
       \,e^{i k_l t }\, (\wbp)^*\,,
\label{e6_06a}
\ee
\be
(\wbp )^*_t = -i\Omega (\wbp)^* -  
 \frac{i}{N} \sum_{l} 
 \left( \widehat{Q_{0}}^{\,*}_{\;l}  + \widehat{Q_{1}}^{\,*}_{\;l}  \right)
   \,e^{-i k_l t }\, \wap + 
 \widehat{P_{0}}^{\,*}_{\;l}  \,e^{-i k_l t }\, (\wbp)^*\,,
\label{e6_06b}
\ee
\label{e6_06}
\esube
where we have used that $P_1\equiv 0$. 
In deriving \eqref{e6_06} we have also used the first two of identities \eqref{e5_06} 
as well as the identity
\be
e^{i\sone k_l t } \, \bm{e}_{(\pm)} \,=\, e^{\pm ik_l t}\, \bm{e}_{(\pm)} \,,
\label{e6_07}
\ee
which follows from \eqref{e4_05} and \eqref{e5_06}. Note that unlike in 
\eqref{e5_07}, terms $\wap,\wbp$ are not coupled with $\wam,\wbm$ due to the
absence of matrices $\stwo,\,\sthr$ in \eqref{e6_04}.

Here comes the next {\em key step} in this analysis: we recognize the sums in \eqref{e6_06}
as the inverse Fourier transform \eqref{e4_02} (but in $t$, not in $x$), 
upon which we rewrite Eqs.~\eqref{e6_06}
and their counterparts obtained for $\wam,\wbm$ as:
\bsube
\be
\left( \bm{c}_{\,(\pm)} \right)_t = \bm{R}_{(\pm)} \, \bm{c}_{\,(\pm)}, 
\label{e6_08a}
\ee
where:
\be
\bm{c}_{\,(\pm)} = \left( 
            \ba{c} \widehat{a}_{(\pm)} \\ (\widehat{b}_{(\pm)} )^* \ea \right), 
\qquad
\bm{R}_{(\pm)} \,= \, i\sthr \left( \Omega\szer + 
 \left( \ba{cc} P_0(\pm t) & Q_0(\pm t) \pm Q_1(\pm t) \\ 
                Q_0^*(\pm t) \pm Q_1^*(\pm t) & P_0^*(\pm t) \ea \right)
								 \, \right) \,.
\label{e6_08b}
\ee
Two remarks about the entries of the last matrix in \eqref{e6_08b} are in order, 
both of which are consequence
of the difference emphasized after Eqs.~\eqref{e6_02}.
First, while $P_0$ etc.~were defined in \eqref{e2_04} as functions
of $x$, in the above system they are functions of time. 
Second, due to the periodicity of discrete Fourier transform, these entries
are {\em periodic}, with the period being $L$, the length of the computational domain.
Thus, although in \eqref{e6_08a}, $t\in[0,\,\infty)$, one also requires that
\be
P_0(t+L)=P_0(t), \qquad Q_{0,1}(t+L)=Q_{0,1}(t) \qquad \forall\, t\,.
\label{e6_08c}
\ee
\label{e6_08}
\esube

Before we discuss the solution of \eqref{e6_08}, 
let us point out that we need to solve {\em only one}, not two, systems. This follows
from the parity properties implied by \eqref{e2_04} and \eqref{e2_01}:
\be
P_0(-t) = P_0(t), \qquad 
Q_0(-t) = Q_0(t), \qquad
Q_1(-t) = -Q_1(t),
\label{e6_09}
\ee
whence $\bm{R}_{(-)}=\bm{R}_{(+)}$. Therefore, below we will omit the subscripts
`$(\pm)$' of $\bm{c}$ and $\bm{R}$. 
As yet another simplification,
we note that $P_0^*=P_0$, $Q_0^*=Q_0$, while $Q_1^*=-Q_1$. Finally, using the
explicit form \eqref{e2_04}(c,d) of $P_0,Q_0,Q_1$, we can rewrite \eqref{e6_08} as:
\be
\bm{c}_t \,= \, i\sthr \left( \Omega\szer + \frac12
 \left( \ba{cc} \Psi_1^2-\Psi_2^2 & (\Psi_1 + \Psi_2)^2 \\ 
                (\Psi_1 - \Psi_2)^2 & \Psi_1^2-\Psi_2^2 \ea \right)
								 \, \right) \,\bm{c}\,.
\label{e6_10}
\ee
Recall that here, $\Psi_{1,2}$ are $L$-periodic functions of time. 
Then, 
\bsube
\be
\|\bm{c}(t)\| \le  \|\bm{\Phi}(L) \|^{t/L}\; \|\bm{c}(0)\|\,,
\qquad  t = \mbox{integer}\cdot L,
\label{e6_11a}
\ee
where $\|\cdots\|$ denotes the ${\ell}^2$-norm, and 
the fundamental solution $\Phi(L)$ of \eqref{e6_10} satisfies
\be
\bm{c}(L) = \bm{\Phi}(L)\,\bm{c}(0)\,.
\label{e6_11b}
\ee
\label{e6_11}
\esube

\begin{figure}[!ht]
\begin{minipage}{5.2cm}
\hspace*{-0.1cm} 
\includegraphics[height=4.2cm,width=5.1cm,angle=0]{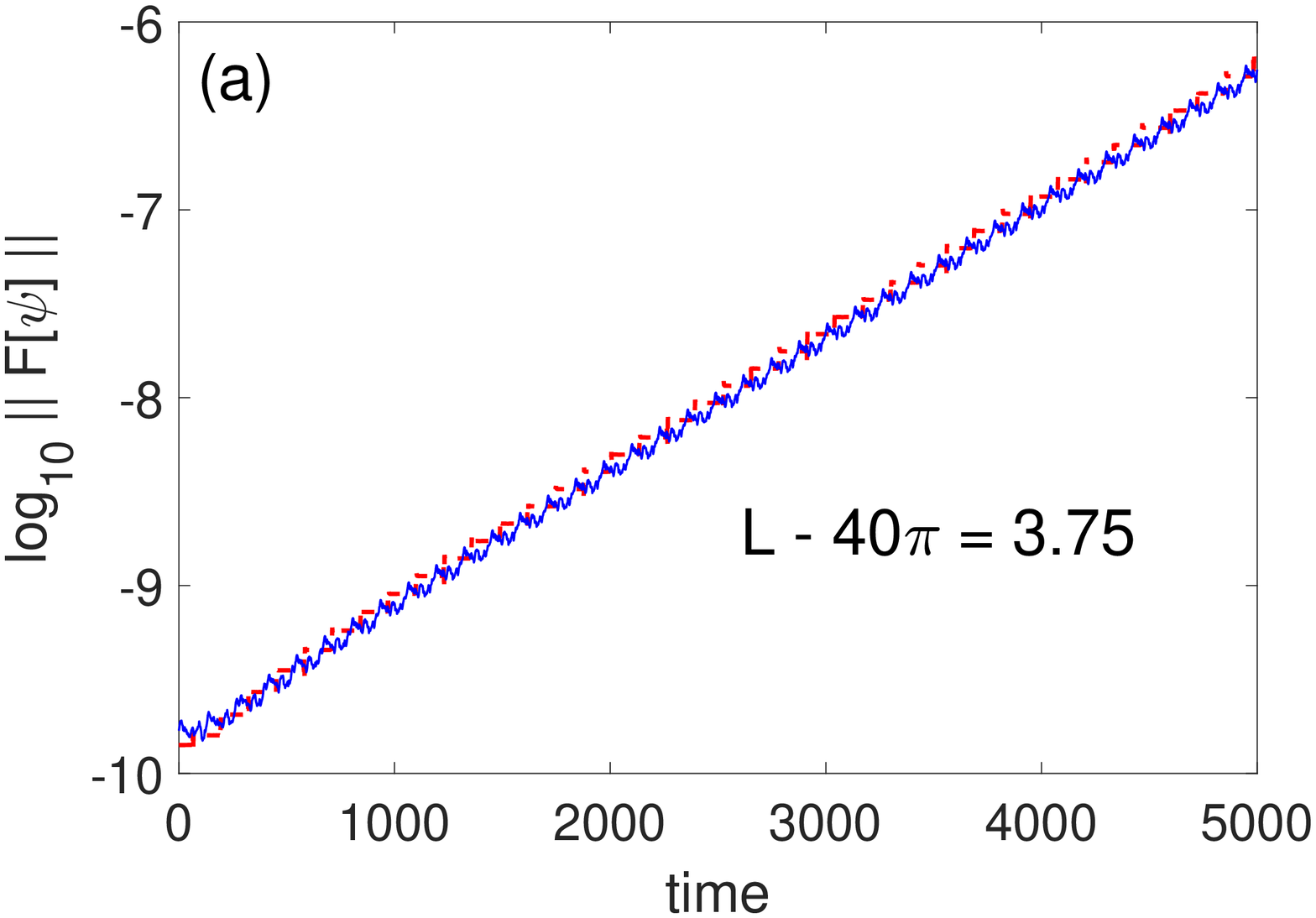}
\end{minipage}
\hspace{0.2cm}
\begin{minipage}{5.2cm}
\hspace*{-0.1cm} 
\includegraphics[height=4.3cm,width=5.1cm,angle=0]{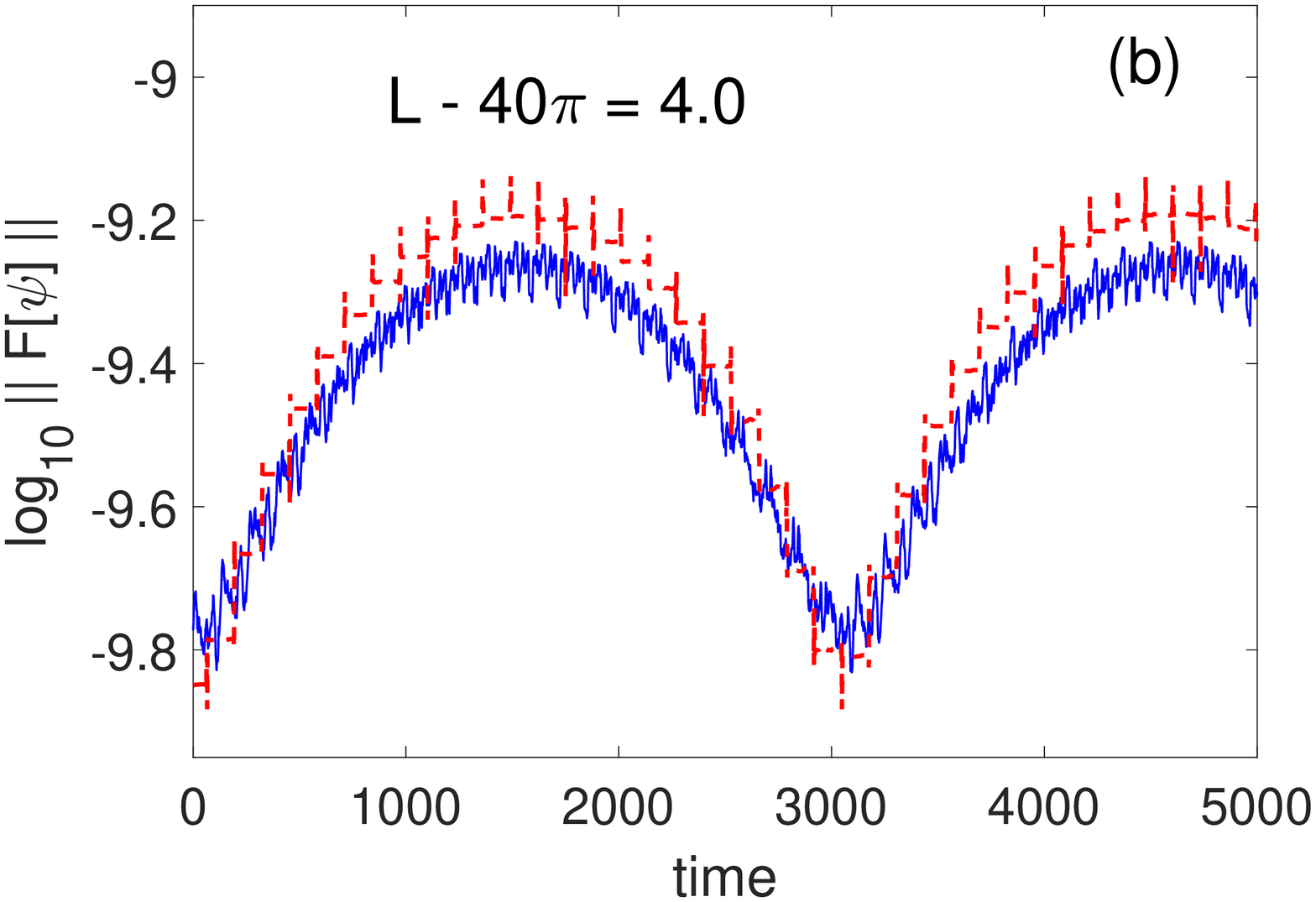}
\end{minipage}
\hspace{0.2cm}
\begin{minipage}{5.2cm}
\hspace*{-0.1cm} 
\includegraphics[height=4.3cm,width=5.1cm,angle=0]{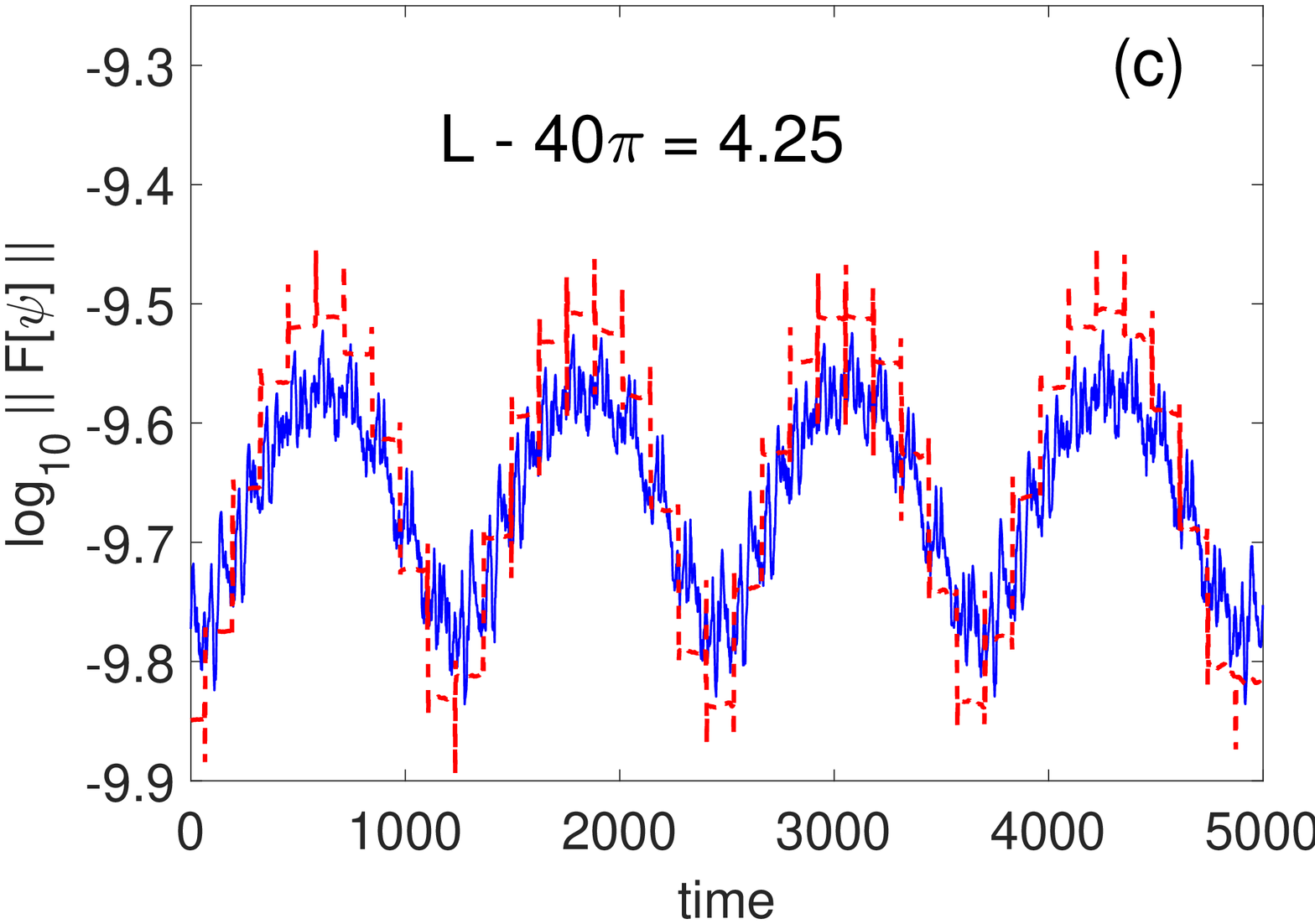}
\end{minipage}
\caption{Evolution of the amplitude of harmonics of the ``noise floor",
similar to that shown in Fig.~\ref{fig_6}(b), but for $\Omega=0.50$ and the
three values of $L$ indicated in the panels. Solid
and dashed lines correspond to a small white noise or a constant, respectively, 
being added to each Fourier mode in the initial condition.
}
\label{fig_10}
\end{figure}

We are now ready to interpret the meaning of ansatz \eqref{e6_05}. A perturbation
\eqref{e6_01} where amplitudes of the harmonics satisfy \eqref{e6_05} is 
approximately\footnote{This would have been exact if the summation had extended for
                all $j\in[1,\,N/2]$.}
the sum of {\em two} delta functions in space, 
where the `two' occurs due to \eqref{e5_05}
having two  contributions, from $\bm{e}_{(+)}$ and $\bm{e}_{(-)}$. 
By virtue of \eqref{e6_02}, these two spikes move with speed 1
in opposite directions. Due to the periodicity of the boundary conditions, they 
repeatedly leave and re-enter the computational domain. Their amplitudes are changed
when they pass through the soliton and remain constant far from the soliton; this
follows from \eqref{e6_10}. If those amplitude changes from consecutive passages 
are accumulated, the amplitude of the ``noise floor" increases; 
see the dashed line in Fig.~\ref{fig_6}(b) 
and Fig.~\ref{fig_10}(a). 
When the initial ``noise floor" consists of white noise as opposed to 
the simplified ansatz \eqref{e6_05},
the above interpretation no longer applies in the exact sense. 
However, the mechanism of the instability
of the ``noise floor" is the same: the perturbation \eqref{e6_01}
repeatedly passes through the soliton, and when $L$ is such that 
changes of its amplitude over
consecutive passages accumulate, the perturbation grows on average 
exponentially. 
When they do not, ``noise floor"'s amplitude oscillates in time
(Figs.~\ref{fig_10}(b,c)). 
The fact that the evolution of the ``simplified" perturbation 
\eqref{e6_01}, \eqref{e6_05} is predictive of that of the generic perturbation
is seen in Figs.~\ref{fig_6}(b) and \ref{fig_10}.

\begin{figure}[!ht]
\begin{minipage}{5.2cm}
\hspace*{-0.1cm} 
\includegraphics[height=4.2cm,width=5.1cm,angle=0]{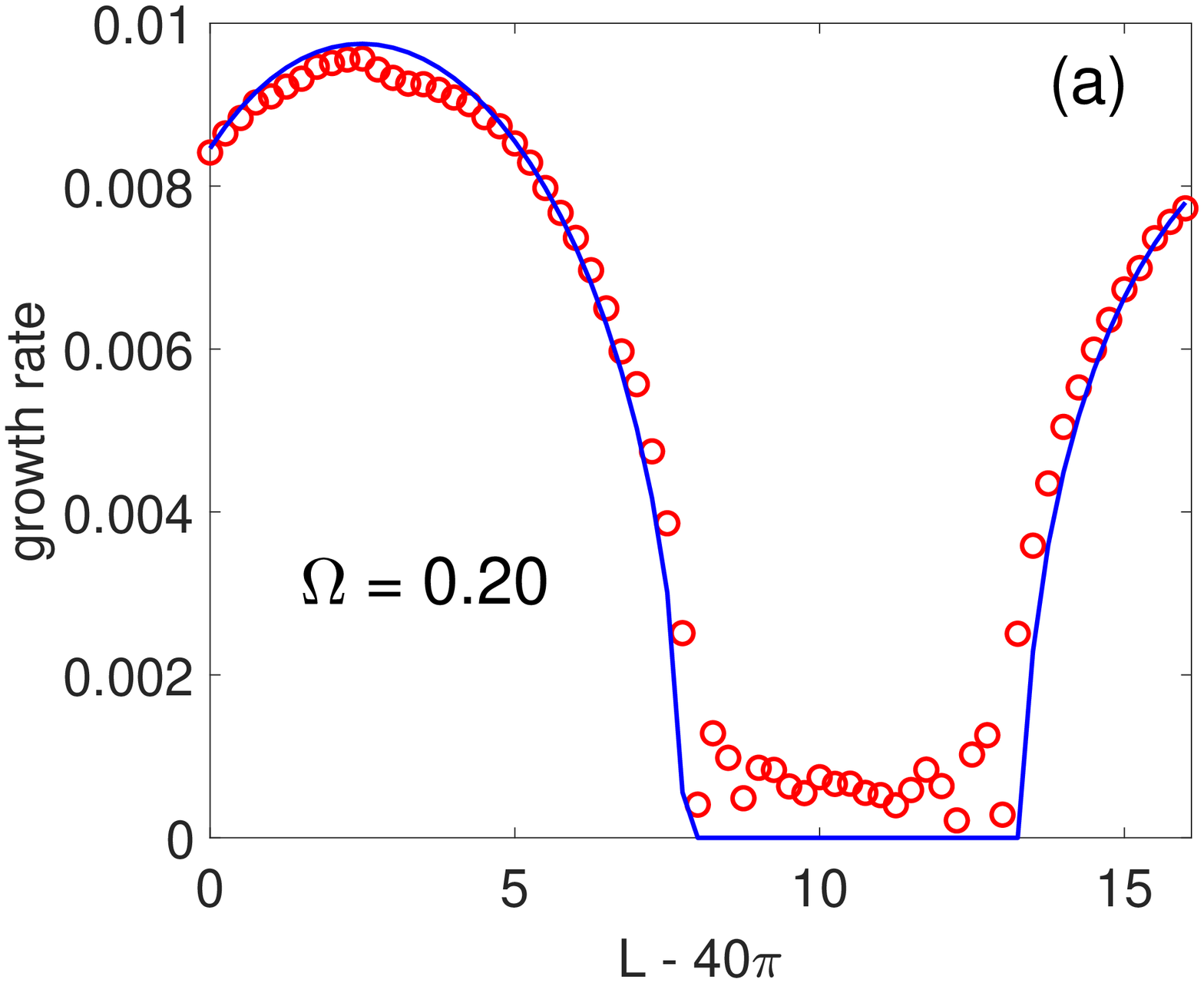}
\end{minipage}
\hspace{0.2cm}
\begin{minipage}{5.2cm}
\hspace*{-0.1cm} 
\includegraphics[height=4.3cm,width=5.1cm,angle=0]{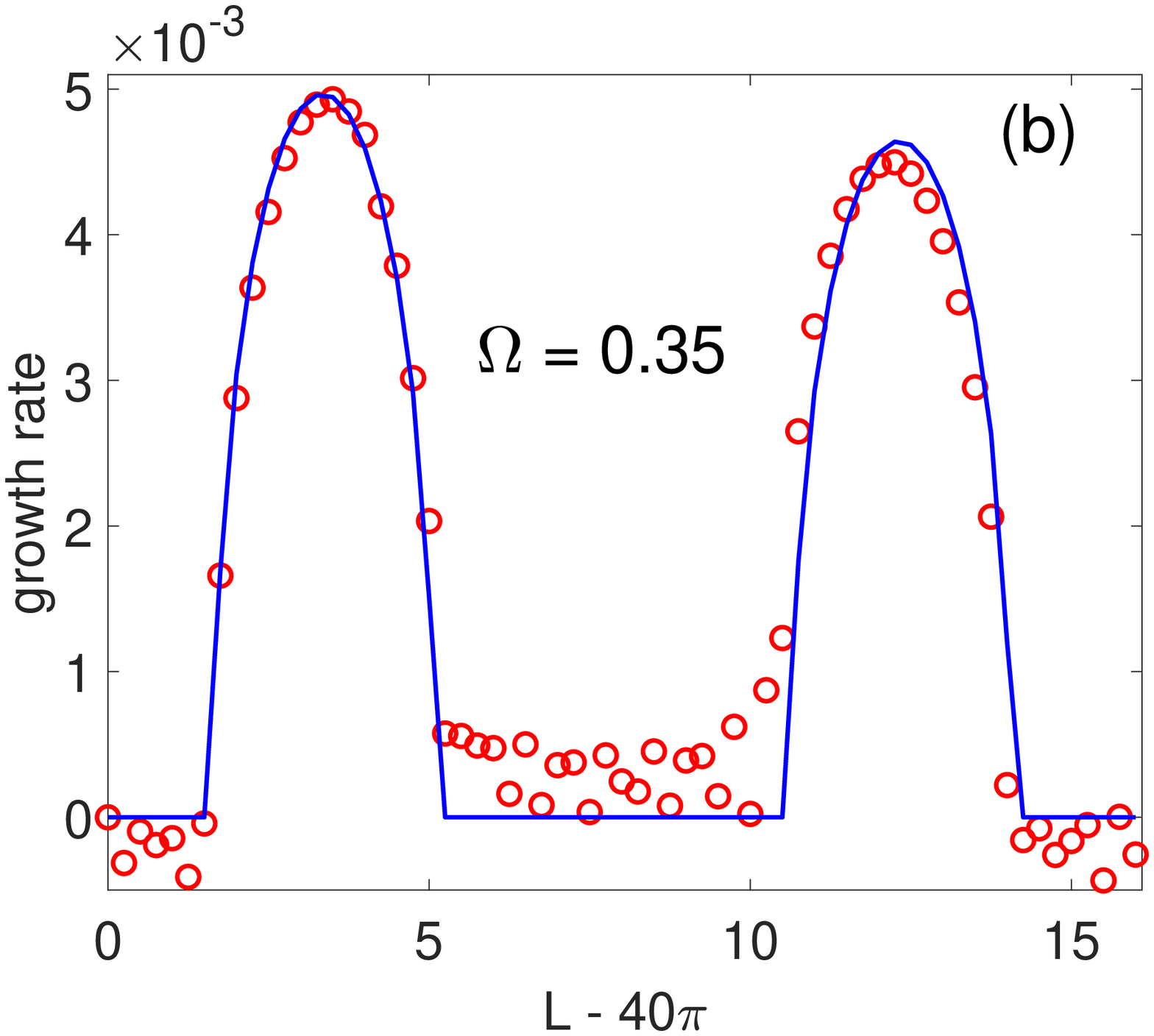}
\end{minipage}
\hspace{0.2cm}
\begin{minipage}{5.2cm}
\hspace*{-0.1cm} 
\includegraphics[height=4.3cm,width=5.1cm,angle=0]{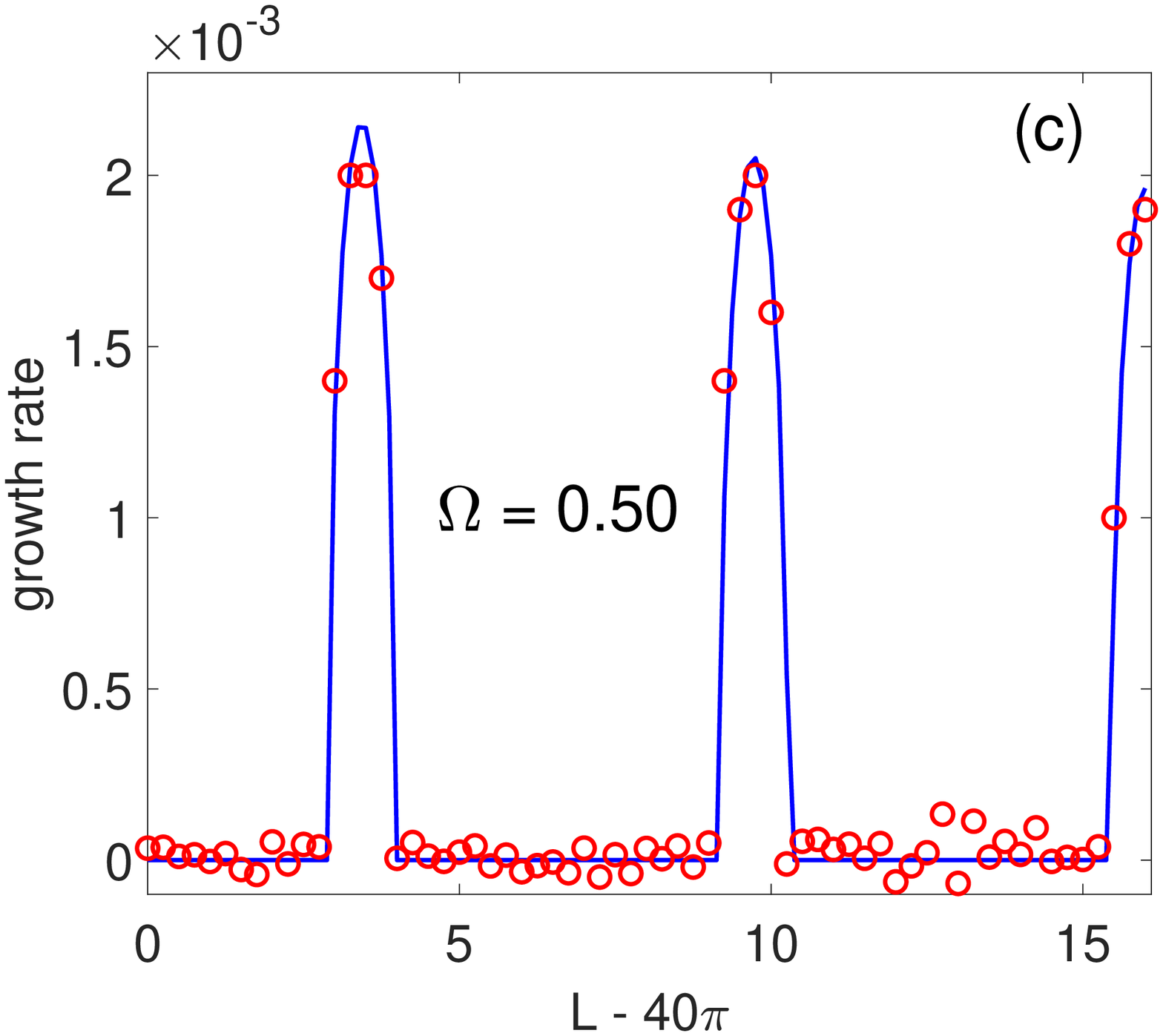}
\end{minipage}
\caption{Growth rate of the ``noise floor" 
as the function of length $L$ of the computational domain
for the same three values of $\Omega$ as in Fig.~\ref{fig_8}.
Note different vertical scales in the three panels.
Solid lines correspond to the analytical solution of \eqref{e6_10}, \eqref{e6_13}.
Circles are the result of SSM simulation, where we measured the amplitude
of the Fourier harmonic at $k=\kmax/2$ and followed the procedure described
at the end of Section 5. 
The simulation times for panels (b) and (c) are the same as those
in Figs.~\ref{fig_7}(b,c) and \ref{fig_8}(b,c): $t=1500$ and $t=5000$. 
For panel (a), the simulation time is $t=1000$. This larger time than in
Figs.~\ref{fig_7}(a) and \ref{fig_8}(a) had to be used to decrease the
effect of the transient behavior (see Figs.~\ref{fig_10}(b,c)) on the
computed growth rate (see text). Moreover, since for $\Omega=0.2$,
the NI at $k\approx \pm \kmax$ 
is so strong that it would destroy the numerical solution at $t=1000$
(see Fig.~\ref{fig_8}(a)), we had to filter out harmonics near the edges
of the computational spectral domain.
}
\label{fig_11}
\end{figure}

In Fig.~\ref{fig_11} we plot the growth rate of the ``noise floor" versus $L$.
The theoretical values are inferred from the spectral radius of $\bm{\Phi}(L)$,
which was found by the above analysis, via the relation (see \eqref{e6_11a}):
\be
\mbox{growth rate} \,=\, \big(\,\ln \rho(\bm{\Phi}(L))\,\big)/L\,.
\label{e6_13}
\ee
The agreement between our analysis and direct numerics is seen to be quite good. 

The numerical solution of system \eqref{e6_10} leaves it unclear why 
$\rho(\bm{\Phi}(L))$
depends on $L$. Moreover, such a dependence may even seem counter-intuitive given
that the amplitude of the perturbation \eqref{e6_01} changes only in the vicinity
of the soliton (which does not depend on $L$)
and remains intact in the rest of the computational domain. 
In Appendix D we show that while $\|\bm{\Phi}(L)\|$ does not depend on $L$,
$\rho(\bm{\Phi}(L))$ varies with $L$ periodically, the period being $2\pi/(2\Omega)$,
which is confirmed by Fig.~\ref{fig_11}. Note that these different dependences of
$\|\bm{\Phi}(L)\|$ and $\rho(\bm{\Phi}(L))$ on $L$ are consistent with the well-known
(see, e.g., \cite{Stewart_Matrix}) result: 
\be
\rho(\bm{\Phi}(L)) \le \|\bm{\Phi}(L)\| = \sigma_{\max} (\bm{\Phi}(L)) \,,
\label{e6_12}
\ee
where $\sigma_{\max}$ denotes the largest singular value.

To conclude this section, we note that while the ``noise floor" NI decreases
when $\Omega$ increases, it is still present even for a non-fragile soliton
 with $\Omega=0.75$. In that case, Eqs.~\eqref{e6_10} and \eqref{e6_13} predict
that the NI growth rate peaks to about $2.9\cdot 10^{-4}$,
i.e., slightly lower than the spectral edge NI, near $L=40\pi + 3.1$. 
Simulations by the SSM support this analytical result. Moreover, they 
reveal a feature of this NI that was not observed for $\Omega \le 0.5$. 
Namely, for certain $L$ values, only part of the ``noise floor" would become 
unstable. Analysis of this phenomenon would be more complicated than that based on
the $k$-independent ansatz \eqref{e6_05} and therefore is not considered here.

\section{Generalizations}


Here we will show that the two types of NI considered in Sections 5 and 6
occur in more general situations than in simulations of a single Gross--Neveu 
soliton. Thus, these types of NI 
appear to be engendered not by a specific model, its solution,
or even the numerical method, but by a combination of various factors. 
 Namely, we will first show that the same phenomena occur for more general
solutions of the same model. Second, we will show that one of them occurs
for the soliton of a different, well-known model in the relativistic field theory. 
Third, we will explain why, and show that, the same types of NI occur in 
other popular numerical methods applied to the Gross--Neveu model. 
Finally, we will show that by changing the 
boundary conditions of the numerical method, one can strongly diminish both types
of NI.

\subsection{More general solutions of the Gross--Neveu model} 

\begin{figure}[!ht]
\begin{minipage}{7.5cm}
\hspace*{-0.1cm} 
\includegraphics[height=5.6cm,width=7.5cm,angle=0]{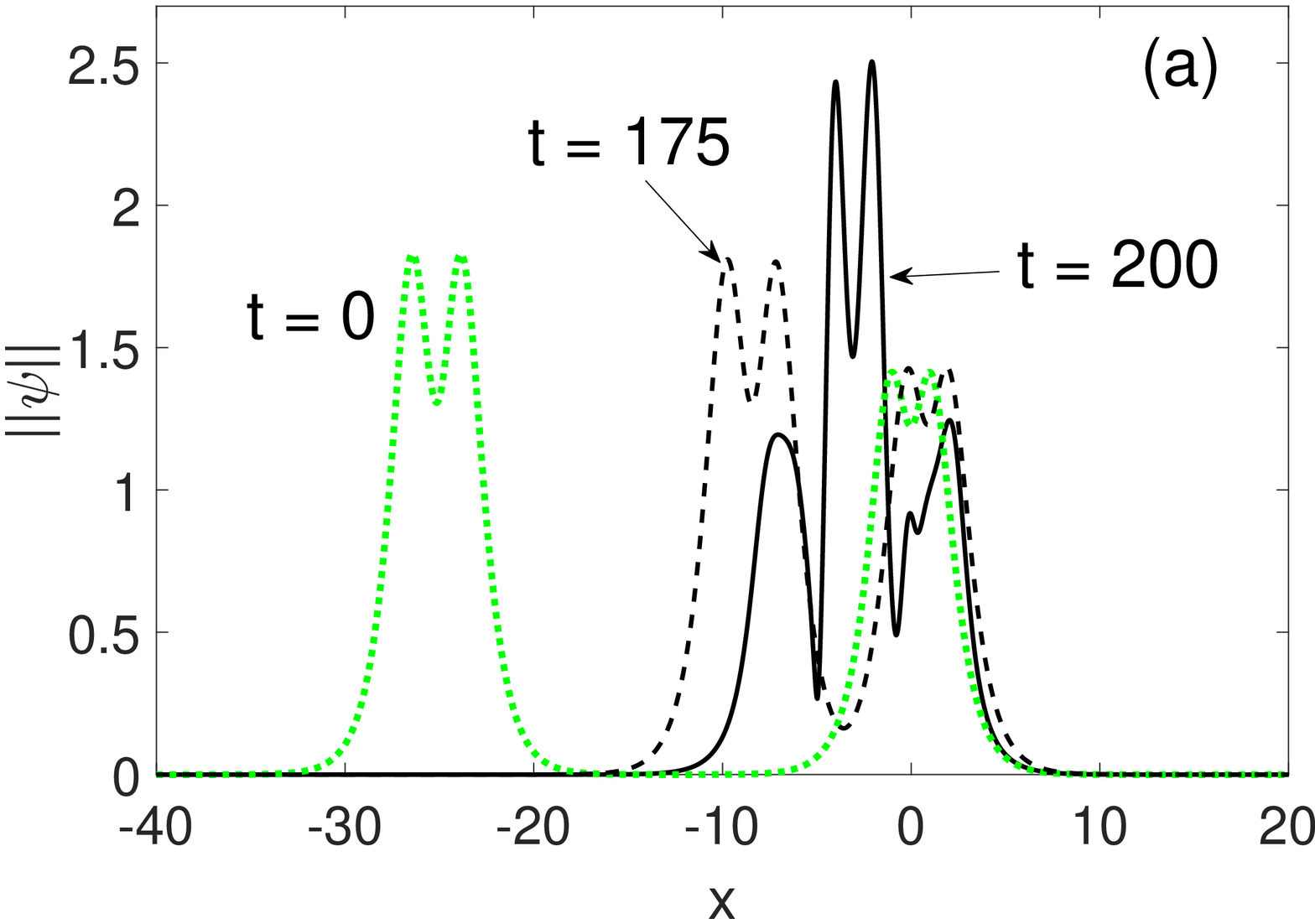}
\end{minipage}
\hspace{0.5cm}
\begin{minipage}{7.5cm}
\hspace*{-0.1cm} 
\includegraphics[height=5.6cm,width=7.5cm,angle=0]{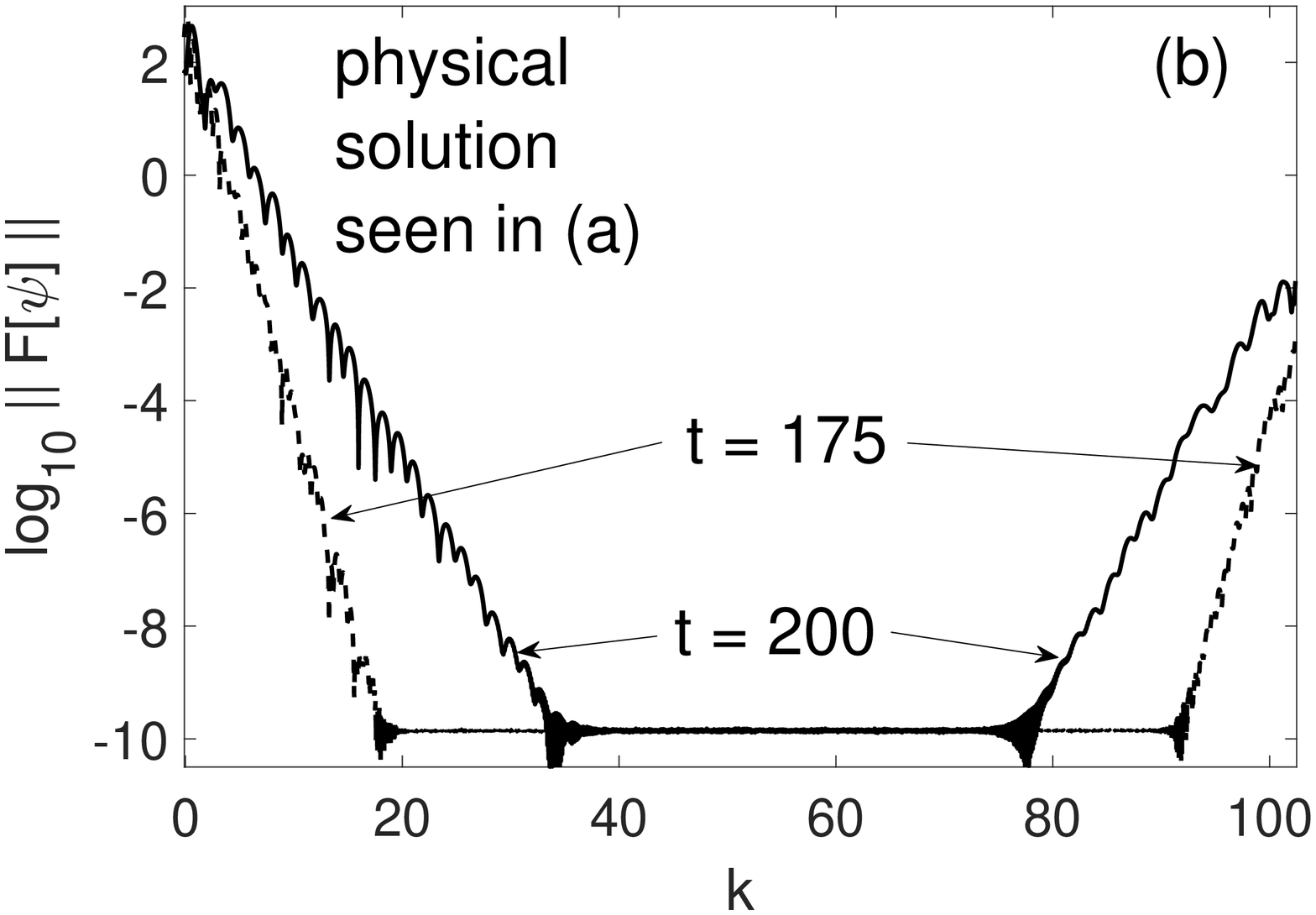}
\end{minipage}
\caption{Colliding solitons described in Section 7.1. Panel (a): Dotted (green), 
dashed, and  solid lines correspond to the solution at $t=0$, $175$, and $200$. 
Only part of the computational domain is shown for better visibility. \ Panel (b): 
Spectra (for $k\ge 0$) 
of the numerical solution at $t=175$ (dashed) and $t=200$ (solid). 
}
\label{fig_12}
\end{figure}

In Fig.~\ref{fig_12} we show the result of simulation of two colliding solitons
\eqref{e2_02}. The parameters of the solitons are: $\Omega_1=0.25$, $V_1=0$, 
$(x_0)_1=0$; \ $\Omega_2=0.15$, $V_2=0.1$, $(x_0)_2=-8\pi$. 
Other simulation parameters
are: $L=160\pi$, $N=2^{14}$, $\dt=\dx/5$, where the computational domain 
was chosen to be larger than that in the rest of this paper to minimize the
effect of radiation re-entering the domain and possibly corrupting the solution. 
The NI at the edges of the spectral domain is clearly visible. Moreover, 
one can see that the spectral support of unstable modes greatly expands 
between $t=175$ and $t=200$, which is when the moving soliton comes in close
proximity with the standing one.

We also considered the evolution of an initial pulse both of whose components
are $20$\% greater than 
those of
the standing soliton \eqref{e2_01} with $\Omega=0.2$. 
The NI at the edges of the computational domain destroys 
  the resulting near-soliton
 solution by
$t\gtrsim 300$. It should be noted that some other combinations of the components
of the initial pulse, --- e.g., where one is 20\% greater and the other is 20\% smaller 
than those of the soliton, --- do {\em not} lead to a strong NI. This occurs because
such an initial pulse evolves towards a soliton with a greater value of $\Omega$,
for which the NI is considerably weaker.

\subsection{Numerical instability for the massive Thirring soliton} 

The massive Thirring model in laboratory coordinates:
\be
u_t + u_x = i(v + u|v|^2), \qquad 
v_t - v_x = i(u + v|u|^2),
\label{e7_01}
\ee
has the one-soliton solution \citep[]{75_MTMsol_1, 75_MTMsol_2, 75_Exact} which, 
for zero velocity, can be written as: 
\be
U_{\rm sol} = \frac{\sin Q\,\exp[-it\cos Q]}{\cosh(x\sin Q - iQ/2)}\,, \qquad
V_{\rm sol} = \frac{-\sin Q\,\exp[-it\cos Q]}{\cosh(x\sin Q + iQ/2)}\,,
\qquad Q\in[0,\pi]\,.
\label{e7_02}
\ee
It may be noted that the Gross--Neveu and massive Thirring models both belong to 
the more general class of fermionic (nonlinear Dirac) field-theoretic models,
corresponding to the cases of scalar--scalar and vector--vector interactions,
respectively \citep[]{75_Exact}.
The massive Thirring model is integrable by the Inverse Scattering Transform 
\citep[]{77_KaupNewell, 77_KuzMikh}, a consequence of which is that the soliton
solution \eqref{e7_02} is physically stable \cite{96_KaupL}. 
Incidentally, a model similar to
\eqref{e7_01} occurs in a different field --- that of nonlinear light propagation 
in optical fibers with a periodic refractive index  
\citep[]{85_Winful}. However, its soliton (known as the Bragg, or gap, soliton)
is physically unstable in a certain range of its parameters \citep[]{98_Peli}.
Therefore, we chose to consider only the stable soliton \eqref{e7_02} of model
\eqref{e7_01} to avoid any issue of possible coexistence of physical and 
numerical instabilities.

To apply the SSM to \eqref{e7_01}, one first solves the linear part of those
equations in the Fourier domain:
\bsube
\be
\left( \ba{c} \widehat{u} \\ \widehat{v} \ea \right)_{\rm lin} = \frac1{1+\delta^2} 
 \left( \ba{cc} e^{i\gamma \dt} + \delta^2 \,e^{-i\gamma \dt} & 
                \delta\, \left( e^{i\gamma \dt} - e^{-i\gamma \dt} \right) \\
								\delta\, \left( e^{i\gamma \dt} - e^{-i\gamma \dt} \right) & 
								\delta^2 \,e^{i\gamma \dt} + e^{-i\gamma \dt} 
\ea \right) \, \left( \ba{c} \widehat{u} \\ \widehat{v} \ea \right)_n\,,
\label{e7_03a}
\ee
where $\gamma=\sqrt{k^2+1}$, $\delta=k+\gamma$. Then the nonlinear substep is:
\be
\left( \ba{c} \widehat{u} \\ \widehat{v} \ea \right)_{n+1} = 
\left( \ba{c} u_{\rm lin}\, \exp[i|v_{\rm lin}|^2\dt] \\
              v_{\rm lin}\, \exp[i|u_{\rm lin}|^2\dt] \ea \right)\,.
\label{e7_03b}
\ee
\label{e7_03}
\esube

Using the SSM \eqref{e7_03} to simulate \eqref{e7_01} with the initial condition
consisting of the soliton \eqref{e7_02} with $Q=0.35\pi$ and white noise of magnitude
on the order of $10^{-12}$, we have observed the numerically unstable modes at the
edges of the spectrum grew by 7 orders of magnitude in $t=1000$. The simulation
parameters were: $L=40\pi$, $N=2^{12}$, and $\dt=\dx/5$. 
The spectrum of the numerical solution looks qualitatively similar to that shown
in Fig.~\ref{fig_4}(a) and therefore is not shown here. 
For greater $Q$, this NI
developed even faster.\footnote{
   The limit $Q\To \pi$ for the massive Thirring soliton is known to 
	 have similarities with the 
	 limit $\Omega\To 0$ for the Gross--Neveu soliton.
	}
Conversely, for smaller $Q$, we found that this NI decreases. 
For example, for $Q=0.30\pi$ we found that the unstable modes at the spectral edges
grow at most by two orders of magnitude when we performed this simulation for a
variety of values of $L$.

We were unable to observe the ``noise floor" NI for the Massive
Thirring soliton in our numerics for any values of $Q$ and $L$
and initially were surprised. However, an analysis similar to
that presented in Section 6 revealed that for this model, this type of NI does not
occur. Below we present a summary of this analysis. 
Although it is possible to apply it directly to Eqs.~\eqref{e7_01}, it is more
convenient to cast those equation in a form with the l.h.s.~identical to that
of the Gross--Neveu model \eqref{e1add1_02} so as to follow the analysis of Section 6
as closely as possible. To that end, defining
\be
\psi_1=(u+v)/\sqrt{2}, \qquad \psi_2=(u-v)/\sqrt{2},
\label{e7_04} 
\ee
we transform \eqref{e7_01} into
%
\be
\psi_{1,t}+ \psi_{2,x} \,=\,  
            \frac{i}2\left( \psi_1|\psi_1|^2 - \psi_1^*\psi_2^2\right) + i\psi_1,
\qquad
\psi_{2,t}+ \psi_{1,x} \,=\,  
             \frac{i}2\left( \psi_2|\psi_2|^2 - \psi_2^*\psi_1^2\right)-i\psi_2,
\label{e7_05}
\ee
Linearizing \eqref{e7_05}, we obtain equations of the form \eqref{e2_04a} with
$\Omega$ on the l.h.s.~being replaced by $\cos Q$ and with:
\bsube
\be
P_0 = \frac12\left(|\Psi_1|^2 + |\Psi_2|^2\right), \qquad
P_1 \equiv 0, \qquad
Q_0 \equiv 0, \qquad
Q_1 \equiv 0, \qquad
\label{e7_06a}
\ee
\be
P_2 = -{\rm Im}\,\left( \Psi_1\Psi_2^* \right), \qquad
P_3 = \frac12\left(|\Psi_1|^2 - |\Psi_2|^2\right) \, +1\,,\qquad
Q_2 \equiv 0, \qquad 
Q_3 = \frac12\left(\Psi_1^2 - \Psi_2^2\right)\,,
\label{e7_06b}
\ee
\label{e7_06}
\esube
where $\Psi_{1,2}$ are the exact one-soliton solutions obtained from
\eqref{e7_02} and \eqref{e7_04}. 
In deriving $P_1\equiv 0$, we used the specific form  \eqref{e7_02} of the 
soliton.
The information about the possible growth of
the ``noise floor" perturbation is found from Eqs.~\eqref{e6_08a}, where
the matrix on the r.h.s.~is now:
\be
\bm{R}_{(\pm)} \,= \, i\sthr \left( \cos Q\szer + 
 \left( \ba{cc} P_0(t)  & 0 \\ 
                0 & P_0(t)  \ea \right)
								 \, \right) \,,
\label{e7_07}
\ee
where we have used that $P_0$ is a real-valued and even function. 
The key point to note is that due to the absence of off-diagonal terms in 
\eqref{e7_07}, the evolution of the perturbation $\bm{c}$ is unitary,
and hence there is no ``noise floor" NI in this case.

\subsection{Numerical instability of other methods applied to Gross--Neveu soliton} 

In our analysis of the NI in Sections 5 and 6 there was nothing that would
explicitly refer to the SSM as opposed to any other numerical method. 
Indeed, all we did was obtain a differential (in time) equation for the error
in the limit $\dt\To 0$. There was, however, an {\em implicit} assumption:
that the numerical scheme does not change the linear dispersion relation 
\eqref{e1_03} (for $|k|\gg 1$). Therefore, we expect that similar NI should
occur for any other numerical scheme that preserves the linear dispersion relation
of models \eqref{e1add1_02} and \eqref{e7_01}.

One such family of schemes is the Exponential Time Differencing
(ETD) and Integrating Factor methods (see, e.g., \cite{02_CoxMatt}).
It should be noted that for the Gross--Neveu model, such methods were
first proposed in \cite{10_IFNLDE} and have recently been considered in
\cite{16_SciChina}. We implemented the ETD method based on
the 4th-order explicit Runge--Kutta (RK) solver, referred to as ETD4RK in 
\cite{02_CoxMatt} and given by Eqs.~(26)--(29) in that paper. For the implementation
of the ETD methods, it is convenient to rewrite the simulated equation in a form
where the matrix multiplying the spatial derivative terms is diagonal. For 
Eqs.~\eqref{e1add1_02} this is achieved via the transformation inverse to \eqref{e7_04}, 
upon which they take on the  form:
\be
u_t = -u_x + i(u|v|^2 + u^* v^2 - v), \qquad 
v_t = \; v_x + i(v|u|^2 +v^* u^2 - u)\,.
\label{e7_08}
\ee
In simulating the Gross--Neveu equations in this form with the ETD4RK, we observed both
types of NI --- at the spectral edges and of the ``noise floor", --- with their
growth rates being practically the same (for the selected values of $L$ that we tested) to those
reported in Sections 5 and 6, as long as $\dt < \dx$ (see \eqref{e1_04}).

The other method for which we tested the presence of NI in the Gross--Neveu model
 is the pseudo-spectral 4th-order RK method. To implement it, 
one solves Eqs.~\eqref{e7_08} (or \eqref{e1add1_02}) by the classical RK method in time, 
with the
spatial derivatives being computed by the direct and inverse Fourier transform
\eqref{e4_02}. This method preserves the dispersion relation \eqref{e1_03} 
{\em only} in the limit $\kmax\dt\,\equiv\, (\pi/\dx)\dt \ll 1$, where the 4th-degree
polynomial in $k\dt$, which results from the 4th-order RK method, approximates
$\exp[ik\dt]$ sufficiently closely for all wavenumbers. However, recall that we are
interested in demonstrating the {\em unconditional} NI, which persists
in the limit $\dt\To 0$. Therefore, at least for sufficiently small $\dt$, the NI in
the pseudo-spectral method is expected to develop similarly to that in the SSM.
We confirmed this to indeed be the case. For example, for $\Omega=0.2$, the NI at the spectral
edge was suppressed by the numerical diffusion of the pseudo-spectral method for 
$\dt=\dx/5$; however, for $\dt=\dx/10$, the numerical diffusion became weak enough to
allow this NI to develop almost as fast as in the SSM. For $\Omega=0.1$, even
the relatively strong numerical diffusion was not able to prevent the NI at the spectral
edge from destroying the numerical solution around $t=200$. Also, the ``noise floor" NI
in the pseudo-spectral method was similar to that in the SSM for both
$\dt=\dx/5$ and $\dt=\dx/10$.

\subsection{Suppression of the numerical instability for the Gross--Neveu soliton}

The previous subsection illustrated the fact that the appearance of NI in simulations
of the Gross--Neveu soliton with sufficiently
small $\Omega$ occurs not just for the Fourier SSM, but
for a variety of numerical methods. In fact, NI for this problem was earlier reported in
\cite{14_NLDE_numerics, 15_NLDE_numerics}, although no details about its nature were
investigated in those studies. On the other hand, we showed in \cite{17_jaNLDE} that
merely imposing {\em nonreflecting} boundary conditions (BC):
\be
u(-L/2)=0, \qquad v(L/2)=0\,,
\label{e7_09}
\ee
   which is done by using the numerical method of characteristics,
 allows the small-$\Omega$ soliton to survive
over several thousands of time units. 
Imposing additional absorption at the boundaries was shown \cite{17_jaNLDE} to further
increase soliton's survival time, thus making it numerically stable even in ultra-long
simulations, in accordance with the theoretical prediction of \cite{12_BerkCom}.

Let us mention in passing that nonreflecting BC \eqref{e7_09} alone, i.e.,
without the additional absorber at the boundaries, do not entirely eliminate NI
for sufficiently small $\Omega$ (e.g., for $\Omega=0.1$). 
   Indeed, as shown in \cite{17_jaNLDE}, the key to suppress NI is to let
	 any radiation {\em completely} leave the computational domain without being partially
	 reflected inside it by the boundaries. (In the case of periodic BC, such a
	 reflection is replaced by mere re-entrance of the radiation into the domain.)
	 Nonreflecting BC \eqref{e7_09} allow the radiation to completely leave the
	 computational domain, without being reflected back (hence the name), {\em only}
	 in the absence of terms other than $u_x$, $v_x$
   on the r.h.s. of \eqref{e7_08}. 
   {\em With} the other terms, the amount of radiation reflected back into the domain
   is proportional to $1/k$. Hence higher harmonics are suppressed more than the
   lower ones; yet it is the lower harmonics that appear to be ``responsible" for
	 soliton's ``fragility" (see Section 2). 
	 This is why the additional absorber, that would equally absorb {\em all} harmonics,
	 was needed in \cite{17_jaNLDE} to avert destruction of the soliton.


\section{Conclusions and Discussion}

In this work, we showed that 
the (Fourier) SSM, \eqref{e3add1_01} or  \eqref{e4_01},
 for the Gross--Neveu model \eqref{e1add1_02}
may exhibit NI, and analytically studied three distinct mechanisms that can lead to it. 
Two of these mechanisms lead to {\em unconditional} NI, which, to our knowledge,
have never been analyzed previously.

The first type of NI,  analyzed in Section 4, 
{\em may} occur when the time step
exceeds the ``threshold" \eqref{e1_04} 
set by the Courant--Friedrichs--Lewy condition. However, unlike in other
schemes for hyperbolic equations, this high-$k$ NI is observed {\em only}
if the simulated solution of the model also exhibits low-$k$ instability. 
In practice, such conditional NI (unlike that for the NLS!) 
is inconsequential for the outcome of 
the simulations. Indeed, if a low-$k$ instability is present, it will destroy 
the solution long before this high-$k$ NI will. 
On the other hand, if there is no (or too weak) low-$k$ instability,
then harmonics near the ``resonance" wavenumber $k_{\pi}$ will also be stable,
even when the time step exceeds the ``threshold" \eqref{e1_04}.

The NI of the second type, analyzed in Section 5, 
occurs near the edges of the computational spectrum. 
It is {\em unconditional}, i.e.~persists for arbitrarily small $\dt$. This 
type of NI becomes stronger as $\Omega$ of the soliton decreases, i.e.,
as the soliton becomes more ``fragile". While the corresponding growth rate
depends on the length $L$ of the computational domain, 
this dependence diminishes as $\Omega$ decreases: 
see Figs.~\ref{fig_8} and \ref{fig_9}(b). 
There is no ``simple" qualitative reason that would unambiguously explain the
origin of this NI. 
Analysis of this NI required numerical solution of a relatively large
eigenvalue problem, but this is still several orders of magnitude faster than 
direct numerical simulations, especially for non-fragile solitons.

The NI of the third type, analyzed in Section 6, 
occurs for Fourier harmonics of the ``noise floor". This NI is also unconditional,
and it also becomes stronger as $\Omega$ of the soliton decreases. The growth rate
is essentially periodic in $L$, but also decreases in inverse proportion to it (i.e.,
for a sufficiently large $L$, slowly). Unlike for the second type of NI, there {\em is}
a qualitative explanation for the third one. Namely, it can occur as perturbations
 travelling in the opposite directions interact in the vicinity of the soliton
via the ``potential" created by it.  
This process can be amplified when the perturbations do so repeatedly, which is enabled
by their staying in the computational domain due to periodic boundary conditions. 
Analysis of this NI requires numerical solution of only two
coupled ordinary differential equations; see Eqs.~\eqref{e6_08} and \eqref{e6_10}.
Remarkably, the same analysis was able to explain the {\em absence} of the ``noise floor" NI
in the massive Thirring model (Section 7.2).

In Section 7, we demonstrated that these two types of unconditional NI 
can occur in more general situations.
First, they occur in multi-soliton solutions, as long as some of the consitituent 
solitons have a sufficiently small $\Omega$. Second, they can occur for other
models that involve solitons in asymptotically dispersionless coupled-mode equations; 
an example is the massive Thirring soliton. Third, they can occur for methods other than the 
SSM.

Based on these generalizations, we propose that there may be only two essential 
conditions that need to be met for these unconditional NIs to be observed. 
The first is 
that the simulated soliton's parameters must be in a
certain range (e.g., a sufficiently small $\Omega$ for the Gross--Neveu soliton or a
sufficiently large $Q$ for the massive Thirring soliton). Although these NIs 
 were found even for Gross--Neveu solitons with $\Omega$ as large as $0.75$,\footnote{
   and we have no reason to think that they would not occur for even greater $\Omega$
	}
their growth rates were too small to be observed in any but the ultra-long simulations.

The second essential condition for observing these NIs is the 
boundary conditions that permit a substantial part of the radiation 
to re-enter the computational domain. The simplest such BC
are periodic; they are automatically imposed when the numerical method involves
the discrete Fourier transform. However, there are indications in \cite{15_NLDE_numerics}
 that other BC,
such as homogeneous Dirichlet, may also lead to similar NIs. As we showed
in \cite{17_jaNLDE} and stressed in Section 7.4, only a combination of nonreflecting
BC and an additional absorber (or a more sophisticated technique \cite{19_DiracPML})
could suppress NI for solitons with arbitrarily small $\Omega$.

We now relate our results with those of \cite{14_NLDE_numerics, 15_NLDE_numerics},
where numerical instabilities of the Gross--Neveu soliton were reported. 
The authors of \cite{14_NLDE_numerics} used a 4th-order non-Fourier SSM, where
the linear substep was computed by the method of characteristics, subject to 
nonreflecting BC \eqref{e7_09}. The evidence of high-$k$ NI is reported there 
in the captions to Figs.~5 and 6. This may be surprising given that we argued
in Section 7.4 that the nonreflecting BC tend to suppress high-$k$ NI
(and so a low-$k$ NI would be observed before any high-$k$ could be seen). 
A detailed analysis of this conundrum would require a separate study and is clearly 
outside the scope of this one. Below we present only a plausible resolution,
which consists of two ingredients. First, 
it should be noted that the 4th-order method of \cite{14_NLDE_numerics} required
a relation $\dt = 12\dx$. Thus, the stability ``threshold" \eqref{e1_04} was
exceeded significantly. Had the method been 2nd-order, as here, this would not
have caused a high-$k$ NI to destroy the soliton. 
However, it was shown in \cite{12_ja} that the equation satisfied by a high-$k$
numerical error in a 4th-order method is different from that for the 2nd-order
one. 
This fact makes it possible for the
modes near the ``resonant" wavenumbers
$k_{n\pi}$, $n=1,\ldots,12$  to become unstable regardless of any low-$k$ NI.
In fact,  the numerical solution reported in \cite{14_NLDE_numerics} for the
soliton with $\Omega=0.1$ at $t=100$, was provided to this author in a private
communication by the first co-author of \cite{14_NLDE_numerics}, and it
{\em confirmed} the above hypothesis about an NI developing near the
``resonant" wavenumbers $k_{n\pi}$, $n=1,\ldots,12$.

The authors of \cite{15_NLDE_numerics} employed three methods to simulate the
dynamics of the soliton. Their best-performing method is, in fact, the pseudo-spectral
method considered in Section 7.3. As we showed there, this method is subject to
the same NIs as the SSM. Moreover, the NIs would become more pronounced as one 
decreases $\dt/\dx$, because it is then suppressed less by numerical diffusion.
In  \cite{15_NLDE_numerics}, the authors conspicuously stated that the NIs 
should go away as $L\To \infty$ and $\dx \To 0$. However, our analysis indicates 
that both of these statements are incorrect. Namely, Figs.~\ref{fig_8}(a) and
 \ref{fig_9}(b) show
that the growth rate of the spectral edge NI for sufficiently small $\Omega$ 
persists almost unchanged from $L=40\pi$ to $L=1600\pi\approx 5,000$, which is more than
30 times the largest length considered in \cite{15_NLDE_numerics}. We also verified that
for the same length as used in most of the simulations in this work, $L=40\pi$,
and for $N=2^{17}$ ($\dx\approx 10^{-3}$, i.e. 32 times smaller than in the rest
of the work) the growth rate of the $\Omega=0.2$ soliton was the same as for 
$N=2^{12}$ ($\dx \approx 0.03$). Thus, the spectral edge NI is not perceptibly
affected by either the $L\To\infty$ or $\dx\To 0$ limits. 
The ``noise floor" NI is not affected by the $\dx\To 0$ limit, but its growth rate
does indeed vanish as $L\To \infty$ since it scales as $1/L$.

Finally, we note that results of our study should alert researchers who 
numerically study new models
similar to the Gross--Neveu model (see, e.g.,  \cite{17_coupledDirac, 19_IgorPTDirac}) 
about possible NI that may occur in simulations of those models.
Our analysis may also serve as a prototype for that of 
the Zakharov equations, describing interaction of short and long waves \cite{72_Zakharov}. 
They consist of an NLS coupled to the hyperbolic nonlinear wave equation. 
The dispersion relation of the latter equation,
which for high wavenumbers coincides with \eqref{e1_03}, may potentially
enable an NI similar to those considered in this work. In fact, an indication
of a possible NI in simulations of the Zakharov equations by an SSM 
is reported in Fig.~10 of \cite{03_ZakharovEqs}.


\section*{Acknowledgments}

This work was supported in part by the NSF grant
 DMS-1217006. The author acknowledges a useful conversation with J.~Yang
 and thanks S.~Shao for providing the numerical results mentioned three
 paragraphs above.


\renewcommand{\theequation}{A.\arabic{equation}}
\section*{Appendix A: \ Code of 2nd-order SSM for Gross--Neveu soliton}
 \setcounter{equation}{0}

\small

\begin{lstlisting}[
  style      = Matlab-editor,
  basicstyle = \mlttfamily,
]
clear all; 
% --- Simulation parameters:  
L = 40*pi;  N = 2^12;  dx = L/N;  x = [ -L/2 :dx: L/2-dx ]; 
dk = 2*pi/L;  k = [ 0 : N/2-1   -N/2 : -1 ]*dk; 
dt = dx/2;  tmax = 1500;  tplot = 50; 
% --- Initial conditions:       
Omega = 0.35; beta = sqrt( 1 - Omega^2 );  mu = (1 - Omega)/(1 + Omega);     
u0 =    sqrt(2*(1 - Omega))./( (1 - mu*tanh(beta*x).^2).*cosh(beta*x) );
v0 = 1i*u0.*sqrt(mu).*tanh(beta*x);
u = u0 + 10^(-12)*randn(size(x));  v = v0 + 10^(-12)*randn(size(x)); 
% --- Auxiliary parameters to be used at each step in the loop:
idt = 1i*dt;  ckdt = cos(k*dt);  mi_skdt = -1i*sin(k*dt);  
ckdto2 = cos(k*dt/2);  mi_skdto2 = -1i*sin(k*dt/2);
% --- Main calculation:  
  for nn = 1 : round(tmax/dt)  
      fftu = fft(u);  fftv = fft(v);
      if nn == 1
          u = ifft( ckdto2.*fftu + mi_skdto2.*fftv );  
          v = ifft( ckdto2.*fftv + mi_skdto2.*fftu );
      else
          u = ifft( ckdt.*fftu + mi_skdt.*fftv );  
          v = ifft( ckdt.*fftv + mi_skdt.*fftu );
      end
      expaux = exp( idt*( (abs(u).^2-abs(v).^2) - 1) );
      u = u.*expaux;  v = v./expaux;
      if nn == round(tmax/dt)
          u = ifft( ckdto2.*fftu + mi_skdto2.*fftv );  
          v = ifft( ckdto2.*fftv + mi_skdto2.*fftu );
      end
      % --- Plot the results:  
      if rem(nn, round(tplot/dt)) == 0
          figure(123);
          subplot(211); plot(x, abs(u), 'b', x, abs(v), 'r');  
          xlabel('x'); ylabel('{|u|, |v|}');
          subplot(212); plot(k, log10(abs(fftu)), 'b', k, log10(abs(fftv)), 'r'); 
          xlabel('k, wavenumber'); ylabel('log_{10}(|F[u, v](k)|)');
          title(['t=', num2str(dt*nn)], 'fontsize', 13); pause(0.3) 	    
      end
  end
  

\end{lstlisting}

\normalsize


\renewcommand{\theequation}{B.\arabic{equation}}
\section*{Appendix B: \ Weak fragility of solitons with $\Omega > 0.6$}
 \setcounter{equation}{0}

Simulations reported in Section 3.1 illustrated our general observation
that if the soliton is fragile, then exponentially growing modes 
will appear around $k_{n\pi}$-peaks whenever
$\dt > \dx$. One can conjecture that the converse statement
may also be true, namely: If one detects such
 modes around $k_{n\pi}$-peaks, then the soliton must be (weakly) fragile, even 
though it has not yet exhibited any fragile behavior for the same 
simulation time. In other words, can a high-$k$ instability
signal the occurrence of a low-$k$ one?

\begin{figure}[!ht]
\begin{minipage}{5.2cm}
\hspace*{-0.1cm} 
\includegraphics[height=4.2cm,width=5.1cm,angle=0]{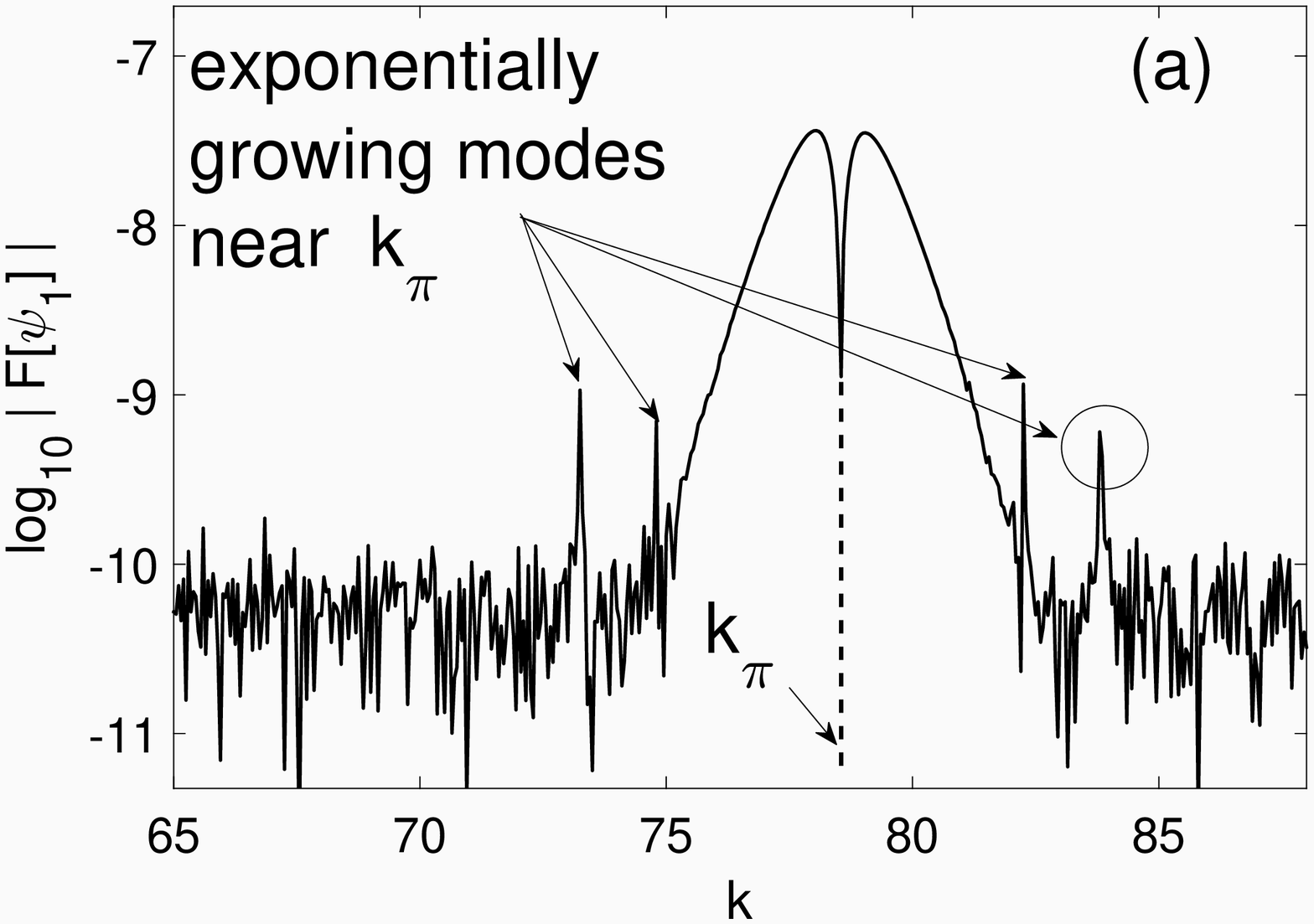}
\end{minipage}
\hspace{0.2cm}
\begin{minipage}{5.2cm}
\hspace*{-0.1cm} 
\includegraphics[height=4.2cm,width=5.1cm,angle=0]{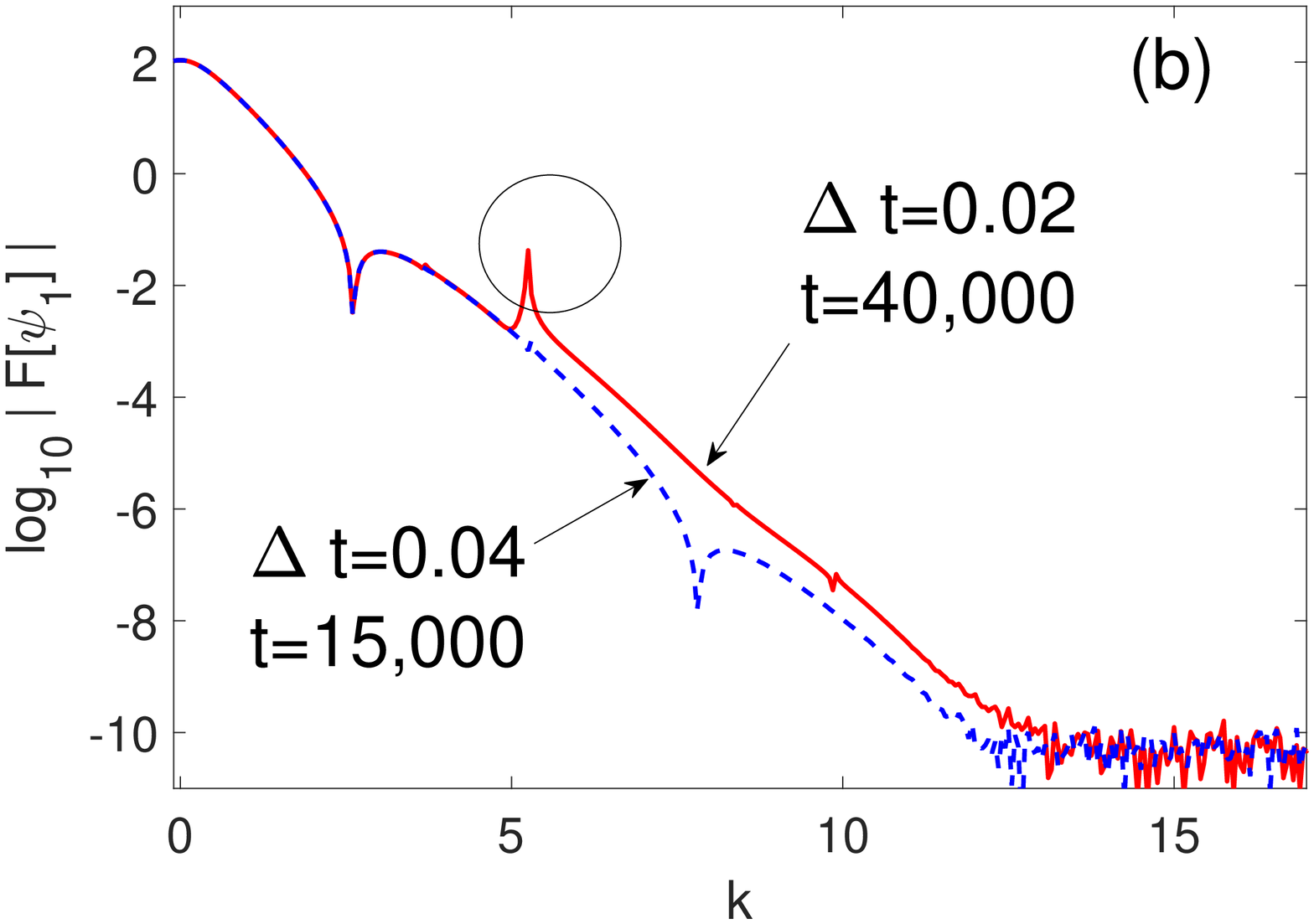}
\end{minipage}
\hspace{0.2cm} 
\begin{minipage}{5.2cm}
\hspace*{-0.1cm} 
\includegraphics[height=4.2cm,width=5.1cm,angle=0]{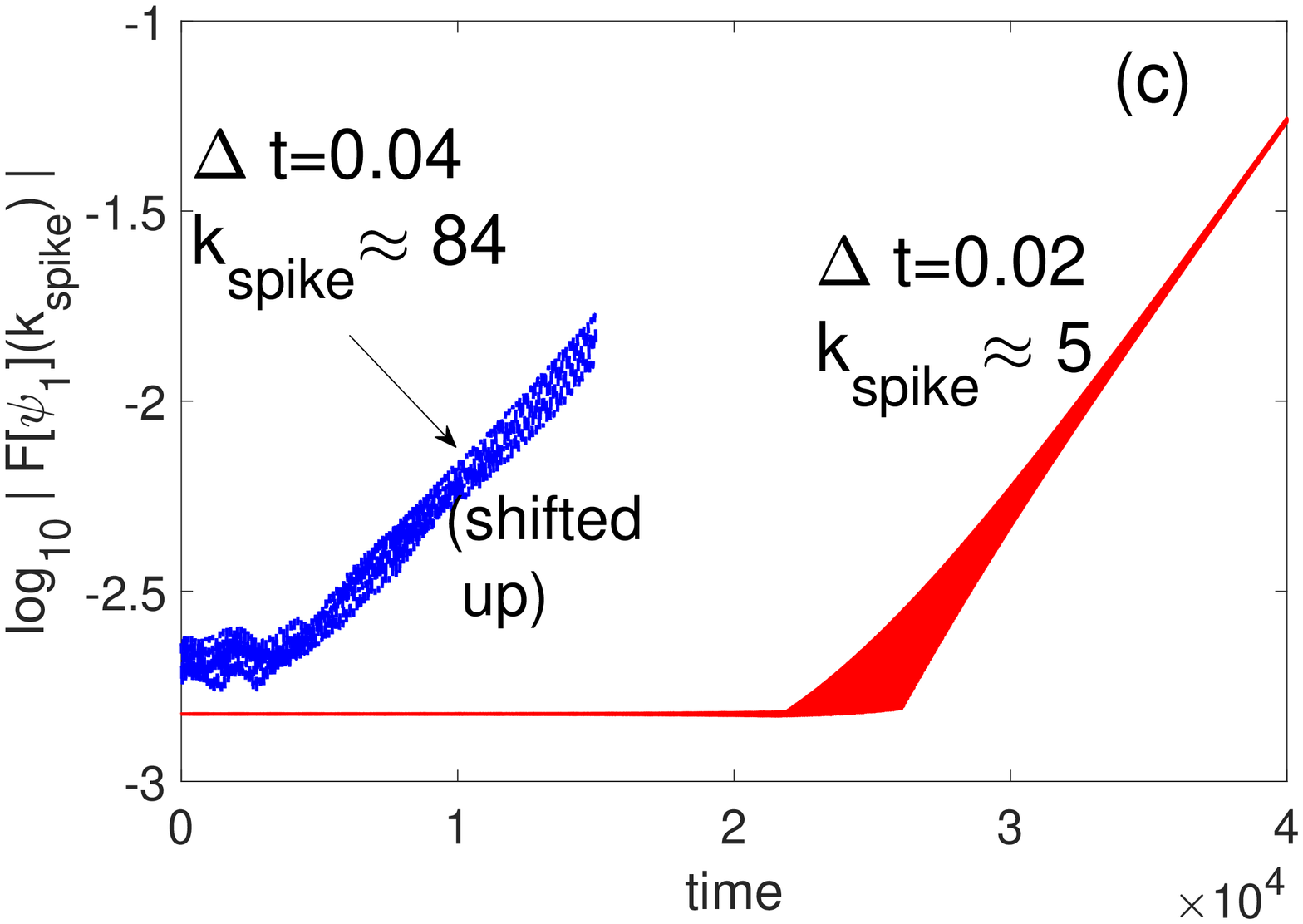}
\end{minipage}
\caption{Ultra-long-time simulations of $\Omega=0.75$-soliton (non-fragile). 
Simulation parameters are: $L=40\pi$, $N=2^{12}$ ($\dx=\dtthresh\approx 0.031$).
Results for the $\psi_2$-component
 are similar to those for the $\psi_1$-component and hence are not shown.
\ (a) Close-up on the vicinity of the $k_{\pi}$-peak for $\dt=0.04 
 \approx 1.3\dx$ and $t=15,000$. \ 
(b) Solitonic part of the spectrum of the numerical solution obtained
at different times and with different $\dt$. The lines are indistinguishable 
for $k<5$. \ 
(c) Evolution of the amplitudes of the ``spikes" circled in panels (a) and (b).
The left curve is shifted up by 7.2 units for better visibility.  
}
\label{figA_1}
\end{figure}

In Fig.~\ref{figA_1} we demonstrate that this is indeed the case. 
In panel (a) we show a close-up on the vicinity of the $k_{\pi}$-peak
for the $\Omega=0.75$-soliton at $t=15,000$, where exponentially
growing ``spikes" near the wider $k_{\pi}$-peak are clearly discernible.
For the same $t$, there is no sign  of fragile behavior of this numerical
solution in either the Fourier domain (panel (b)) or $x$-domain (not shown). 
However, at a much greater time, a small ``spike" becomes visible
within the soliton's own spectrum (around $k=5$); see Fig.~\ref{figA_1}(b). 
Finally, Fig.~\ref{figA_1}(c) shows that one can detect the
exponentially growing modes around the 
$k_{\pi}$-peaks much earlier than similarly growing modes around the soliton.

Two clarifications regarding the above result are in order. 
First, the reason that exponentially growing modes are observed sooner
around the $k_{\pi}$-peak is that the ``noise floor", from which these modes
arise, is several orders of magnitude closer to the $k_{\pi}$-peak than to
the soliton's spectral maximum. Thus, it takes less time for such modes to
become visible relative to the $k_{\pi}$-peak than relative to the soliton. 
Second, the reason why the emergence of the exponential instability is
significantly delayed from the start of the simulations 
(see Fig.~\ref{figA_1}(c))
was explained in \cite{16_ja2}. It is related to
the fact that both (i) the overlap of the actual unstable mode with any one 
Fourier harmonic and (ii) the instability growth rate, are small.


\renewcommand{\theequation}{C.\arabic{equation}}
\section*{Appendix C: \ Reduced form of Eqs.~\eqref{e5_04}, and 
matrix blocks in \eqref{e5_07c}}
 \setcounter{equation}{0}

Substitution of \eqref{e5_05} and \eqref{e5_06} into Eqs.~\eqref{e5_04}
and taking complex conjugate of two of the four equations 
yields the following system: 
\bsube
\be
{\widehat{a}_{(+),\,t}} \,= \,\hspace{12.5cm}
 \label{B_01a}
\ee
$$
 i \left( \Omega - \dkm \right) \wap
 + i\, {\mathcal F} \left[ 
	       \left[ P_0\, \pap  + (iP_2 +P_3) \pbm + 
	        Q_3\, \pam^{\,*} + (Q_0+Q_1) \pbp^{\,*} 
					 \right]^{(<0)} 		\, \right] \,, 
$$
\be
{\widehat{b}_{(-),\,t}} \,= \,\hspace{12.5cm}
 \label{B_01b}
\ee
$$
 i \left( \Omega + \dkp \right) \wbm
 + i\, {\mathcal F} \left[ 
	       \left[ (-iP_2+P_3) \pap  + P_0\, \pbm + 
	        (Q_0 - Q_1) \pam^{\,*} + Q_3\, \pbp^{\,*}
					 \right]^{(\ge 0)} 		\, \right] \,.
$$
\be
{\widehat{a}_{(-),\,t}}^{\,*} \,=\, \hspace{12.5cm}
 \label{B_01c}
\ee
$$
-i \left( \Omega + \dkm \right) \wam^{\,*}
 - i\, \left( {\mathcal F} \left[ 
	       \left[ ( Q_3\, \pap^{\,*}  + (Q_0-Q_1) \pbm^{\,*} + 
	        P_0\, \pam + (-iP_2+P_3) \pbp 
					 \right]^{(<0)} 		\, \right] \,\right)^* , 
$$
\be
{\widehat{b}_{(+),\,t}}^{\,*} \,= \,\hspace{12.5cm}
 \label{B_01d}
\ee
$$
 -i \left( \Omega - \dkp \right) \wbp^{\,*}
 - i\, \left( {\mathcal F} \left[ 
	       \left[ (Q_0+Q_1) \pap^{\,*}  + Q_3\, \pbm^{\,*} + 
	        (iP_2+P_3) \pam + P_0\, \pbp
					 \right]^{(\ge 0)} 		\, \right] \,\right)^* .
$$
\label{B_01}
\esube
Here we have also used that $P_1=Q_2=0$ from \eqref{e2_04}.

The $M\times M$ 
matrix blocks $\mathbb{C}_{1\,m}$, $m=1,\ldots,4$,
in \eqref{e5_07c}  are obtained from
the respective four convolution-like terms in \eqref{B_01a}. Similarly,
the other blocks are obtained from the respective terms in 
\eqref{B_01b}--\eqref{B_01d}. Here we will present a derivation of 
$\mathbb{C}_{11}$ and $\mathbb{C}_{12}$ and will state the results for
the other $\mathbb{C}_{jm}$, which are derived analogously. 

To obtain the form of $\mathbb{C}_{11}$, we substitute the first relation
from \eqref{e5_01} into the first ${\mathcal F}$-term in \eqref{B_01a}
and use the definitions of discrete Fourier transform and its inverse
\eqref{e4_02} to obtain the following expression for the harmonic
with wavenumber $k_{-j}=-j\dk$, where $1\le j \le M$:
\bea
{\mathcal F} \left[ P_0\, \pap  \, \right]_{-j}   & = &
\sum_{m=-N/2}^{N/2-1} e^{-ik_{-j} x_m} 
\frac1N \sum_{n=-N/2}^{N/2-1} \widehat{P_0}_{\,n}\; e^{ik_n x_m} 
\frac1N \sum_{l=1}^{M} \widehat{a}_{(+)\,l} \; e^{-ik_l x_m} 
\label{C_01}  \\
 & = & \frac1{N^2} \sum_{n=-N/2}^{N/2-1} \widehat{P_0}_{\,n}\; 
       \sum_{l=1}^{M} \widehat{a}_{(+)\,l} \;  \delta_{n,\;l-j} 
\; = \; \frac1N \sum_{l=1}^{M} \widehat{P_0}_{\;l-j}\;\widehat{a}_{(+)\,l};
 \nonumber
\eea
in the last line, $\delta$ is the Kroneker symbol. Therefore,
\bsube
\be
\mathbb{C}_{11} =  \left( 
 \ba{cccc} 
 \widehat{P_0}_{\,0} & \widehat{P_0}_{\,1} & \cdots & \widehat{P_0}_{\,M-1} \\
 \widehat{P_0}_{\,-1} & \widehat{P_0}_{\,0} & \cdots & \widehat{P_0}_{\,M-2} \\
 \ddots  & \ddots & \ddots & \ddots \\
 \widehat{P_0}_{\,-(M-1)} & \widehat{P_0}_{\,-(M-2)} & \cdots & \widehat{P_0}_{\,0}
 \ea \right).
\label{C_02a}
\ee
This is a Toeplitz matrix, and half of the entries in \eqref{e5_07c}
will also be Toeplitz. Therefore, we introduce a notation:
\be
\mathbb{C}_{11} \,\equiv  \, {\mathcal T}\left[ 
 P_0,\, -(M-1),\, M-1\,\right],
\label{C_02b}
\ee
\label{C_02}
\esube
where: \ ${\mathcal T}$ stands for a Toeplitz matrix, the first argument
indicates the function whose harmonics make up the entries of the matrix,
and the third and fourth entries indicate the harmonic's indices of the
lower-left and upper-right entries of the matrix, respectively. 

Similarly, and denoting $iP_2+P_3\equiv P_{i23}$ for brevity, one has:
\bea
{\mathcal F} \left[ P_{i23}\, \pbm  \, \right]_{-j}   & = &
\sum_{m=-N/2}^{N/2-1} e^{-ik_{-j} x_m} 
\frac1N \sum_{n=-N/2}^{N/2-1} \widehat{P_{i23}}_{\,n}\; e^{ik_n x_m} 
\frac1N e^{-i\dk x_m} \, 
        \sum_{l=1}^{M} \widehat{b}_{(-)\,l} \; e^{ik_l x_m} 
\label{C_03}  \\
  & = & \frac1N 
   \sum_{l=1}^{M} \widehat{P_{i23}}_{\;1-l-j}\;\widehat{b}_{(-)\,l};
 \nonumber
\eea
whence
\bsube
\be
\mathbb{C}_{12} \,=\,  \left( 
 \ba{cccc} 
 \widehat{P_{i23}}_{\,-1} & \widehat{P_{i23}}_{\,-2} 
                         & \cdots & \widehat{P_{i23}}_{\,-M} \\
 \widehat{P_{i23}}_{\,-2} & \widehat{P_{i23}}_{\,-3} 
                         & \cdots & \widehat{P_{i23}}_{\,-M-1} \\
 \ddots  & \ddots & \ddots & \ddots \\
 \widehat{P_{i23}}_{\,-M} & \widehat{P_{i23}}_{\,-M-1} 
                         & \cdots & \widehat{P_{i23}}_{\,-2M+1}
 \ea \right).
\label{C_04a}
\ee
This is a Hankel (i.e. an ``upside-down" Toeplitz) matrix, and the
remaining half of the entries in \eqref{e5_07c} are
also Hankel. Therefore, we introduce a notation:
\be
\mathbb{C}_{12} \, \equiv  \, {\mathcal H}\left[ 
 iP_2+P_3,\, -1,\, -2M+1\,\right],
\label{C_04b}
\ee
\label{C_04}
\esube
where, similarly to \eqref{C_02b}: \ 
${\mathcal H}$ stands for a Hankel matrix, the first argument
indicates the function whose harmonics make up the entries of the matrix,
and the third and fourth entries indicate the harmonic's indices of the
upper-left and lower-right entries of the matrix, respectively. 

Similarly, one has:
\bsube
\be
\mathbb{C}_{13} \, =  \, {\mathcal H}\left[ 
            Q_3,\, -2,\, -2M\,\right],
\qquad
\mathbb{C}_{14} \, =  \, {\mathcal T}\left[ 
            Q_0+Q_1,\, -M,\, M-2\,\right];
\label{C_05a}
\ee
\be
\ba{ll}
\mathbb{C}_{21} \, =  \, {\mathcal H}\left[ 
            -iP_2+P_3, \, 1, 2M-1\, \right],
	& 
\mathbb{C}_{22} \, =  \, {\mathcal T}\left[ 
            P_0,\, M-1,\, -(M-1)\,\right],
		\vspace{0.2cm}  \\
\mathbb{C}_{23} \, =  \, {\mathcal T}\left[ 
            Q_0-Q_1,\, M-2,\, -M\,\right],
  & 
\mathbb{C}_{24} \, =  \, {\mathcal H}\left[ 
            Q_3, \, 0, 2M-2\, \right];
\ea
\label{C_05b}
\ee
\be
\ba{ll}
\mathbb{C}_{31} \, =  \, -{\mathcal H}\left[ 
            Q_3^*,\, -2,\, -2M\,\right], 
						&
\mathbb{C}_{32} \, =  \, -{\mathcal T}\left[ 
            (Q_0-Q_1)^*,\, -M,\, M-2\,\right],
		\vspace{0.2cm}  \\
\mathbb{C}_{33} \,=  \, -{\mathcal T}\left[ 
 P_0^*,\, -(M-1),\, M-1\,\right],
			&
\mathbb{C}_{34} \, =  \, -{\mathcal H}\left[ 
 (-iP_2+P_3)^*,\, -1,\, -2M+1\,\right];
\ea
\label{C_05c}
\ee
\be
\ba{ll}
\mathbb{C}_{41} \, =  \, -{\mathcal T}\left[ 
            (Q_0+Q_1)^*,\, M-2,\, -M\,\right],
  & 
\mathbb{C}_{42} \, =  \, -{\mathcal H}\left[ 
            Q_3^*, \, 0, 2M-2\, \right],
		\vspace{0.2cm}  \\						
\mathbb{C}_{43} \, =  \, -{\mathcal H}\left[ 
            (iP_2+P_3)^*, \, 1, 2M-1\, \right],
	& 
\mathbb{C}_{44} \, =  \, -{\mathcal T}\left[ 
            P_0^*,\, M-1,\, -(M-1)\,\right].					
\ea
\label{C_05d}
\ee
\label{C_05}
\esube

These matrix blocks can also be easily coded. For example, 
if \verb+dftP0+ and \verb+dftPi23+ denote the discrete Fourier spectra
\eqref{e4_02} with harmonics limited to \ $l\in [-2M,\,2M]$ 
of $P_0$ and $(iP_2+P_3)$, then the 
respective matrix blocks \eqref{C_02} and \eqref{C_04} 
can be programmed in Matlab as:
\begin{verbatim}
for j = 1 : M
   C_11(j, :) =   dftP0(2*M+1 - (j-1) : 2*M+1 - (j-1) + (M-1));
   C_12(j, :) = dftPi23(2*M - (j-1) :-1:  2*M - (j-1) - (M-1));
end
\end{verbatim}


\renewcommand{\theequation}{D.\arabic{equation}}
\section*{Appendix D: \ Explanation of why and how
the spectral range of ${\bm \Phi}(L)$  depends on $L$}
 \setcounter{equation}{0}

For the purpose of this explanation, it will suffice to replace the soliton
with a box-like profile as follows. Since $\Psi_1(x)$ has a bell-like shape, 
we replace it by some
constant for $x\in[-L_{\rm sol}/2,\,L_{\rm sol}/2]$. On the other hand,
$i\Psi_2(x)$ changes sign at $x=0$ (see Fig.~\ref{fig_1}), and therefore it
needs to be replaced by a piecewise-constant profile on the same interval.
Thus, we approximate:
\bsube
\be
\Psi_1(x) = \left\{ \ba{ll} A, & x\in[-L_{\rm sol}/2,\,L_{\rm sol}/2] \\ 
                            0, & |x| > L_{\rm sol}/2
									  \ea \right.\,,  
	\qquad
\Psi_2(x) = \left\{ \ba{ll} \pm iB, & x\in[0,\,\pm L_{\rm sol}/2] \\ 
                            0, & |x| > L_{\rm sol}/2
									  \ea \right.\,,
\label{E_01a}
\ee
for some $A,B,L_{\rm sol}>0$, with $L_{\rm sol}\ll L$. 
Then the $t$-dependent matrix on the r.h.s.~of \eqref{e6_10}
is replaced with:
\be
\left( \ba{cc} \Psi_1^2-\Psi_2^2 & (\Psi_1+\Psi_2)^2 \\ 
               (\Psi_1-\Psi_2)^2 & \Psi_1^2-\Psi_2^2 \ea \right) \;\equiv\; 
\left( \ba{cc} C & C\,e^{\pm i\phi} \\ 
               C\,e^{\mp i\phi} & C \ea \right), \quad x\in[0,\,\pm L_{\rm sol}/2]
\label{E_01b}
\ee
\label{E_01}
\esube
and the zero matrix outside $[-L_{\rm sol}/2,\,L_{\rm sol}/2]$, where 
$C=A^2+B^2$ and $\phi=2\arctan(B/A)$.

Using the replacement \eqref{E_01b}, we can now calculate ${\bm \Phi}(L)$ in 
\eqref{e6_11} as
\be
{\bm \Phi}(L) \equiv {\bm \Phi}_{L_{\rm sol}} \, {\bm \Phi}_{\rm free}\,,
\qquad 
{\bm \Phi}_{L_{\rm sol}} = \left( 
     \ba{cc} \Phi_{11} & \Phi_{12} \\ -\Phi_{12} & \Phi_{11}^* \ea \right)\,,
\quad 
{\bm \Phi}_{\rm free}=e^{i\Omega (L-L_{\rm sol})\sthr}\,
\label{E_02}
\ee
where the last matrix is obtained
by solving the soliton-free part of \eqref{e6_10} on an interval 
$t\in(L_{\rm sol}/2,\,L-L_{\rm sol}/2)$.
The entries of ${\bm \Phi}_{L_{\rm sol}}$ depend on $A,\,B$, and $L_{\rm sol}$ in 
a complicated way, but their explicit form is not needed for our purpose;
we will only require the result, found by a tedious calculation, that
$\Phi_{12}\in i\mathbb{R}$. 
 The eigenvalues $\lambda$
of ${\bm \Phi}(L)$ in \eqref{E_02} are found from the quadratic equation:
\be
\lambda^2 - 2|\Phi_{11}|\cos\left[ \Omega (L-L_{\rm sol}) + \arg{\Phi_{11}}\right]
 \cdot \lambda + \left(|\Phi_{11}|^2 - |\Phi_{12}|^2\right) \,=\,0.
\label{E_03}
\ee

Since $\rho({\bm \Phi}(L)) = |\lambda|$, one can see that it varies with $L$
periodically, with the period being $2\pi/(2\Omega)$. 
Moreover, as stated in \eqref{e6_12}, 
\be
\|{\bm \Phi}(L)\| = \rho\left(\, {\bm \Phi}(L)\;{\bm \Phi}^{\dagger}(L) \,\right)^{1/2}
 \,=\, 
 \rho\left(\, {\bm \Phi}_{L_{\rm sol}}\;{\bm \Phi}^{\dagger}_{L_{\rm sol}} \,\right)^{1/2}
\,=\, |\Phi_{11}| + |\Phi_{12}|,
\label{E_04}
\ee
where $\dagger$ denotes Hermitian conjugation, and the last equation is found via a
somewhat tedious but straightforward calculation. This result means that even when the
``noise floor" does not grow {\em on average} over time, the amplitude of the ``noise
floor" harmonics can still increase by the factor $\|{\bm \Phi}(L)\|$ over a time
$t=L$; however, it will decrease in one or more of subsequent time intervals of the
same length (see Figs.~\ref{fig_11}(b,c)).




\begin{thebibliography}{99}


%
\bibitem{68_Strang}
G.~Strang, \
On the construction and comparison of difference schemes, \
SIAM J.~Numer.~Anal. 5 (1968) 506--517.
%
\bibitem{90_Suzuki}
M.~Suzuki, \ 
Fractal decomposition of exponential operators with applications
to many-body theories and Monte Carlo simulations, \ 
Phys.~Lett.~A 146 (1990) 319--323.
%
\bibitem{90_Yoshida}
H.~Yoshida,  \ 
Construction of higher order symplectic integrators, \ 
Phys.~Lett.~A 150 (1990) 262--268.
%
\bibitem{91_Glassner}
M.~Glassner, D.~Yevick, B.~Hermansson,  \
High-order generalized propagation techniques, \ 
J.~Opt.~Soc.~Am. B 8 (1991) 413--415.
%
%
\bibitem{93_Bandrauk}
A.D.~Bandrauk, H.~Shen,  \ 
Exponential split operator methods 
for solving coupled time-dependent Schr\"odinger equations, \ 
J.~Chem.~Phys. 99 (1993) 1185--1193.
%
\bibitem{73_HardinTappert}
R.H.~Hardin, F.D.~Tappert,  \ 
Applications of the split-step Fourier method to the numerical 
soltion of nonlinear and variable coefficient wave equations, \ 
SIAM Review (Chronicle) 15 (1973) 423.
%
\bibitem{84_AblowitzTaha}
T.~Taha, M.~Ablowitz,  \ 
Analytical and numerical aspects of certain nonlinear evolution equations.
II. Numerical, Nonlinear Schrodinger equation, \ 
J.~Comput.~Phys. 55 (1984) 203--230.
%
\bibitem{06_Bao}
W.~Bao, H.~Wang,  \
An efficient and spectrally accurate numerical method for computing 
dynamics of rotating Bose--Einstein condensates, \ 
J.~Comput.~Phys. 217 (2006) 612--626. 
%
\bibitem{76_VlasovP}
C.Z.~Cheng, G.~Knorr,  \ 
The integration of the Vlasov equation in configuration space, \ 
J.~Comput.~Phys. 22 (1976) 330--351. 
%
\bibitem{89_DiracSSM}
J.~de Frutos, J.M.~Sanz-Serna,  \ 
Split-step spectral schemes for nonlinear Dirac systems, \ 
J.~Comput.~Phys. 83 (1989) 407--423. 
%
\bibitem{03_ZakharovEqs}
W.~Bao, F.~Sun, G.W.~Wei, \ 
Numerical methods for the generalized Zakharov system, \ 
J.~Comput.~Phys. 190 (2003) 201--228. 
%
\bibitem{04_ZakharovEqs}
S.~Jin, P.A.~Markowich, C.~Zheng, \ 
Numerical simulation of a generalized Zakharov system, \ 
J.~Comput.~Phys. 201 (2004) 376--395. 
%
\bibitem{86_WH}
J.A.C.~Weideman,  B.M.~Herbst,  \ 
  Split-step methods for the solution of the 
    nonlinear Schr\"odinger equation, \
SIAM J.~Numer.~Anal. 23 (1986) 485--507.
%
\bibitem{10_GaucklerLubich}
L.~Gauckler, C.~Lubich,  \
Splitting integrators for nonlinear Schr\"odinger
equations over long times, \ 
Found.~Comput.~Math. 10 (2010) 275--302.
%
\bibitem{10_Faouetal}
E.~Faou, B.~Gr\'ebert, E.~Paturel,  \
Birkhoff normal form for splitting methods applied
to semilinear Hamiltonian PDEs. Part I.
Finite-dimensional discretization, \
Numer.~Math.  114 (2010) 429--458.
%
\bibitem{12_ja}
T.I.~Lakoba, \ 
Instability analysis of the split-step Fourier method on the background 
of a soliton of the nonlinear Schr\"odinger equation, \
 Numer.~Methods Partial Differ.~Equ. 28 (2012) 641--669.
%
\bibitem{13_JCP_methods}
J.~Xu, S.~Shao, H.~Tang, \ 
Numerical methods for nonlinear Dirac equation, \ 
J.~Comput.~Phys. 245 (2013) 131--149. 
%
\bibitem{16_SciChina}
W.~Bao, Y.Y.~Cai, X.W.~Xia, J.~Yin, \ 
Error estimates of numerical methods for the nonlinear Dirac equation
in the nonrelativistic limit regime, \ 
Sci.~China Math. 59 (2016) 1461--1494. 
%
\bibitem{17_NMPDEs}
S.-C.~Li, X.-G.~Li, F.-Y.~Shi, \ 
Time-splitting methods with charge conservation
for the nonlinear Dirac equation, \ 
 Numer.~Methods Partial Differ.~Equ. 33 (2017) 1582--1602. 
%
\bibitem{17_CMS}
M.~Lemou, F.~M\'ehats, X.~Zhao, \ 
Uniformly accurate numerical schemes for the nonlinear Dirac equation 
in the nonrelativistic limit regime, \ 
Commun.~Math.~Sci. 15 (2017) 1107--1128. 
%
\bibitem{18_ESAIM}
Y.Y.~Cai, Y.~Wang, \ 
A uniformly accurate (UA) multiscale time integrator pseudospectral method
for the nonlinear Dirac equation in the nonrelativistic limit regime, \ 
ESAIM: Math.~Model.~Numer.~Anal. 52 (2018) 543--566.
%
\bibitem{18_preprintSchratz}
P.~Kr\"amer, K.~Schratz, X.~Zhao, \ 
Splitting methods for nonlinear Dirac
equations with Thirring type interaction
in the nonrelativistic limit regime, \ 
preprint (https://publikationen.bibliothek.kit.edu/1000087636). 
%
\bibitem{98_JOSAB_SSMforGS}
G.~Lenz, B.J.~Eggleton, N.M.~Litchinitser, \ 
Pulse compression using fiber gratings as highly
dispersive nonlinear elements, \ 
J.~Opt.~Soc.~Am.~B 15 (1998) 715--721. 
%
\bibitem{04_Mak_SSMforGS}
W.C.K.~Mak , B.A.~Malomed, P.L. Chu, \ 
Slowdown and splitting of gap solitons in apodized Bragg gratings, \ 
J.~Mod.~Opt. 51 (2004) 2141--2158. 
%
\bibitem{17_PRE_SSMforGS}
T.~Ahmed, J.~Atai, \ 
Bragg solitons in systems with separated nonuniform Bragg grating 
and nonlinearity, \ 
Phys.~Rev.~E 96 (2017) 032222. 
%
\bibitem{17_SciRep_SSMforGS}
S.A.M.S.~Chowdhury, J.~Atai, \ 
Moving Bragg grating solitons in a semilinear dual-core system with
dispersive reflectivity, \ 
Sci.~Rep. 7 (2017) 4021. 
%
\bibitem{04_DiracMaxwell}
W.~Bao, X.-G.~Li, \ 
An efficient and stable numerical method for the Maxwell--Dirac system, \ 
J.~Comput.~Phys. 199 (2004) 663--687. 
%
\bibitem{17_DiracMaxwell}
W.~Bao, Y.Y.~Cai, X.W.~Xia, Q.~Tang, \ 
Numerical methods and comparison for the Dirac equation in the 
nonrelativistic limit regime, \ 
J.~Sci.~Comput. 71 (2017) 1094--1134. 
%
\bibitem{14_NLDE_numerics}
S.~Shao, N.R.~Quintero, F.G.~Mertens, F.~Cooper, A.~Khare, A.~Saxena, \
Stability of solitary waves in the nonlinear Dirac equation with arbitrary
nonlinearity, \ 
Phys.~Rev.~E 90 (2014) 032915.
%
\bibitem{15_NLDE_numerics}
J.~Cuevas-Maraver, P.G.~Kevrekidis, A.~Saxena, F.~Cooper, F.G.~Mertens, \
Solitary waves in the nonlinear Dirac equation at the continuum limit: 
Stability and dynamics, \
In: Ord.~Part.~Diff.~Eqs., (Nova Science, Boca Raton), Chap.~4.
%
\bibitem{12_BerkCom}
G.~Berkolaiko, A.~Comech, \ 
On spectral stability of solitary waves of nonlinear Dirac equation in 1D, \
Math.~Model.~Nat.~Phenom. 7 (2012) 13--31. 
%
\bibitem{17_jaNLDE}
T.I.~Lakoba, \
Numerical study of solitary wave stability in cubic nonlinear 
Dirac equations in 1D, \ 
Phys.~Lett.~A 382 (2018) 300--308.
%
\bibitem{75_Exact}
S.Y.~Lee, T.K.~Kuo, A.~Gavrielides, \
Exact localized solutions of two-dimensional field theories of massive fermions
with Fermi interactions, \ 
Phys.~Rev.~D  12 (1975) 2249--2253. 
%
\bibitem{81_2solColl}
A.~Alvarez, B.~Carreras, \ 
Interaction dynamics for the solitary waves of a nonlinear Dirac model, \ 
Phys.~Lett.~A 86 (1981) 327--332. 
%
\bibitem{Stewart_Matrix}
G.W.~Stewart, \
Introduction to Matrix Computations, \
(Academic, New York, 1973); Secs.~6.2 and 6.6. 
%
\bibitem{75_MTMsol_1}
S.J.~Orfanidis, R.~Wang, \ 
Soliton solutions of the Massive Thirring model, \ 
Phys.~Lett.~B 57 (1975) 281--283. 
%
\bibitem{75_MTMsol_2}
S.-J.~Chang, S.D.~Ellis,  B.W.~Lee, \
Chiral confinement: An exact solution of the massive Thirring model, \ 
Phys.~Rev.~D 12 (1975)  3572--3582. 
%
\bibitem{77_KaupNewell}
D.J.~Kaup, A.C.~Newell, \ 
On the Coleman correspondence and the solution of the Massive Thirring model, \
Lett.~Nuovo Cimento 20 (1977) 325--331.
%
\bibitem{77_KuzMikh}
E.A.~Kuznetsov, A.V.~Mikhailov, \ 
On the complete integrability of the two-dimensional classical Thirring model, \
Theor.~Math.~Phys. 30 (1977) 193--200.
%
\bibitem{96_KaupL}
D.J.~Kaup, T.I.~Lakoba, \ 
The squared eigenfunctions of the massive Thirring model
in laboratory coordinates, \ 
J.~Math.~Phys. 37 (1996) 308--323. 
%
\bibitem{85_Winful}
H.G.~Winful, \ 
Pulse compression in optical fiber filters, \ 
Appl.~Phys.~Lett. 46 (1985) 527--529. 
%
\bibitem{98_Peli}
I.V.~Barashenkov, D.E.~Pelinovsky, E.V.~Zemlyanaya, \ 
Vibrations and Oscillatory Instabilities of Gap Solitons, \ 
Phys. Rev. Lett. 80 (1998) 5117--5120.
%
\bibitem{02_CoxMatt}
S.M.~Cox, P.C.~Matthews, \ 
Exponential time differencing for stiff systems, \ 
J.~Comput.~Phys. 176 (2002) 430--455.
%
\bibitem{10_IFNLDE}
F.~de la Hoz, F.~Vadillo, \ 
An integrating factor for nonlinear Dirac equations, \
Comput.~Phys.~Commun. 181 (2010) 1195--1203.
%
\bibitem{19_DiracPML}
X.~Antoine, E.~Lorin, \ 
A simple pseudospectral method for the computation of the time-dependent 
Dirac equation with Perfectly Matched Layers, \ 
J.~Comput.~Phys. 395 (2019) 583--601. 
%
\bibitem{17_coupledDirac}
A.~Khare, F.~Cooper, A.~Saxena, \ 
Approximate analytic solutions to coupled nonlinear Dirac equations, \ 
Phys.~Lett.~A 381 (2017) 1081--1086. 
%
\bibitem{19_IgorPTDirac}
N.V.~Alexeeva, I.V.~Barashenkov, A.~Saxena, \ 
Spinor solitons and their PT-symmetric offspring, \ 
Ann.~Phys. 403 (2019) 198--223.
%
\bibitem{72_Zakharov}
V.E.~Zakharov, \ 
Collapse of Langmuir waves, \ 
Sov.~Phys. JETP 35 (1972) 908--914. 
%
\bibitem{16_ja2}
T.I.~Lakoba, \ 
Instability of the finite-difference split-step method applied to 
the nonlinear Schr\"odinger equation. II. moving soliton, \
Numer.~Methods Partial Differ.~Equ. 32 (2016) 1024--1040. 
%


%
%




\end{thebibliography}
\end{document}